\documentclass[10pt]{amsart}

\usepackage{xypic,amscd,amssymb,latexsym}
\usepackage{amsmath}	
\usepackage[breaklinks=true]{hyperref}

\usepackage{mathrsfs}

\usepackage{url}

\usepackage{pdfsync}

\usepackage{cancel}

\begin{document}
\pagenumbering{arabic}

\newtheorem{theorem}{Theorem}[section]
\newtheorem{proposition}[theorem]{Proposition}
\newtheorem{lemma}[theorem]{Lemma}
\newtheorem{corollary}[theorem]{Corollary}
\newtheorem{remark}[theorem]{Remark}
\newtheorem{definition}[theorem]{Definition}
\newtheorem{question}[theorem]{Question}
\newtheorem{claim}[theorem]{Claim}
\newtheorem{conjecture}[theorem]{Conjecture}
\newtheorem{defprop}[theorem]{Definition and Proposition}
\newtheorem{example}[theorem]{Example}
\newtheorem{deflem}[theorem]{Definition and Lemma}

\def\qed{{\quad \vrule height 8pt width 8pt depth 0pt}}

\newcommand{\cplx}[0]{\mathbb{C}}

\newcommand{\vs}[0]{\vspace{2mm}}

\newcommand{\til}[1]{\widetilde{#1}}

\newcommand{\mcal}[1]{\mathcal{#1}}

\newcommand{\ul}[1]{\underline{#1}}

\newcommand{\ol}[1]{\overline{#1}}

\newcommand{\wh}[1]{\widehat{#1}}

\newcommand{\mut}[1]{\stackrel{#1}{\to}}

\author{Hyun Kyu Kim}
\email{hyunkyukim@ewha.ac.kr, hyunkyu87@gmail.com}

\address{Department of Mathematics, Ewha Womans University, 52 Ewhayeodae-gil, Seodaemun-gu, Seoul 03760, Republic of Korea}

\numberwithin{equation}{section}

\title[Phase constants in quantum cluster varieties]{Phase constants in the Fock-Goncharov quantum cluster varieties}

\begin{abstract}
A cluster variety of Fock and Goncharov is a scheme constructed by gluing split algebraic tori, called seed tori, via birational gluing maps called mutations. In quantum theory, the ring of functions on seed tori are deformed to non-commutative rings, represented as operators on Hilbert spaces. Mutations are quantized to unitary maps between the Hilbert spaces intertwining the representations. These unitary intertwiners are described using the quantum dilogarithm function $\Phi^\hbar$. Algebraic relations among classical mutations are satisfied by the intertwiners up to complex constants. The present paper shows that these constants are $1$. So the mapping class group representations resulting from the Chekhov-Fock-Goncharov quantum Teichm\"uller theory are genuine, not projective. During the course, the hexagon and the octagon operator identities for $\Phi^\hbar$ are derived.
\end{abstract}

\maketitle

\tableofcontents

\newpage

\section{Introduction}
\label{sec:introduction}

A cluster variety of Fock and Goncharov \cite{FG09} \cite{FG09b} can be viewed as an algebro-geometric space whose ring of regular functions is an analog of the upper cluster algebra of Fomin and Zelevinsky \cite{FoZ02}. More precisely, it is a scheme over $\mathbb{Z}$ constructed by gluing split algebraic tori $(\mathbb{G}_m)^n$ along certain birational maps, where the tori are enumerated by the data called seeds in the theory of cluster algebras, while the birational maps correspond to the mutations of seeds. There are three different types of cluster varieties, namely $\mcal{A}$, $\mcal{X}$, and $\mcal{D}$, according to how these birational maps are given. In fact, only the cluster $\mcal{A}$-variety provides an honest example of a cluster algebra in the sense of \cite{FoZ02}, while the $\mcal{X}$- and the $\mcal{D}$-varieties are certain generalizations. 

\vs

Fix a rank $n \in \mathbb{Z}_{>0}$. Underlying a \emph{seed} is a combinatorial data, which essentially can be encoded as an $n\times n$ integer matrix $\varepsilon$, called the \emph{exchange matrix}, which is required to be skew-symmetrizable, i.e. must be skew-symmetric when multiplied by some diagonal matrix from the left. In particular, when $\varepsilon$ is skew-symmetric, it can be realized as the signed adjacency matrix of a quiver, i.e. a graph with oriented edges, without cycles of length $1$ or $2$. A seed is also equipped with $n$ algebraically independent formal variables, which can be thought of as being attached to the vertices of a quiver. If one picks any $k \in \{1,\ldots,n\}$, or a vertex of a quiver, one can apply the \emph{mutation} $\mu_k$ to the seed to obtain another seed; the quiver or the exchange matrix changes according to a certain combinatorial rule, while the new $n$ variables are related to the former $n$ variables by certain rational formulas which depend on the exchange matrix and are `positive', i.e. can be written as quotients of polynomials with positive coefficients. The mutation formulas for the exchange matrices and the attached variables might seem ad hoc at a first glance. So it is remarkable that such a structure appears in very many areas of mathematics.

\vs

What is even more remarkable is a Poisson-like structure on a cluster variety, which is written in a simple form on each torus in terms of the exchange matrix and is compatible under the mutation formulas. More precisely, the $\mcal{A}$-, the $\mcal{X}$-, and the $\mcal{D}$-varieties are equipped with a $2$-form (of `$K_2$ class'), a Poisson structure, and a symplectic form, respectively. So it is sensible to investigate whether there exist deformation quantizations of the $\mcal{X}$- and the $\mcal{D}$-varieties. The first step would be to construct a quantum version of the cluster variety, which I explain below. 

\vs

For a cluster variety we associate to each seed a torus $(\mathbb{G}_m)^n$, which can be viewed as being defined as the space whose ring of regular functions is the ring of Laurent polynomials in $n$ variables. So, one can think of this situation as attaching such a commutative ring to each seed. To each mutation we associate a birational map between the two tori, that is, an isomorphism of the fields of fractions of the attached Laurent polynomial rings.

\vs

Establishing a quantum cluster variety means that we would like to associate to each seed a one-parameter family of non-commutative rings which deforms the commutative ring of Laurent polynomials, while to each mutation an isomorphism of the skew fields of fractions of the non-commutative rings associated to the relevant seeds. First, this `quantum' isomorphism of skew fields must be a deformation of the corresponding `classical' isomorphism of fields of fractions of the commutative Laurent polynomial rings. Second, more importantly, there must be a consistency. Namely, there exists sequences of mutations that return the same exchange matrices after their application and whose induced classical isomorphisms of the fields of fractions are the identity maps. It is natural to require that the quantum isomorphisms of the skew fields of fractions induced by these sequences are also identity maps. However, classification of all such sequences for classical setting is still an open problem; see e.g. \cite{KY}. At this point, one can only say that there are some known examples, most prominent ones being the following:
\begin{enumerate}
\item[\rm (S1)] the rank $1$ identity (or, the `twice-flip' identity, or $A_1$-type identity) :  $\mu_k \mu_k = \mathrm{Id}$;

\item[\rm (S2)]  the rank $2$ identities:, 

\begin{enumerate}
\item[\rm (S2-1)] ($A_1\times A_1$-type identity, or commuting identity) $\mu_j \mu_i \mu_j \mu_i={\rm Id}$ in case $\varepsilon_{ij}=\varepsilon_{ji}=0$;

\item[\rm (S2-2)] ($A_2$-type identity, or the {\em pentagon identity}) $\mu_i \mu_j \mu_i \mu_j \mu_i = (i \, j)$ in case $\varepsilon_{ij}=-\varepsilon_{ji} \in \{1,-1\}$, where $(i \, j)$ stands for the relabling $i\leftrightarrow j$.

\item[\rm (S2-3)] ($B_2$-type identity, or the {\em hexagon identity}) $\mu_j \mu_i \mu_j \mu_i \mu_j \mu_i = {\rm Id}$ in case $\varepsilon_{ij}=-2\varepsilon_{ji}\in \{2,-2\}$ or $\varepsilon_{ji}=-2\varepsilon_{ij}\in \{2,-2\}$;

\item[\rm (S2-4)] ($G_2$-type identity, or the {\em octagon identity}) $\mu_j \mu_i \mu_j \mu_i \mu_j \mu_i \mu_j \mu_i = {\rm Id}$ in case $\varepsilon_{ij}=-3\varepsilon_{ji} \in \{3,-3\}$ or $\varepsilon_{ji}=-3\varepsilon_{ij} \in \{3,-3\}$;
\end{enumerate}
\end{enumerate}
Note that the famous pentagon identity arises in many different areas of mathematics, such as for the flips of ideal triangulations of punctured surfaces, or for the consistency axiom for the associativity morphisms in a tensor category; see e.g. \cite{FK12}. Here the identities (S2-3) and (S2-4) are referred to by the names hexagon and octagon, as analogs of the pentagon identity (S2-2). Fock and Goncharov \cite{FG09} constructed quantum isomorphisms of skew fields associated to mutations and checked that they satisfy this consistency for the above identities (S1) and (S2). 

\vs

As is often the case for a `physical' quantization, one would like to represent the non-commutative ring associated to a seed as an algebra of operators, i.e. to find a family of representations of this ring on a Hilbert space parametrized by a real quantum parameter $\hbar$ commonly referred to as the Planck constant, satisfying certain desirable analytic properties. Then, to each mutation $\mu_k$ one would like to associate a unitary map ${\bf K}^\hbar_k$ between the Hilbert spaces assigned to the relevant seeds that  intertwines the representations of the non-commutative rings which are related by the quantum isomorphism. In the Fock-Goncharov quantization \cite{FG09}, one advantage of having this intertwining map is that its formula is more transparent than the quantum isomorphism formulas for the non-commutative algebras. In particular, while the formula for the quantum isomorphism of algebras looks somewhat complicated and not much enlightening, the corresponding intertwiner operator is written quite neatly via functional calculus applied to a famous special function called the {\em non-compact quantum dilogarithm} $\Phi^\hbar(z)$ of Faddeev and Kashaev \cite{F95} \cite{FK94}, defined for a real parameter $\hbar$ and a complex number $z$ living in the strip $|{\rm Im}(z)| < \pi(1+\hbar)$ as
$$
\Phi^\hbar(z) = \exp\left( - \frac{1}{4} \int_\Omega \frac{e^{-{\rm i}pz}}{\sinh(\pi p) \sinh(\pi \hbar p)} \frac{dp}{p} \right),
$$
where ${\rm i} = \sqrt{-1}$, and $\Omega$ is the contour along the real line that avoids the origin by a small half-circle above the origin. In fact, this integral expression appeared 100 yeas ago \cite{B01}, and has gained much interest recently.

\vs

The unitary map ${\bf K}^\hbar_k$ is required not only to just intertwine the representations of the quantum algebras, but also to satisfy the consistency relations corresponding to the special mutation sequences (S1)--(S2) mentioned above. These consistency relations are operator identities, and they look like (S1) ${\bf K}^\hbar_k {\bf K}^\hbar_k = {\rm Id}$, (S2-1) ${\bf K}^\hbar_i {\bf K}^\hbar_j {\bf K}^\hbar_i {\bf K}^\hbar_j = {\rm Id}$ when $\varepsilon_{ij}=0$, (S2-2) ${\bf K}^\hbar_i {\bf K}^\hbar_j {\bf K}^\hbar_i {\bf K}^\hbar_j {\bf K}^\hbar_i = (i j)$ when $\varepsilon_{ij} = - \varepsilon_{ji} \in \{1,-1\}$, etc. Such consistency condition, which in fact is required to hold only up to multiplicative constants, is what makes the quantization problem both more difficult and more interesting. The approach taken by Fock and Goncharov in \cite{FG09} is to deduce this consistency of intertwiners essentially by taking advantage of the `irreducibility' of the representations of rings which are being intertwined, instead of directly checking the operator identities. The situation is roughly like Schur's lemma, which says that a complex vector space endomorphism that commutes with all the operators for an irreducible representation of a group or a ring on that vector space is a scalar operator. In particular, for each sequence of mutations in (S1)--(S2), Fock and Goncharov showed that the composition of corresponding unitary intertwiners is a scalar operator on a relevant Hilbert space \cite[Thm.5.5]{FG09}. One may interpret this result as having a unitary projective representation  on a Hilbert space of the groupoid formed by mutations, where the word `projective' carries the connotation that the relations are satisfied only up to complex scalars. As all operators are unitary, these scalars are of modulus $1$, hence can be viewed as \emph{phases}. The existence of these phase scalars, which can be thought of as a certain anomaly, is not at all a defect of the program of quantization of cluster varieties. As a matter of fact, this projective anomaly alone may already contain a very interesting bit of information. Ones of the prominent examples of cluster varieties are various versions of the \emph{Teichm\"uller spaces} of non-compact Riemann surfaces, where the \emph{mapping class groups} can be embedded into the groupoid of mutations. The theory of quantum cluster varieties is thus an approach to the theory of quantum Teichm\"uller spaces, one of the main results of which is a family of unitary projective representations of mapping class groups on Hilbert spaces. The anomaly can then be interpreted as group $2$-cocycles of mapping class groups, or equivalently, as central extensions of these groups, and this led to the works of Funar-Sergiescu \cite{FS10} and Xu \cite{X14}; see also \cite{K12}.

\vs

However, these phase scalars for the projective anomaly have not been determined as precise complex numbers. The main result of the present paper, Thm.\ref{thm:main}, is that these scalars are all $1$. So, we get \emph{genuine} representations of groupoid of mutations on Hilbert spaces, instead of projective representations. 
To prove this I perform explicit computations involving unitary operators, instead of resorting to the irreducibility of representations and the Schur's-lemma-type philosophy.  While doing so, in order to exploit results from functional analysis, I replace Fock-Goncharov's  `Heisenberg relations' by the corresponding `Weyl-type relations' which are more rigid, and also introduce and study the concept of \emph{special affine shift operators} on $L^2(\mathbb{R}^n)$ (\S\ref{subsec:special_affine_shift_operators}), which are the operators naturally induced by the `special affine' transformations $\mathrm{SL}_\pm(n,\mathbb{R}) \ltimes \mathbb{R}^n$ of $\mathbb{R}^n$, in order to facililate notations and computations for unitary operators, which may be used in future projects too, such as the upcoming work \cite{KS}. 

\vs

A crucial technical point is about the decomposition ${\bf K}^\hbar_k = {\bf K}^{\hbar\sharp}_k \circ {\bf K}'_k$ of the unitary intertwiner for a mutation into the `automorphism part' ${\bf K}^{\hbar\sharp}_k$ and the `linear part' ${\bf K}'_k$, which was already used by Fock and Goncharov \cite{FG09}. The linear part ${\bf K}'_k$ is essentially just a special affine shift operator, while the automorphism part ${\bf K}^{\hbar\sharp}_k = \Phi^\hbar({\bf x}^\hbar_k) \, \Phi^\hbar(\til{\bf x}^\hbar_k)^{-1}$ is expressed via functional calculus for certain basic mutually commuting self-adjoint quantum coordinate operators ${\bf x}^\hbar_k$ and $\til{\bf x}^\hbar_k$ applied to the non-compact quantum dilogarithm function $\Phi^\hbar$. Given a composition of bunch of ${\bf K}^\hbar_k$, write each factor as this decomposition, and move around all the linear parts to group them together. The composition of these linear parts is not canceled in general, but can in principle be explicitly computed as a special affine shift operator. The remaining are the quantum dilogarithm factors, which one should prove to equal a desired special affine shift operator, by using known operator identities for quantum dilogarithm, such as the famous quantum pentagon identity for $\Phi^\hbar$, which is needed for the identity (S2-2). For (S2-3) and (S2-4), we would need analogous operator identities, namely the hexagon and the octagon identities for $\Phi^\hbar$, which are hinted in the literature but have not been rigorously proved. 

\vs

To overcome this situation, I first adapt the idea of a {\em signed} decomposition of the mutation intertwiner ${\bf K}^\hbar_k$; namely, I consider certain two different ways of decomposing into the automorphism and the linear parts, enumerated by a sign$\in\{+,-\}$. Such idea has been developed and used in the cluster algebra literature, e.g. in \cite{KN} \cite{Keller}, and has deeper categorical meaning; in the present paper I enhance and justify this idea to the operator level, for representation of quantum cluster $\mathcal{D}$-variety. What is also known in the algebraic level is that, for a mutation sequence, there is a preferred choice of signs for these mutations, called {\em tropical signs}, which are related to the so-called `sign-coherence conjecture for ${\bf c}$-vectors', proved in general in \cite{GHKK}. If we consider the signed decompositions of the mutation intertwiners for these tropical signs, one can show that the linear parts all cancel one another when grouped together.  

\vs

After removing all the linear parts this way from the composition of mutation intertwiners for a mutation sequence, what is left is composition of the automorphism parts, i.e. quantum dilogarithm factors only; the task is to show that this composition precisely equals the identity operator, without a multiplicative constant. Instead of directly proving this, I invoke Fock-Goncharov's statement in \cite[Thm5.5]{FG09} which implies that this composition is a scalar operator. Note that each automorphism part is composition of two quantum dilogarithm factors which commute. By some basic functional analysis, I `separate' these two commuting factors, for each automorphism part, to obtain two separated operator identities for quantum dilogarithm, each of which hold up to a constant. This way I prove the hexagon and the octagon identities for $\Phi^\hbar$ (Prop.\ref{prop:hexagon_identity_for_ncQD} and Prop.\ref{prop:octagon_identity_for_ncQD}), which I then apply back to the sought-for composition of quantum dilogarithm factors to compute the scalar precisely, as desired. I note that the previously known arguments for these hexagon and octagon identities for the non-compact quantum dilogarithm $\Phi^\hbar$ were only heuristic, based on analogous results for the {\em compact quantum dilogarithm}. A key point of the proof of these identities in the present paper is to extract from their `(symplectic) double' version identities proved in \cite[Thm.5.5]{FG09} the sought-for `single' version identities. 

\vs

One consequence of the triviality of the phase scalars proved in the present paper is that the mapping class group representations coming from the version of  quantum Teichm\"uller theory that is based on the Fock-Goncharov quantization of cluster varieties are genuine, not projective. Since the works \cite{FS10} \cite{X14} are based on the assumption that these representations are projective but not genuine, the triviality result of the present paper unfortunately deprives of the representation theoretic meaning from 
these two works at the moment, leaving only the group theoretic meaning. Meanwhile, the projective anomaly of such representations was further investigated in my previous work \cite{K12} in comparison with another version of quantum Teichm\"uller spaces obtained by Kashaev, which yields projective but not genuine representations. It was observed in \cite{K12} that these two versions of quantum Teichm\"uller spaces are different in a very interesting sense, related to surface braid groups. In my opinion, the result of \cite{K12} suggests that there exists a version of quantum Teichm\"uller spaces in the style of Chekhov-Fock-Goncharov that 1) yields cohomologically non-trivial phase constants, such that 2) these constants are (half-)integer powers of $\zeta = e^{-\frac{\pi{\rm i}}{6}} c_\hbar$, where $c_\hbar=e^{-\frac{\pi{\rm i}}{12}(\hbar+\hbar^{-1})}$ (eq.\eqref{eq:c_hbar}). However, such has not been explicitly constructed yet, hence calls for a future investigation. Another possible topic of further research is about whether the arguments of the present paper can be generalized to {\em any} sequence of mutations that induces the identity maps on the exchange matrices and the fields of fractions of the classical coordinate rings. For example, the shortest known such sequence which is not covered by the rank 2 identities is a certain sequence of 32 mutations, applied to a rank 6 seed, as investigated in \cite{KY}. One might try to check whether the results of the present paper extend to such sequences.

\vs

\noindent{\bf Acknowledgments.} I would like to thank Myungho Kim, Dylan Allegretti, Ivan Chi-ho Ip, Carlos Scarinci, Seung-Jo Jung, Woocheol Choi, Louis Funar, Vladimir V. Fock, and Alexander B. Goncharov for helpful discussions. This work was supported by the Ewha Womans University Research Grant of 2017. This research was supported by Basic Science Research Program through the National Research Foundation of Korea(NRF) funded by the Ministry of Education(grant number 2017R1D1A1B03030230).

\section{Cluster $\mcal{A}$-, $\mcal{X}$- and $\mcal{D}$-varieties}

I shall describe Fock-Goncharov's cluster $\mcal{A}$-varieties, cluster $\mcal{X}$-varieties, and cluster $\mcal{D}$-varieties. Definitions and treatment in this section are from \cite{FG09} and references therein, with certain modifications.

\subsection{Seed and seed tori}

In the present paper, we choose and fix one positive integer $n$, which can be regarded as the \ul{\em rank} of the relevant cluster algebras/varieties.
\begin{definition}
\label{def:skew_symmetrizable_matrix}
Let ${\bf k}$ be a field. We say that an $n\times n$ matrix $\varepsilon = (\varepsilon_{ij})_{i,j=1,\ldots,n}$ with entries in ${\bf k}$ is \ul{\em skew-symmetrizable} if there exist $d_1,\ldots,d_n \in {\bf k}^*$ such that the matrix $(\wh{\varepsilon}_{ij})_{i,j}$ defined as
\begin{align}
\label{eq:wh_varepsilon}
\wh{\varepsilon}_{ij} := \varepsilon_{ij} \, d_j^{-1}
\end{align}
is skew-symmetric, i.e. $\wh{\varepsilon}_{ij} = - \wh{\varepsilon}_{ji}$, $\forall i,j$. We call such a collection $d = (d_i)_{i=1,\ldots,n}$ a \ul{\em skew-symmetrizer} of the matrix $\varepsilon=(\varepsilon_{ij})_{i,j}$.
\end{definition}
For any skew-symmetrizable matrix over ${\bf k}=\mathbb{Q}$, one can always find a skew-symmetrizer consisting only of integers.

\begin{definition}[seed]
\label{def:seed}
An \ul{\em $\mcal{A}$-seed} $\Gamma$ is a triple $(\varepsilon, d, \{A_i\}_{i=1}^n)$ of an $n\times n$ skew-symmetrizable integer matrix $\varepsilon$, a skew-symmetrizer $d$ consisting of {\em positive} integers, and an ordered set of formal variables $A_1,A_2,\ldots,A_n$ called the \ul{\em cluster $\mcal{A}$-variables}.

\vs

An \ul{\em $\mcal{X}$-seed} $\Gamma$ is a triple $(\varepsilon, d, \{X_i\}_{i=1}^n)$ of $\varepsilon$, $d$ as above, and an ordered set of variables $X_1,\ldots,X_n$ called the \ul{\em cluster $\mcal{X}$-variables}.
\vs

An \ul{\em $\mcal{D}$-seed} $\Gamma$ is a triple $(\varepsilon, d,\{B_i,X_i\}_{i=1}^n)$ of $\varepsilon,d$ as above, and an ordered set of variables $B_1,\ldots,B_n,X_1,\ldots,X_n$ called the \ul{\em cluster $\mcal{D}$-variables}.

\vs

We call any of these a \ul{\em seed}, $\varepsilon$ the \ul{\em exchange matrix} of the seed, the relevant formal variables the \ul{\em cluster variables} of the seed, and the relevant symbols $\mcal{A}$, $\mcal{X}$, $\mcal{D}$ the \ul{\em kind} of the seed.  We write a seed by $\Gamma = (\varepsilon, d, *)$ if it is clear from the context what the cluster variables are. 
\end{definition}

Some words must be put in order. What Fock and Goncharov \cite{FG09} call a `feed' is essentially a seed without the cluster variables, i.e. the exchange matrix $\varepsilon$ and the skew-symmetrizer $d$. I think it makes more sense to use seeds as defined in Def.\ref{def:seed} instead of feeds, in the construction of cluster varieties, whence I do so in the present paper. An $\mcal{A}$-seed, with the data of $d$ left out, is what is called a `labeled seed' with `trivial coefficients' in the theory of cluster algebras initiated by Fomin and Zelevinsky \cite{FoZ02}.  A $\mcal{D}$-seed minus the data of $d$ is what is called a `labeled seed' with `coefficients' in the `universal semifield' in \cite{FoZ07}; $\{B_i\}_{i=1}^n$ is their `cluster' and $X_i$'s their coefficients. Readers will be able to verify these comparison when they get to the  formulas for the {\em mutations} in \S\ref{subsec:seed_mutations}. On the other hand, the cluster $\mcal{X}$-variables are not the usual `cluster variables' in the sense of Fomin-Zelevinsky. 
So, a {\em seed} as defined in the above Def.\ref{def:seed} is, say,  a `generalized seed together with the choice of a skew-symmetrizer'. The symbol $\Gamma$ has nothing to do with the cluster mapping class group (cluster modular group)  in \cite{FG09}.

\vs

For a $\mcal{D}$-seed there is an interesting set of redundant variables: 
\begin{definition}[tilde variables for $\mcal{D}$-seed]
\label{def:tilde_variables}
For a $\mcal{D}$-seed $\Gamma = (\varepsilon, d, \{B_i,X_i\}_{i=1}^n)$, define new variables $\til{X}_i$, $i=1,\ldots,n$, by
$$
\til{X}_i := X_i \, \prod_{j = 1}^n B_j^{\varepsilon_{ij}}, \qquad \forall i = 1,\ldots,n.
$$
\end{definition}

For each seed, we consider an affine scheme whose coordinate functions are identified with the cluster variables.
\begin{definition}[seed tori]
\label{def:seed_tori}
To the seeds in Def.\ref{def:seed}, all denoted by $\Gamma$ by abuse of notations, we assign respectively the split algebraic tori
\begin{align*}
\mcal{A}_\Gamma := (\mathbb{G}_m)^n, \qquad
\mcal{X}_\Gamma = (\mathbb{G}_m)^n, \qquad
\mcal{D}_\Gamma = (\mathbb{G}_m)^{2n},
\end{align*}
called the \ul{\em seed $\mcal{A}$-torus}, the \ul{\em seed $\mcal{X}$-torus}, and the \ul{\em seed $\mcal{D}$-torus}, where $\mathbb{G}_m$ is the multiplicative algebraic group, satisfying $\mathbb{G}_m ({\bf k}) = {\bf k}^*$ for any field ${\bf k}$. We call these \ul{\em seed tori}, collectively. We identify the canonical coordinate functions of $\mcal{A}_\Gamma$ with the seed's cluster variables $A_i$'s, those of $\mcal{X}_\Gamma$ with $X_i$'s, and those of $\mcal{D}_\Gamma$ with $B_i$'s and $X_i$'s.
\end{definition}

So, one $\mcal{A}$-seed $\Gamma$ gives us one affine scheme $\mcal{A}_\Gamma$. We shall consider $\mcal{A}_\Gamma$ for different $\Gamma$'s, and glue them together by certain birational maps, to construct a {\em cluster $\mcal{A}$-variety}. We do likewise for the spaces $\mcal{X}_\Gamma$ for different $\Gamma$'s to construct a {\em cluster $\mcal{X}$-variety}, and similarly for $\mcal{D}_\Gamma$'s to construct a {\em cluster $\mcal{D}$-variety}. I will explain these gluings in the subsequent subsections. Note that we never glue seed tori of different kinds, whence the abuse of notations committed by labeling the different kinds of seeds by the same letter $\Gamma$ will not be too harmful. The ring of regular functions on each seed torus is the ring of all Laurent polynomials over $\mathbb{Z}$ in $n$ (or $2n$) coordinate variables.

\vs

What make the cluster varieties a lot richer and more interesting are the remarkable geometric structures, which are defined as follows on their local patches, i.e. on the seed tori.
\begin{definition}[geometric structures on seed tori]
\label{def:geometric_structures_on_seed_tori}
We equip the following $2$-form on $\mcal{A}_\Gamma$ 
\begin{align}
\label{eq:2-form_on_A_Gamma}
\Omega_\Gamma = \sum_{i,j \in \{1,\ldots,n\}} \, \til{\varepsilon}_{ij} \, d \log A_i \wedge d\log A_j, \quad \mbox{where} \quad \til{\varepsilon}_{ij} := d_i \varepsilon_{ij},
\end{align}
and the following Poisson structure on $\mcal{X}_\Gamma$
$$
\{X_i, X_j\} = \wh{\varepsilon}_{ij} \, X_i X_j, \qquad \forall i,j \in \{1,\ldots,n\}.
$$
On $\mcal{D}_\Gamma$ we consider the Poisson structure defined by
\begin{align}
\label{eq:Poisson_on_D}
\{B_i, B_j\} = 0, \qquad
\{X_i, B_j\} = d_i^{-1} \delta_{ij} X_i B_j, \qquad
\{X_i, X_j\} = \wh{\varepsilon}_{ij} X_i X_j,
\end{align}
where $\delta_{ij}$ is the Kronecker delta, as well as the following $2$-form
\begin{align}
\label{eq:2-form_on_D}
- \frac{1}{2} \sum_{i,j\in \{1,\ldots,n\}} \til{\varepsilon}_{ij} \, d\log B_i \wedge d \log B_j - \sum_{i \in \{1,\ldots,n\}} d_i \, d\log B_i \wedge d\log X_i.
\end{align}
\end{definition}

\begin{lemma}[\cite{FG09}]
The $2$-form \eqref{eq:2-form_on_D} on $\mcal{D}_\Gamma$ is a symplectic form, and is compatible with the Poisson structure \eqref{eq:Poisson_on_D} on it.
\end{lemma}

The tilde variables defined in Def.\ref{def:tilde_variables} satisfy
\begin{align}
\label{eq:Poisson_on_D_tilde}
\{ \til{X}_i, B_j\} = d_i^{-1} \delta_{ij} \til{X}_i B_j, \qquad
\{\til{X}_i, \til{X}_j\} = -\wh{\varepsilon}_{ij} \til{X}_i \til{X}_j, \qquad
\{\til{X}_i, X_j\}=0.
\end{align}

\subsection{Seed mutations and seed automorphisms}
\label{subsec:seed_mutations}

In order to study the gluing of the seed tori, we first need to investigate the following transformation rules for seeds, which we can also understand as a way of recursively creating new seeds from previously constructed ones.
\begin{definition}[seed mutation]
\label{def:mutation}
For $k\in\{1,2,\ldots,n\}$, a seed $\Gamma' = (\varepsilon', d', *')$  is said to be obtained by applying the \ul{\em seed mutation $\mu_k$ in the direction $k$} to a seed $\Gamma = (\varepsilon, d, *)$, if all of the following hold: 

\begin{itemize}
\item the exchange matrices are related by
\begin{align}
\label{eq:varepsilon_prime_formula}
\varepsilon'_{i j} = \left\{
\begin{array}{ll}
-\varepsilon_{ij} & \mbox{if $i=k$ or $j=k$,} \\
\varepsilon_{ij} + \frac{1}{2}( |\varepsilon_{ik}| \varepsilon_{kj} + \varepsilon_{ik} | \varepsilon_{kj} | )& \mbox{otherwise}.
\end{array}
\right.
\end{align}

\item the skew-symmetrizers are related by
\begin{align}
\label{eq:d_prime_formula}
d'_i = d_i, \qquad \forall i.
\end{align}

\item the cluster variables are related by:
\begin{enumerate}
\item[\rm 1)] In the case of $\mcal{A}$-seeds: 
\begin{align}
\label{eq:mu_k_on_A}
A'_i = \left\{
\begin{array}{ll}
A_i & \mbox{if $i\neq k$}, \\
\displaystyle A_k^{-1} \left( \prod_{j \, | \, \varepsilon_{kj}>0 } A_j^{\varepsilon_{kj}} + \prod_{j \, | \, \varepsilon_{kj}<0} A_j^{-\varepsilon_{kj}} \right) & \mbox{if $i=k$},
\end{array}
\right.
\end{align}
where $\Gamma = (\varepsilon, d, \{A_i\}_{i=1}^n)$ and $\Gamma' = (\varepsilon', d', \{A_i'\}_{i=1}^n)$.

\vs

\item[\rm 2)] In the case of $\mcal{X}$-seeds: 
\begin{align}
\label{eq:mu_k_on_X}
X'_i = \left\{
\begin{array}{ll}
X_k^{-1} & \mbox{if $i=k$}, \\
X_i \left( 1 + X_k^{{\rm sgn}(-\varepsilon_{ik})} \right)^{-\varepsilon_{ik}} & \mbox{if $i\neq k$,}
\end{array}
\right.
\end{align}
where $\Gamma = (\varepsilon, d, \{X_i\}_{i=1}^n)$ and $\Gamma' = (\varepsilon', d', \{X_i'\}_{i=1}^n)$, and
$$
{\rm sgn}(a)= \left\{
\begin{array}{ll}
1 & \mbox{if $a\ge 0$} \\
-1 & \mbox{if $a<0$}.
\end{array}
\right.
$$

\vs

\item[\rm 3)] In the case of $\mcal{D}$-seeds: the formulas in  \eqref{eq:mu_k_on_X}, together with
\begin{align}
\label{eq:mu_k_on_B}
B'_i = \left\{
\begin{array}{ll}
B_i & \mbox{if $i\neq k$}, \\
\displaystyle \frac{ \left(\prod_{j \, | \, \varepsilon_{kj}<0 } B_j^{-\varepsilon_{kj}}\right) + X_k \left(\prod_{j\, | \, \varepsilon_{kj}>0} B_j^{\varepsilon_{kj}}\right) }{B_k(1+X_k)} & \mbox{if $i=k$,} 
\end{array}
\right.
\end{align}
where $\Gamma = (\varepsilon, d, \{B_i,X_i\}_{i=1}^n)$ and $\Gamma' = (\varepsilon', d', \{B_i',X_i'\}_{i=1}^n)$.
\end{enumerate}

\end{itemize}

We denote such a situation by $\mu_k(\Gamma) = \Gamma'$, or $\Gamma \stackrel{k}{\to} \Gamma'$. We call this procedure a \ul{\em mutation}.
\end{definition}

The new cluster variables are viewed as elements of the ambient field, the field of all rational functions in the previous cluster variables over $\mathbb{Q}$. Each of \eqref{eq:mu_k_on_A}--\eqref{eq:mu_k_on_B} should be thought of as an equality in the ambient field.

\begin{definition}[seed automorphism]
\label{def:seed_automorphism}
A seed $\Gamma' = (\varepsilon',d',*')$ is said to be obtained by applying the \\
\ul{\em seed automorphism $P_\sigma$} to a seed $\Gamma = (\varepsilon,d,*)$ for a permutation $\sigma$ of $\{1,2,\ldots,n\}$, if
$$
d'_{\sigma(i)} = d_i, \qquad \varepsilon'_{\sigma(i) \, \sigma(j)} = \varepsilon_{ij}, \qquad \forall i,j,
$$
and
$$
A'_{\sigma(i)} = A_i, \quad(\mbox{for $\mcal{A}$-seeds}), \qquad
X'_{\sigma(i)} = X_i, \quad(\mbox{for $\mcal{X}$-, $\mcal{D}$-seeds}), \qquad
B'_{\sigma(i)} = B_i, \quad(\mbox{for $\mcal{D}$-seeds}),
$$
for all $i=1,2,\ldots,n$. We denote such a situation by $P_\sigma(\Gamma) = \Gamma'$.
\end{definition}
We can consider the mutations and seed automorphisms as being applied to seeds from the left, and use the usual notation $\circ$ for the composition of them. So, when we apply a finite sequence of mutations and permutations to a seed, we read the sequence from right to left.

\begin{definition}[cluster transformations]
A \ul{\em cluster transformation} is a finite sequence of mutations and seed automorphisms, together with a seed $\Gamma$ on which the sequence is to be applied to. If the resulting seed is $\Gamma'$, we say that this cluster transformation \ul{\em connects} $\Gamma$ \ul{\em to} $\Gamma'$. We call mutations and seed automorphisms \ul{\em elementary cluster transformations}.
\end{definition}

Note that, each of the symbols $\mu_k$ and $P_\sigma$ stands for many different elementary cluster transformations, for it can be thought of as being applied to different seeds, which could be of any of the three kinds.

\subsection{Cluster modular groupoids}
\label{subsec:cluster_modular_groupoids}

In practice, we start from a single seed, which we often call an {\em initial seed}, and produce new seeds by applying cluster transformations to it. Then the set of all seeds created this way would be in correspondence with the set of all finite sequences of mutations and permutations. However, we will identify two seeds whenever they are essentially the same, in the sense I explain now.  Given any cluster transformation, say, connecting a seed $\Gamma$ to a seed $\Gamma'$, one obtains a well-defined identification of the cluster variables for $\Gamma'$ as rational functions in the cluster variables for $\Gamma$, by `composing' the formulas in Definitions \ref{def:mutation} and \ref{def:seed_automorphism}.
\begin{definition}[trivial cluster transformations]
\label{def:trivial_cluster_transformations}
A cluster transformation connecting a seed $\Gamma = (\varepsilon,d,*)$ to a seed $\Gamma'=(\varepsilon',d',*')$ is said to be \ul{\em weakly trivial} if $\varepsilon'=\varepsilon$ and $d'=d$. A cluster transformation is said to be \ul{\em $\mcal{A}$-trivial} if it is weakly trivial and induces the identity map between the cluster $\mcal{A}$-variables, i.e. the $i$-th cluster $\mcal{A}$-variable $A_i$ for $\Gamma$ equals the $i$-th cluster $\mcal{A}$-variable $A_i'$ for $\Gamma'$. The notions \ul{\em $\mcal{X}$-trivial} and \ul{\em $\mcal{D}$-trivial} are defined analogously. A cluster transformation is said to be \ul{\em trivial} if it is $\mcal{A}$-trivial, $\mcal{X}$-trivial, or $\mcal{D}$-trivial.
\end{definition}

\begin{remark}
In \cite{FG09} a `(feed) cluster transformation' means a sequence of mutations and permutations applied to a feed $(\varepsilon,d)$, and it is said to be `trivial' if it is $\mcal{A}$-trivial and $\mcal{X}$-trivial at the same time, when applied to $\mcal{A}$-seeds and $\mcal{X}$-seeds whose underlying feeds are $(\varepsilon,d)$.
\end{remark}

\begin{definition}[identification of seeds]
\label{def:identification_of_seeds}
If a cluster transformation connecting an $\mcal{A}$-seed $\Gamma$ to an $\mcal{A}$-seed $\Gamma'$ is $\mcal{A}$-trivial, we identify $\Gamma$ and $\Gamma'$ as $\mcal{A}$-seeds, and write $\Gamma = \Gamma'$. Likewise for $\mcal{X}$-seeds and $\mcal{D}$-seeds.
\end{definition}
Keeping this identification of seeds in mind, we consider:
\begin{definition}[equivalence of seeds]
\label{def:equivalence_of_seeds}
The two seeds are \ul{\em equivalent} if they are connected by a cluster transformation. For a seed $\Gamma$, denote by $|\Gamma|$ the equivalence class of $\Gamma$. 
\end{definition}

\begin{remark}
Identifying the seeds connected to each other by seed automorphisms amounts to considering Fomin-Zelevinsky's notion of `unlabeled' seeds.
\end{remark}

Let us now investigate some examples of trivial cluster transformations. It is a standard and straightforward exercise to show that the following example is indeed a trivial cluster transformation.
\begin{lemma}[involution identity of a mutation]
\label{lem:involution_identity}
$\mu_k \circ \mu_k$ is a trivial cluster transformation on any seed of any kind. That is, if we write $\Gamma' = \mu_k(\Gamma)$ and $\Gamma'' = \mu_k(\Gamma')$, with $\Gamma'' = (\varepsilon'',d'',*'')$ and $\Gamma = (\varepsilon,d,*)$, then $\varepsilon''=\varepsilon$, $d''=d$, and $A_i'' = A_i$, $\forall i$ (for $\mcal{A}$-seeds), $X_i'' = X_i'$, $\forall i$ (for $\mcal{X}$- and $\mcal{D}$-seeds), and $B_i'' = B_i$, $\forall i$ (for $\mcal{D}$-seeds). \qed
\end{lemma}
In particular, we identify $(\mu_k\circ\mu_k)(\Gamma)=\mu_k(\mu_k(\Gamma))$ and $\Gamma$. Some more examples of trivial cluster transformations involving the seed automorphisms are as follows, which I also omit a proof of.
\begin{lemma}[permutation identities]
\label{lem:permutation_identities}
For any permutations $\sigma,\gamma$ of $\{1,\ldots,n\}$ and any $k\in \{1,\ldots,n\}$,
\begin{align*}
P_\sigma \circ P_\gamma = P_{\sigma \circ \gamma}, \qquad
P_\sigma \circ \mu_k \circ P_{\sigma^{-1}} = \mu_{\sigma(k)}, \qquad
P_{{\rm Id}} = {\rm Id},
\end{align*}
all of which hold as identities when applied to any seed of any kind, where $\mathrm{Id}$ stands either for the identity permutation or the identity cluster transformation. \qed
\end{lemma}
The first two identities in this statement can be translated into saying that $P_\sigma \circ P_\gamma \circ P_{\sigma\circ \gamma}^{-1}$ and $P_\sigma \circ \mu_k \circ P_{\sigma^{-1}} \circ \mu_{\sigma(k)}^{-1}$ are trivial cluster transformations on any seed of any kind. 

\vs

There is another set of interesting trivial cluster transformations, pointed out e.g. in eq.(20) in of \cite[p.238]{FG09}; also including the `opposite' cases (with $i$ and $j$ switched), they are:
\begin{lemma}[The $(h+2)$-gon relations]
\label{lem:h_plus_2-gon_relations}
Suppose that a seed $\Gamma = (\varepsilon,d,*)$ of any kind satisfies
$$
\varepsilon_{ij} = - p \, \varepsilon_{ji} \in \{p,-p\} \qquad\mbox{or}\qquad
\varepsilon_{ji} = - p \, \varepsilon_{ij} \in \{p,-p\}
$$
for some $i,j\in \{1,\ldots,n\}$ with $i\neq j$ and $p\in \{0,1,2,3\}$. Let $h=2,3,4,6$ for $p=0,1,2,3$, respectively. Denote by $(i\, j)$ the transposition permutation of $\{1,\ldots,n\}$ interchanging $i$ and $j$. Then
\begin{align}
\mbox{$(P_{(i\, j)} \circ \mu_i)^{h+2}$ applied to $\Gamma$ is a trivial cluster transformation.}
\end{align}
\end{lemma}

The involution identity in Lem.\ref{lem:involution_identity} can be thought of as coming from rank one cluster algebras, namely, of Dynkin type $A_1$. The more serious identities in Lem.\ref{lem:h_plus_2-gon_relations} can be thought of as coming from rank two cluster algebras of Dynkin types $A_1\times A_1$, $A_2$, $B_2$, $G_2$, respectively, as pointed out in \cite{FG09}. Fock and Goncharov say in \cite{FG09} that they do not know how to find more examples of trivial cluster transformations that are not consequences of the ones already discussed, let alone how to describe the complete list of trivial cluster transformations. Thus they take only these known relations to formulate and solve the problem of establishing a quantum cluster variety.

\vs

An idealistic structure formed by the equivalent seeds is the following:
\begin{definition}
The \ul{\em $\mcal{A}$-cluster modular groupoid $\mcal{G}^\mcal{A} = \mcal{G}^\mcal{A}_{|\Gamma|}$} associated to an equivalence class $|\Gamma|$ of $\mcal{A}$-seeds is a category whose set of objects is $|\Gamma|$, and whose set of morphisms ${\rm Hom}_{\mcal{G}^\mcal{A}} (\Gamma, \Gamma')$ from an object $\Gamma$ to $\Gamma'$ is the set of all cluster transformations from $\Gamma$ to $\Gamma'$ modulo $\mcal{A}$-trivial cluster transformations.

\vs

Likewise for \ul{\em the $\mcal{X}$-cluster modular groupoid $\mcal{G}^\mcal{X}$} and \ul{\em the $\mcal{D}$-cluster modular groupoid $\mcal{G}^\mcal{D}$}.

\end{definition}
Here `modulo' means the following. Consider cluster transformations ${\bf m}_1, {\bf m}_2, {\bf m}_3$, where some of them may be empty sequences of elementary cluster transformations. If ${\bf m}_2$ is a $\mcal{A}$-trivial cluster transformation, then ${\bf m}_1 \circ {\bf m}_2 \circ {\bf m}_3$ and ${\bf m}_1 \circ {\bf m}_3$ are viewed as the same morphism in the category $\mcal{G}^\mcal{A}$. Likewise for $\mcal{G}^\mcal{X}$ and $\mcal{G}^\mcal{D}$. 

\vs

A slightly less idealistic version used in \cite{FG09} is the following, taking into account only the known relations.

\begin{definition}
\label{eq:saturated-modular_groupoid}
The \ul{\em saturated $\mcal{A}$-cluster modular groupoid $\wh{\mcal{G}}^\mcal{A} = \wh{\mcal{G}}^\mcal{A}_{|\Gamma|}$} associated to an equivalence class $|\Gamma|$ of $\mcal{A}$-seeds is a category whose set of objects is $|\Gamma|$, and whose set of morphisms ${\rm Hom}_{\wh{\mcal{G}}^\mcal{A}}(\Gamma, \Gamma')$ from an object $\Gamma$ to $\Gamma'$ is the set of all cluster transformations from $\Gamma$ to $\Gamma'$ modulo only the trivial cluster transformations that are described in Lemmas \ref{lem:involution_identity}, \ref{lem:permutation_identities} and \ref{lem:h_plus_2-gon_relations}. Likewise for $\wh{\mcal{G}}^\mcal{X}$ and $\wh{\mcal{G}}^\mcal{D}$.

\end{definition}

\begin{definition}
In these groupoids, an \ul{\em elementary morphism} is a morphism representing an elementary cluster transformation.
\end{definition}

As we shall soon see, these groupoids provide a handy way to formulate the construction of a cluster variety and its quantum counterpart. However, Fock and Goncharov \cite{FG09} define groupoids based on `feeds' $(\varepsilon,d)$, rather than seeds. Using theirs has an advantage as it allows to define the notion of (saturated) cluster modular group, or (saturated) cluster mapping class group\footnote{in the case of cluster variety coming from a punctured surface, this group coincides with the mapping class group of the surface.}, as the group of automorphisms of one feed; in particular, upon quantization, one gets projective representations of the cluster mapping class group.

\subsection{Cluster varieties}

Each elementary cluster transformation, say from a seed $\Gamma$ to a seed $\Gamma'$, induces a rational map from the seed torus for $\Gamma$ to that for $\Gamma'$, denoted by the same symbol as the elementary cluster transformation itself, defined in the level of functions by the identification formulas in Definitions \ref{def:mutation} and \ref{def:seed_automorphism} for the cluster variables. For example, we define the rational map $\mu_k : \mcal{X}_\Gamma \to \mcal{X}_{\Gamma'}$ between the seed $\mcal{X}$-tori by describing what the pullbacks of the coordinate functions of the torus $\mcal{X}_{\Gamma'}$ are, in terms of coordinate functions of $\mcal{X}_\Gamma$:
\begin{align}
\label{eq:mu_k_in_coorindates}
\mu_k^*(X'_i) = \left\{
\begin{array}{ll}
X_k^{-1} & \mbox{if $i=k$}, \\
X_i \left( 1 + X_k^{{\rm sgn}(-\varepsilon_{ik})} \right)^{-\varepsilon_{ik}} & \mbox{if $i\neq k$,}
\end{array}
\right.
\end{align}
and $\mu_k^*( {X_i'}^{-1} ) = {\mu_k^*(X_i')}^{-1}$, $\forall i$. As this map $\mu_k$ has a rational inverse map, namely what is denoted again by $\mu_k$ (Lem.\ref{lem:involution_identity}), we see that $\mu_k$ is in fact a birational map. It is easy to see that the map $P_\sigma$ between two seed tori is an isomorphism of varieties.  Now that each of the symbols $\mu_k$ and $P_\sigma$ represents a birational map, so does a sequence of such symbols, i.e. so does any cluster transformation. We finally construct cluster varieties by gluing the seed tori along these birational maps.
\begin{definition}[cluster varieties]
\label{def:cluster_varieties}
The \ul{\em cluster $\mcal{A}$-variety} for an equivalence class $|\Gamma|$ of $\mcal{A}$-seeds, is a scheme obtained by gluing all the seed tori $\mcal{A}_{\Gamma'}$ for $\Gamma'$ in this equivalence class $|\Gamma|$ using the birational maps associated to cluster transformations. We denote it by $\mcal{A}_{|\Gamma|}$.

\vs

Define the \ul{\em cluster $\mcal{X}$-variety} $\mcal{X}_{|\Gamma|}$ and the \ul{\em cluster $\mcal{D}$-variety} $\mcal{D}_{|\Gamma|}$ in a similar fashion.
\end{definition}
In order for $\mcal{A}_{|\Gamma|}$ to be well-defined, for each pair of seeds $(\Gamma', \Gamma'')$ in $|\Gamma|$ we must have a unique birational map $\mcal{A}_{\Gamma'} \to \mcal{A}_{\Gamma''}$ along which we glue the two tori. There may be many different sequences of mutations and permutations that connect $\Gamma$ to $\Gamma'$, and we should make sure that they induce the same birational map on the tori; it is an easy consequence of Def.\ref{def:identification_of_seeds}.

\vs

A cluster variety defined in Def.\ref{def:cluster_varieties} can be formulated as a contravariant functor
\begin{align}
\label{eq:eta}
\eta : \mbox{the cluster modular groupoid} ~ \mcal{G}_{|\Gamma|} ~ \longrightarrow ~ \mbox{certain category of commutative rings},
\end{align}
where $\mcal{G}_{|\Gamma|}$ stands for one of $\mcal{G}^\mcal{A}_{|\Gamma|}$, $\mcal{G}^\mcal{X}_{|\Gamma|}$, and $\mcal{G}^\mcal{D}_{|\Gamma|}$, and one may want to use the saturated version $\wh{\mcal{G}}_{|\Gamma|}$ instead of $\mathcal{G}_{|\Gamma|}$. In the target category, morphisms between objects are homomorphisms between the fields of fractions of the objects. For example, for a cluster $\mathcal{A}$-variety, the object $\eta(\Gamma)$ is the Laurent polynomial ring $\mathbb{Z}[\{A_i, A_i^{-1} : i=1,\ldots,n \}]$, and the morphism $\eta(\mu_k)$ is the map $\mu_k : {\rm Frac}(\eta(\Gamma')) \to {\rm Frac}(\eta(\Gamma))$ described in eq.\eqref{eq:mu_k_in_coorindates}, where ${\rm Frac}$ means the field of fractions. This groupoid formulation of the cluster varieties is not just an unnecessary luxury to have, for such a formulation in terms of groupoids and rings is the only known way to describe the quantum versions of cluster varieties; there is no actual topological space in the quantum world.

\vs

The following geometric structures, which are extra data on the cluster varieties, are crucial in the story of quantization.
\begin{lemma}[geometric structures on cluster varieties]
\label{lem:geometric_structures_on_cluster_varieties}
The geometric structures on the seed tori defined in Def.\ref{def:geometric_structures_on_seed_tori} induce well-defined corresponding geometric structures on the respective cluster varieties.
\end{lemma}
One checks this lemma by verifying that $\mu_k$ and $P_\sigma$ preserve the geometric structures on the seed tori.

\section{Representations of quantum cluster varieties}

\subsection{Quantum cluster varieties}

Let us first briefly review what the words quantum and quantization mean. Let $M$ be a Poisson manifold; so $M$ is a smooth manifold, and there is a Poisson bracket $\{\cdot,\cdot\}$ on the ring of real-valued smooth functions $C^\infty(M)$. A {\em deformation quantization} of $M$ refers to a set of data consisting of a complex Hilbert space $\mathscr{H}$, and an assignment
$$
\rho^\hbar : C^\infty(M) \to \{\mbox{self-adjoint operators $\mathscr{H} \to \mathscr{H}$}\}
$$
that depends real analytically on a real parameter $\hbar$, such that
\begin{enumerate}
\item[\rm (Q1)] $\rho^\hbar$ is $\mathbb{R}$-linear,

\item[\rm (Q2)] $\rho^\hbar(1) = {\rm Id}_\mathscr{H}$,

\item[\rm (Q3)] $[\rho^\hbar(f), \rho^\hbar(g)] = 2\pi{\rm i} \hbar \, \rho^\hbar(\{f,g\}) + o(\hbar)$ as $\hbar\to 0$, \qquad $\forall f,g\in C^\infty(M)$,
\end{enumerate}
often also with certain irreducibility condition. The commutative ring $C^\infty(M)$ equipped with the Poisson bracket $\{\cdot,\cdot\}$ is viewed as representing a classical system of observables, while the operators on the Hilbert space $\mathscr{H}$ equipped with the operator commutator $[\cdot,\cdot]$ is viewed as representing a quantum system of observables. A quantization is then a map from classical to quantum in the above sense, and a quantization problem is to find such a map.

\vs

Let us now break down the quantization problem into easier ones. As usual, only some Poisson subalgebra $R$ of $C^\infty(M)$ is dealt with; in our case, the Laurent polynomial ring plays the role of $R$. Then, the sought-for map $\rho^\hbar$ is factored into the composition of a map $R\to R^\hbar$ and $R^\hbar \to \{\mbox{operators on $\mathscr{H}$}\}$, where $R^\hbar$ is a non-commutative ring that deforms the classical commutative ring $R$ in a certain sense. In particular, $R^\hbar$ should recover the original ring $R$ when $\hbar=0$, and the former map $R\to R^\hbar$ should recover the identity map when $\hbar=0$, and must satisfy some algebraic analogs of the conditions (Q1), (Q2), and (Q3), involving the Poisson bracket on $R$. The latter map $R^\hbar \to \{\mbox{operators on $\mathscr{H}$}\}$ is just a representation of the ring $R^\hbar$ on $\mathscr{H}$. So, establishing a non-commutative ring $R^\hbar$ and its representation amounts to establishing a quantum system, while finding a map $R\to R^\hbar$ amounts to finding a quantization. As in most of the literature on quantum Teichm\"uller theory and quantum cluster varieties including \cite{FG09}, we focus only on building a quantum system in the present paper, instead of finding a quantization map. 

\vs

When it comes to quantizing a Poisson variety, we do not just have one ring $R$, but bunch of rings $R_\Gamma$ labeled by the data $\Gamma$. For each $\Gamma$, one would establish a non-commutative ring $R^\hbar_\Gamma$ and its representation. Meanwhile, for each pair $(\Gamma,\Gamma')$, the respective rings $R_\Gamma$ and $R_{\Gamma'}$ are related through the isomorphism $\mu^*_{\Gamma,\Gamma'} : {\rm Frac}(R_{\Gamma'}) \to {\rm Frac}(R_\Gamma)$; one would like to enhance this to a quantum version, i.e. an isomorphism $\mu^\hbar_{\Gamma,\Gamma'} : {\rm Frac}(R^\hbar_{\Gamma'}) \to {\rm Frac}(R^\hbar_\Gamma)$ between skew-fields of fractions of the quantum rings. The quantum isomorphism $\mu^\hbar_{\Gamma,\Gamma'}$ must recover the classical isomorphism $\mu^*_{\Gamma,\Gamma'}$ when $\hbar=0$, and must satisfy the consistency equations, i.e.
$$
\mu^\hbar_{\Gamma,\Gamma'} \circ \mu^\hbar_{\Gamma',\Gamma''}  = \mu^\hbar_{\Gamma,\Gamma''}
$$
for every triple $(\Gamma,\Gamma',\Gamma'')$. Thus, a \ul{\em quantum cluster variety} for a cluster variety defined in eq.\eqref{eq:eta} can be formulated as a contravariant functor
\begin{align}
\nonumber
\eta^\hbar : \mbox{the cluster modular groupoid} ~ \mcal{G}_{|\Gamma|} ~ \longrightarrow ~ \mbox{certain category of non-commutative rings}.
\end{align}
A saturated version would be a contravariant functor
\begin{align}
\nonumber
\eta^\hbar : \mbox{the saturated cluster modular groupoid} ~ \wh{\mcal{G}}_{|\Gamma|} ~ \longrightarrow ~ \mbox{certain category of non-commutative rings},
\end{align}
and this is how it is formulated in \cite{FG09}. Finding such a functor $\eta^\hbar$ is already a highly non-trivial task, but still one more step remains, namely building a representation. For each $\Gamma$, one would like to construct a Hilbert space $\mathscr{H}_\Gamma$, and a representation $\pi^\hbar_\Gamma$ of the algebra $\eta^\hbar(\Gamma)=R^\hbar_\Gamma$ on $\mathscr{H}_\Gamma$. Per each pair $(\Gamma,\Gamma')$, one would like to build a unitary map ${\bf K}^\hbar_{\Gamma,\Gamma'} : \mathscr{H}_{\Gamma'} \to \mathscr{H}_\Gamma$ between the respective Hilbert spaces that intertwine the two representations $\pi^\hbar_{\Gamma'}$ and $\pi^\hbar_\Gamma$, in the sense that the diagram
$$
\xymatrix@C+10mm{
\mathscr{H}_{\Gamma'} \ar[r]^-{\pi^\hbar_{\Gamma'}(u)} \ar[d]_{{\bf K}^\hbar_{\Gamma,\Gamma'}} & \mathscr{H}_{\Gamma'} \ar[d]^{{\bf K}^\hbar_{\Gamma,\Gamma'}} \\
\mathscr{H}_\Gamma \ar[r]^-{\pi^\hbar_\Gamma(\mu^\hbar_{\Gamma,\Gamma'}(u))} & \mathscr{H}_\Gamma
}
$$
commutes for all $u\in R^\hbar_{\Gamma'}$, i.e. the {\em intertwining equations} hold
$$
\pi^\hbar_\Gamma ( \mu^\hbar_{\Gamma,\Gamma'}(u) ) \circ {\bf K}^\hbar_{\Gamma,\Gamma'}  = {\bf K}^\hbar_{\Gamma,\Gamma'} \circ \pi^\hbar_{\Gamma'}(u)
$$
for all $u\in R^\hbar_{\Gamma'}$. In fact, a lot more care is needed to make sense of the intertwining equations; first, the equation makes sense only when both $u\in R^\hbar_{\Gamma'}$ and $\mu^\hbar_{\Gamma,\Gamma'} (u) \in R^\hbar_{\Gamma}$ hold, and second, one should carefully keep track of the domain of each densely-defined operator appearing in the equations. Let us not go into such issues here, as they are not relevant for the purpose of the present paper; see \cite{FG09} for details. Besides the intertwining equations, the intertwiner operators ${\bf K}^\hbar_{\Gamma,\Gamma'}$ must also satisfy the consistency equations, this time allowed up to constants:
$$
{\bf K}^\hbar_{\Gamma,\Gamma'} \circ {\bf K}^\hbar_{\Gamma',\Gamma''} = c^\hbar_{\Gamma,\Gamma',\Gamma''} \, {\bf K}^\hbar_{\Gamma,\Gamma''}
$$
for each triple $(\Gamma,\Gamma',\Gamma'')$, for some constants $c^\hbar_{\Gamma,\Gamma',\Gamma''} \in {\rm U}(1)\subset \mathbb{C}$; these constants are referred to as the \ul{\em phase constants} in the present paper. Thus, a {\em representation of a quantum cluster variety} can be formulated as a contravariant projective functor
\begin{align}
\label{eq:K_functor}
\pi^\hbar : \mbox{the cluster modular groupoid} ~ \mcal{G}_{|\Gamma|} ~ \longrightarrow ~ {\rm Hilb},
\end{align}
or a saturated version
\begin{align}
\label{eq:K_functor2}
\pi^\hbar : \mbox{the saturated cluster modular groupoid} ~ \wh{\mcal{G}}_{|\Gamma|} ~ \longrightarrow ~ {\rm Hilb},
\end{align}
where ${\rm Hilb}$ is the category of Hilbert spaces, whose morphisms are unitary maps.

\vs

The purpose of the present paper is to show that, for the Fock-Goncharov representation of a quantum cluster $\mathcal{D}$-variety constructed in \cite{FG09}, the phase constants $c^\hbar_{\Gamma,\Gamma',\Gamma''}$ are all $1$. So we will recall the construction of the Hilbert spaces $\mathscr{H}_\Gamma$ and the unitary intertwiners ${\bf K}^\hbar_{\Gamma,\Gamma'}$; for more details such as the construction of the quantum cluster variety functor $\eta^\hbar$, a reader can consult the original paper \cite{FG09} and references therein.

\vs

The Hilbert space $\mathscr{H}_\Gamma$ associated to a $\mathcal{D}$-seed $\Gamma$ is easy to describe. As shall be seen soon, it is convenient to also define a nice dense subspace of it.

\begin{definition}
\label{def:H_Gamma_and_D_Gamma}
For a $\mcal{D}$-seed $\Gamma = (\varepsilon,d,\{B_i,X_i\}_{i=1}^n)$, define the Hilbert space
\begin{align}
\label{eq:H_Gamma}
\mathscr{H}_\Gamma := L^2(\mathbb{R}^n, \, da_1\, da_2 \, \ldots \, da_n),
\end{align}
where the measure on $\mathbb{R}^n$ is the product of Lebesgue measures $da_i$ on $\mathbb{R}$, and the inner product $\langle \,\, , \, \rangle_{\mathscr{H}_\Gamma}$ is as usual:
$$
\langle f, g \rangle_{\mathscr{H}_\Gamma} = \int_{\mathbb{R}^n} f \, \ol{g} \, da_1 \ldots da_n.
$$
Define $D_\Gamma$ as the following $\cplx$-vector subspace of $\mathscr{H}_\Gamma$: the row vector of variables being denoted as $\vec{a} = (a_1 \, a_2 \, \cdots \, a_n)$,
\begin{align}
\label{eq:D_Gamma}
\hspace{-3mm}
D_\Gamma := {\rm span}_\cplx\left\{ ~ e^{ \vec{a} M \vec{a}^{\rm t} + \vec{v} \cdot \vec{a} } \, P(\vec{a}) ~ \, \left| \, ~ \begin{array}{l} 
\mbox{$M \in {\rm Mat}_{n\times n}(\mathbb{C})$ with the real part ${\rm Re}(M)$ being negative-definite,} \\
\mbox{$\vec{v} = (v_1\cdots v_n) \in \mathbb{C}^n$, and $P$ is a polynomial in $a_1,\ldots,a_n$ over $\mathbb{C}$}
\end{array} \right. \right\}.
\end{align}

\end{definition}
\begin{remark}
The symbols $a_1,\ldots,a_n$ are used in \cite{FG09} as if they are ${\rm log}$ versions of the $\mathcal{A}$-variables $A_1,\ldots,A_n$; it is not clear whether this viewpoint is useful or valid, but let us live with it at the moment.
\end{remark}

\begin{remark}
The above version of $D_\Gamma$ is used in \cite{KS}; this is enlarged from Fock-Goncharov's subspace $W_{\bf i}$ \cite{FG09}, in order to make it invariant under more examples of  unitary operators used in proofs.
\end{remark}

The strategy to construct the intertwiners ${\bf K}^\hbar_{\Gamma,\Gamma'}$ is to construct them only for the elementary morphisms of the cluster groupoid, i.e. when $\Gamma'$ is either $\mu_k(\Gamma)$ or $P_\sigma(\Gamma)$. For a general morphism, i.e. a general cluster transformation from $\Gamma$ to $\Gamma'$, the operator ${\bf K}^\hbar_{\Gamma,\Gamma'}$ is constructed as the composition of the operators for the elementary cluster transformations. Then, what must be checked is whether trivial cluster transformations correspond to scalar operators under this construction; that is, we should check whether the intertwiner operators for elementary cluster transformations satisfy the relations analogous to the ones in Lemmas \ref{lem:involution_identity}, \ref{lem:permutation_identities}, and \ref{lem:h_plus_2-gon_relations}.

\subsection{Self-adjoint operators for ${\rm log}$ coordinate functions}

As the first step towards the construction of the unitary intertwiner operator ${\bf K}^\hbar_{\Gamma,\Gamma'}$, I will describe the self-adjoint operators ${\bf b}^\hbar_i$ and ${\bf x}^\hbar_i$ on $\mathscr{H}_\Gamma$, whose exponentials would be the quantum operators associated to the classical coordinate functions $B_i$ and $X_i$ of a $\mathcal{D}$-seed $\Gamma$.

\begin{definition}[the basic self-adjoint operators]
\label{def:Schrodinger_representation}
Let $\Gamma$, $\mathscr{H}_\Gamma$ and $D_\Gamma$ be as in Def.\ref{def:H_Gamma_and_D_Gamma}. Let $\hbar \in \mathbb{R}_{>0}$. Define the operators ${\bf q}^\hbar_i$, ${\bf p}^\hbar_i$ : $D_\Gamma \to \mathscr{H}_\Gamma$, $i=1,\ldots,n$, by the formulas
\begin{align}
\label{eq:Schrodinger_representation}
{\bf q}^\hbar_i := a_i \qquad\mbox{and}\qquad {\bf p}^\hbar_i := 2\pi {\rm i} \, \hbar \, \frac{\partial}{\partial a_i} \qquad \mbox{on}\quad D_\Gamma,
\end{align}
for each $i=1,\ldots,n$, i.e.
$$
({\bf q}^\hbar_i f)(a_1,\ldots,a_n) = a_i \, f(a_1,\ldots,a_n), \qquad
({\bf p}^\hbar_i f)(a_1,\ldots,a_n) = 2\pi {\rm i} \, \hbar \, \frac{\partial f}{\partial a_i} (a_1,\ldots,a_n)
$$
for all $f(a_1,\ldots,a_n) \in D_\Gamma$.  Denote their unique self-adjoint extensions in $\mathscr{H}_\Gamma$ also by the same symbols ${\bf q}^\hbar_i$'s and ${\bf p}^\hbar_i$'s respectively, by abuse of notation.
\end{definition}
The above definition makes sense because the operators ${\bf q}^\hbar_i$ and ${\bf p}^\hbar_i$ defined on $D_\Gamma$ are essentially self-adjoint in $\mathscr{H}_\Gamma$. Various well-known facts about these operators can be easily found in the literature, e.g. in \cite{RS70} \cite{Hall}. They satisfy the relations
\begin{align}
\label{eq:Heisenberg_relations_of_p_and_q}
[{\bf p}^\hbar_i, \, {\bf q}^\hbar_j] = 2\pi {\rm i} \hbar \, \delta_{ij} \cdot \mathrm{Id}, \qquad
[{\bf p}^\hbar_i, \, {\bf p}^\hbar_j]=0, \qquad
[{\bf q}^\hbar_i, \, {\bf q}^\hbar_j]=0,
\end{align}
say, as operators $D_\Gamma \to D_\Gamma$. In fact, much stronger relations than eq.\eqref{eq:Heisenberg_relations_of_p_and_q} hold. 
\begin{definition}
Let $A,B$ be self-adjoint operators on a Hilbert space $\mathscr{H}$, and let $c\in \mathbb{R}$. Consider the relation
$$
[A,B] = {\rm i} c \cdot {\rm Id}
$$
understood only formally; such a relation is called a \ul{\em Heisenberg relation}. The \ul{\em Weyl relations} for this Heisenberg relation refer to the following family of identities of unitary operators
$$
e^{ {\rm i} \alpha A} \, e^{{\rm i} \beta B} = e^{-{\rm i} c \alpha\beta} \, e^{{\rm i} \beta B} \, e^{ {\rm i} \alpha A}, \qquad \forall \alpha, \beta\in \mathbb{R},
$$
where $e^{ {\rm i} \alpha A}$ and $e^{{\rm i} \beta B}$ are defined via the functional calculus for $A$ and $B$.
\end{definition}

\begin{lemma}
The Weyl relations for the Heisenberg relations in eq.\eqref{eq:Heisenberg_relations_of_p_and_q} are satisfied in $\mathscr{H}_\Gamma$. \qed
\end{lemma}
The Weyl relations are so strong that they actually characterize these self-adjoint operators.
\begin{proposition}[Stone-von Neumann theorem; {\cite{vN31} \cite[Thm.VIII.14]{RS70} \cite[Thm.14.8]{Hall}}]
\label{prop:SvN}
Let $\mathscr{H}$ be a separable complex Hilbert space, and $A_1,\ldots,A_n, B_1,\ldots,B_n$ be possibly unbounded self-adjoint operators on $\mathscr{H}$. Let $\hbar \in \mathbb{R}$. Suppose that the Weyl relations for the Heisenberg relations
$$
[A_i, B_j] = 2\pi {\rm i} \hbar \, \delta_{ij} \cdot {\rm Id}, \qquad
[B_i, B_j ] = 0, \qquad
[A_i, A_j] = 0
$$
hold, for all $i,j$. Then there are closed subspaces $\mathscr{H}_\ell$, $\ell = 1,\ldots,N$, where $N$ is a positive integer or $\infty$, such that
\begin{enumerate}
\item[\rm (1)] $\mathscr{H}$ is the orthogonal direct sum $\mathscr{H} = \bigoplus_{\ell=1}^N \mathscr{H}_\ell$,

\item[\rm (2)] For every $\alpha,\beta\in\mathbb{R}$ and every $i,j,\ell$, the unitary operators $e^{{\rm i}\alpha A_i}$ and $e^{{\rm i}\beta B_j}$ preserve $\mathscr{H}_\ell$, and the restrictions of the self-adjoint operators $A_i$ and $B_j$ to $\mathscr{H}_\ell$ also preserve $\mathscr{H}_\ell$,

\item[\rm (3)] For each $i,j,\ell$, the restrictions $A_i \restriction \mathscr{H}_\ell$ and $B_j \restriction \mathscr{H}_\ell$ are self-adjoint in $\mathscr{H}_\ell$, and there exists a unitary operator $U_\ell : \mathscr{H}_\ell \to L^2(\mathbb{R}^n, da_1 \cdots da_n)$ such that
$$
U_\ell \, (A_i\restriction \mathscr{H}_\ell) \, U_\ell^{-1} = {\bf p}^\hbar_i, \qquad
U_\ell \, (B_j\restriction \mathscr{H}_\ell) \, U_\ell^{-1} = {\bf q}^\hbar_j, \qquad \forall i,j=1,\ldots,n,
$$
hold as equalities of self-adjoint operators. Such $U_\ell$ is unique up to multiplication by a scalar in ${\rm U}(1)$, i.e. a complex number of modulus $1$. 

\item[\rm (4)] If $\mathscr{H}$ is irreducible, in the sense that the only closed subspaces of $\mathscr{H}$ that are invariant under every $e^{{\rm i} \alpha A_i}$ and every $e^{{\rm i}\beta B_j}$ (for $i,j=1,\ldots,n$, $\alpha,\beta \in \mathbb{R}$) are $0$ and $\mathscr{H}$, then $N=1$.
\end{enumerate}

\end{proposition}

It is well-known that any $\mathbb{R}$-linear combination of ${\bf q}^\hbar_i$'s and ${\bf p}^\hbar_i$'s is an essentially self-adjoint operator on $D_\Gamma$, hence yields a unique self-adjoint extension in $\mathscr{H}_\Gamma$; see e.g. \cite[\S14]{Hall}. Any two such self-adjoint operators satisfy the Weyl relations for a suitable Heisenberg relation; such a Heisenberg relation can be computed from the basic ones in eq.\eqref{eq:Heisenberg_relations_of_p_and_q} using bilinearity of the commutator. 

\begin{definition}[the log quantum operators for coordinate functions; \cite{FG07}]
\label{def:old_representation}
Let $\Gamma$, $\mathscr{H}_\Gamma$, $D_\Gamma$, and $\hbar$ be as in Def.\ref{def:Schrodinger_representation}. Let
\begin{align}
\label{eq:hbar_i}
\hbar_i := \hbar/d_i, \qquad \forall i=1,\ldots,n.
\end{align}
Define the self-adjoint operators ${\bf b}^\hbar_i$ and ${\bf x}^\hbar_i$ in $\mathscr{H}_\Gamma$ as the following $\mathbb{R}$-linear combinations of ${\bf q}^\hbar_i$'s and ${\bf p}^\hbar_i$'s
\begin{align}
\label{eq:old_representation}
{\bf b}^\hbar_i := 2{\bf q}_i^\hbar, \qquad
\displaystyle {\bf x}^\hbar_i := \frac{1}{2} d_i^{-1} \, {\bf p}^\hbar_i - \sum_{j=1}^n \varepsilon_{ij} \, {\bf q}^\hbar_j, \qquad \forall i=1,\ldots,n,
\end{align}
and define the tilde version operators $\til{\bf x}^\hbar_i$ as
\begin{align}
\label{eq:old_representation_tilde}
\displaystyle \til{\bf x}^\hbar_i := {\bf x}^\hbar_i + \sum_{j=1}^n \varepsilon_{ij} \, {\bf b}^\hbar_j
= \frac{1}{2} d_i^{-1} \, {\bf p}^\hbar_i + \sum_{j=1}^n \varepsilon_{ij} \, {\bf q}^\hbar_j, \qquad \forall i =1,\ldots,n.
\end{align} 
\end{definition}

Then the Weyl relations for the following Heisenberg relations hold:
\begin{align}
\label{eq:Heisenberg_relations}
[{\bf b}^\hbar_i, {\bf b}^\hbar_j] = 0, \qquad
[{\bf x}^\hbar_i, {\bf b}^\hbar_j] = 2\pi {\rm i} \hbar_i \, \delta_{i,j} \cdot {\rm Id}, \qquad
[{\bf x}^\hbar_i, {\bf x}^\hbar_j] = 2 \pi {\rm i} \hbar \, \wh{\varepsilon}_{ij} \cdot {\rm Id};
\end{align}
notice that these relations are quantum versions of the classical Poisson bracket relations in eq.\eqref{eq:Poisson_on_D}. For the tilde version operators, the Weyl relations for the following Heisenberg relations hold:
\begin{align}
\label{eq:Heisenberg_relations_involving_tildes}
[\til{\bf x}^\hbar_i, {\bf b}^\hbar_j] = 2\pi {\rm i} \hbar_i \, \delta_{i,j} \cdot {\rm Id}, \qquad
[\til{\bf x}^\hbar_i, \til{\bf x}^\hbar_j] = - 2 \pi {\rm i} \hbar \, \wh{\varepsilon}_{ij} \cdot {\rm Id}, \qquad
[\til{\bf x}^\hbar_i, {\bf x}^\hbar_j] = 0,
\end{align}
which correspond to eq.\eqref{eq:Poisson_on_D_tilde}. Notice that the operators constructed so far are for each single chosen $\mathcal{D}$-seed $\Gamma$. So, when necessary, I will write
$$
{\bf x}^\hbar_{\Gamma;i} = {\bf x}^\hbar_i, \quad
{\bf b}^\hbar_{\Gamma;i} = {\bf b}^\hbar_i, \quad
\til{\bf x}^\hbar_{\Gamma;i} = \til{\bf x}^\hbar_i, \quad
{\bf q}^\hbar_{\Gamma;i} = {\bf q}^\hbar_i, \quad
{\bf p}^\hbar_{\Gamma;i} = {\bf p}^\hbar_i,
$$
to emphasize which $\Gamma$ we are working with. It is wise to keep track of the superscript $\hbar$, for we will eventually be dealing with operators associated to different $\hbar$'s at the same time.

\subsection{Special affine shift operators} 
\label{subsec:special_affine_shift_operators}

For construction of the sought-for mutation intertwiner operator, we need some ingredients, reviewed in the present and the next subsections. Let $(M,\mu)$ and $(N,\nu)$ be measure spaces, and $\phi : M \to N$ be an invertible measure-preserving map whose inverse is also measure-preserving. It induces a natural unitary isomorphism of Hilbert spaces $L^2(M,d\mu) \to L^2(N,d\nu)$ given by $f \mapsto (\phi^{-1})^* f = f \circ \phi^{-1}$. In the present subsection we investigate some special cases of these operators, when both measure spaces are the usual Euclidean space $(\mathbb{R}^n, da_1 \, \ldots \, da_n)$ and $\phi: \mathbb{R}^n \to \mathbb{R}^n$ is of the following type.
\begin{definition}
Write each element of $\mathbb{R}^n$ as a row vector $\vec{a} = (a_1 \, a_2 \, \cdots \, a_n)$ with real entries.  A bijective map $\phi : \mathbb{R}^n \to \mathbb{R}^n$ is called an \ul{\emph{special affine transformation}} if 
\begin{align}
\label{eq:phi_inverse}
\phi^{-1}(\vec{a}) = \vec{a} \, \mathbf{c} + \vec{t}, \quad \forall \vec{a} \in \mathbb{R}^n.
\end{align}
holds for some matrix $\mathbf{c} = (c_{ij})_{i,j\in\{1,\ldots,n\}} \in \mathrm{SL}_\pm(n,\mathbb{R})$ and $\vec{t} \in \mathbb{R}^n$, where
\begin{align}
\mathrm{SL}_\pm(n,\mathbb{R}) := \{ \, \mathbf{c} \in \mathrm{GL}(n,\mathbb{R}) \, : \, |\det \mathbf{c}|=1 \, \}.
\end{align}

\end{definition}

The notions and results in the present subsection are nothing fundamentally new, but I find it convenient to establish a notation and some lemmas.

\begin{lemma}
The above correspondence $\phi \leadsto (\mathbf{c},\vec{t}\,)$ is a group isomorphism between the group of all special affine transformations of $\mathbb{R}^n$ and the following semi-direct product group
\begin{align}
\mathrm{SL}_\pm(n,\mathbb{R}) \ltimes \mathbb{R}^n = \{ (\mathbf{c}, \vec{t}\,) \, | \, \mathbf{c}\in \mathrm{SL}_\pm(n,\mathbb{R}), \,\, \vec{t} \in \mathbb{R}^n \},
\end{align}
whose multiplication is given by
\begin{align}
\label{eq:multiplication_of_semidirect_product}
(\mathbf{c}, \vec{t}\,) \, (\mathbf{c}', \vec{t}\,') = (\mathbf{c} \, \mathbf{c}', \vec{t} \, \mathbf{c}' + \vec{t}\,'). \qed
\end{align}
\end{lemma}

Let us adapt the notation for the Hilbert space $\mathscr{H}_\Gamma$ and its dense subspace $D_\Gamma$, defined in eq.\eqref{eq:H_Gamma}--\eqref{eq:D_Gamma} in Def.\ref{def:H_Gamma_and_D_Gamma}. Note that the $\mathcal{D}$-seed data $\Gamma$, except for its rank $n$, is not playing any role in the definitions of $\mathscr{H}_\Gamma$ and $D_\Gamma$; it does only when we consider multiple $\mathcal{D}$-seeds at the same time, which is not the case in the present subsection. So, in principle we can drop the subscript $\Gamma$ for now, but I keep it in order to avoid any unnecessary confusion.

\begin{lemma}
For any special affine transformation $\phi$, the map $f \mapsto (\phi^{-1})^* f = f \circ \phi^{-1}$ is a unitary isomorphism from the Hilbert space $\mathscr{H}_\Gamma$ \eqref{eq:H_Gamma} to itself. \qed
\end{lemma}

\begin{definition}
\label{def:special_affine_shift_operator}
For a special affine transformation $\phi$, given by $\phi^{-1}(\vec{a}) = \vec{a} \, \mathbf{c} + \vec{t}$, denote the unitary automorphism $f\mapsto (\phi^{-1})^*f$ of the Hilbert space $\mathscr{H}_\Gamma$ \eqref{eq:H_Gamma} by  $\mathbf{S}_{\Gamma;(\mathbf{c},\vec{t}\,)}$ or $\mathbf{S}_{(\mathbf{c},\vec {t}\,)}$, and call it a \ul{\em special affine shift operator} on $\mathscr{H}_\Gamma$.
\end{definition}
More explicitly, $\mathbf{S}_{(\mathbf{c},\vec{t}\,)}$ is given by
\begin{align}
\label{eq:def_of_S_c_t}
  (\mathbf{S}_{(\mathbf{c},\vec{t}\,)}f)(\vec{a}) = f(\vec{a} \, \mathbf{c} + \vec{t}\,), \qquad \forall f\in \mathscr{H}_\Gamma=L^2(\mathbb{R}^n, da_1\cdots da_n), \quad \forall \vec{a}\in \mathbb{R}^n.
\end{align}

\begin{lemma}
\label{lem:special_affine_shift_operators_multiplicative}
The correspondence $(\mathbf{c}, \vec{t}\,) \leadsto \mathbf{S}_{(\mathbf{c}, \vec{t}\,)}$ is a group homomorphism, i.e.
\begin{align}
\label{eq:bf_S_is_group_homomorphism}
  \mathbf{S}_{(\mathbf{c}, \vec{t}\,) (\mathbf{c}',\vec{t}\,'\,)} = \mathbf{S}_{(\mathbf{c}, \vec{t}\,)} \, \mathbf{S}_{(\mathbf{c}', \vec{t}\,')}, \qquad\quad \mathbf{S}_{(\mathbf{Id}, \vec{0}\,)} = \mathrm{Id}_{\mathscr{H}_\Gamma},
\end{align}
and is injective. \qed
\end{lemma}

\begin{corollary}
\label{cor:special_affine_shift_operators_form_a_group}
The special affine shift operators on $\mathscr{H}_\Gamma$ \eqref{eq:H_Gamma} form a group, isomorphic to $\mathrm{SL}_\pm(n,\mathbb{R})\ltimes \mathbb{R}^n$. \qed
\end{corollary}

From the definition in eq.\eqref{eq:def_of_S_c_t} we immediately obtain the following observation, which becomes useful in the proof of the main result of the present paper:
\begin{lemma}[scalar special affine shift operator is identity]
\label{lem:scalar_special_affine_shift_operator_is_identity}
If $\mathbf{S}_{(\mathbf{c},\vec{t}\, )} = c \cdot \mathrm{Id}_{\mathscr{H}_\Gamma}$ for some complex scalar $c$, then $c=1$. \qed
\end{lemma}
The following lemma is also immediate:
\begin{lemma}
\label{lem:bf_S_preserves_D}
$\mathbf{S}_{(\mathbf{c},\vec{t}\,)} (D_\Gamma) = D_\Gamma$.  \qed
\end{lemma}

An easy example of a special affine shift operator is:
\begin{definition}
\label{def:P_sigma}
For a permutation $\sigma$ of $\{1,\ldots,n\}$, define the \ul{\em permutation operator} $\mathbf{P}_\sigma = {\bf P}_{\Gamma;\sigma}$ on $\mathscr{H}_\Gamma$ \eqref{eq:H_Gamma} as 
$$
(\mathbf{P}_\sigma f)(a_1,\ldots,a_n) = f(a_{\sigma(1)},\ldots,a_{\sigma(n)}), \qquad \forall f\in \mathscr{H}_\Gamma.
$$
\end{definition}

Another class of examples is given by the usual shift operators $\mathbf{S}_{(\mathbf{Id},\vec{t}\,)}$: for any $\vec{t} = (t_1 \, t_2 \, \cdots \, t_n)\in \mathbb{R}^n$, observe
\begin{align}
\label{eq:bf_S_as_shift_operator}
(\mathbf{S}_{(\mathbf{Id}, \vec{t}\,)}f)(\vec{a}) = f(\vec{a}+\vec{t}\,). 
\end{align}
In view of the Taylor series expansion, one may expect that we could write
\begin{align}
\label{eq:bf_S_as_Taylor_series}
\mathbf{S}_{(\mathbf{Id},\vec{t}\,)} = e^{\sum_{i=1}^n t_i \frac{\partial}{\partial a_i}}, 
\end{align}
which a priori makes sense only formally, but can also be made precise by functional calculus of several mutually commuting self-adjoint operators; we need not prove it here, as eq.\eqref{eq:bf_S_as_Taylor_series} will only be used to give intuition. In fact, we will mostly use examples where $\vec{t} = \vec{0}$, hence let us introduce the following abbreviated notation:
\begin{align}
\label{eq:S_bf_c}
{\bf S}_{\Gamma;{\bf c}} = {\bf S}_{\bf c} := {\bf S}_{({\bf c},\vec{0}\,)}.
\end{align}
The following basic lemma on the conjugation action of ${\bf S}_{\bf c}$ on ${\bf p}^\hbar_i$'s and ${\bf q}^\hbar_i$'s will be useful.
\begin{lemma}
\label{lem:bf_S_conjugation_on_bf_p_and_bf_q}
Let $\mathbf{c} \in \mathrm{SL}_\pm(n,\mathbb{R})$ and consider ${\bf S}_{\Gamma;{\bf c}} = \mathbf{S}_{\bf c} : \mathscr{H}_\Gamma \to \mathscr{H}_\Gamma$. Then, for each $i=1,\ldots,n$,
\begin{align}
\label{eq:conjugation_of_bf_S_on_bf_p}
  \mathbf{S}_{\bf c} \, \mathbf{p}^\hbar_i \,\, \mathbf{S}_{\bf c}^{-1} & = \sum_{j=1}^n c^{ij} \, \mathbf{p}^\hbar_j, \quad \mbox{where} \quad \mathbf{c}^{-1} = (c^{ij})_{i,j\in\{1,\ldots,n\}}, \\
\label{eq:conjugation_of_bf_S_on_bf_q}
  \mathbf{S}_{\bf c} \,\, \mathbf{q}^\hbar_i \,\, \mathbf{S}_{\bf c}^{-1} & = \sum_{j=1}^n c_{ji} \, \mathbf{q}^\hbar_j, \quad \mbox{where} \quad \mathbf{c} = (c_{ij})_{i,j\in\{1,\ldots,n\}},
\end{align}
hold as equalities of self-adjoint operators. \qed
\end{lemma}
For example, these equations can be directly checked when applied to elements $\psi\in D_\Gamma$; viewing both sides of the equations as operators $D_\Gamma \to D_\Gamma$, one obtains the desired equalities by taking the unique self-adjoint extensions. Another way of proving these equations of self-adjoint operators is to consider and prove the exponentiated version, using the strongly continuous one-parameter unitary groups generated by the self-adjoint operators; see e.g. \cite{Hall}.

\subsection{Quantum dilogarithm function}
\label{subsec:non-compact_QD}

Another crucial ingredient of the construction of intertwiners is the special function called the quantum dilogarithm studied by Faddeev and Kashaev \cite{FK94} \cite{F95}. I first recall the `compact' version.
\begin{definition}[the compact quantum dilogarithm]
For any complex number $\mathbf{q}$ with $|\mathbf{q}|<1$ define a merophorphic function $\Psi^\mathbf{q}(z)$ on the complex plane $\mathbb{C}$, called the \ul{\em compact quantum dilogarithm}, as
\begin{align}
\label{eq:Psi_q}
\Psi^\mathbf{q}(z) := \prod_{i=1}^\infty (1 + \mathbf{q}^{2i-1} z)^{-1} = \frac{1}{(1+\mathbf{q}z)(1+\mathbf{q}^3 z)(1+\mathbf{q}^5 z) \cdots}
\end{align}
\end{definition}
\begin{lemma}
The above infinite product absolutely converges for any $z$ except at the poles $z = \mathbf{q}^{-(2i-1)}$, $i=1,2,3,\ldots$. \qed
\end{lemma}
\begin{lemma}
One has the functional equation
\begin{align}
\label{eq:functional_equation_of_compact_QD}
  \Psi^\mathbf{q}(\mathbf{q}^2 z) = (1+\mathbf{q}z) \, \Psi^\mathbf{q}(z). \qed
\end{align}
\end{lemma}
This compact quantum dilogarithm function lies behind the description of the image of morphisms under the (algebraic) quantum cluster variety functor $\eta^\hbar$, where the functional equation eq.\eqref{eq:functional_equation_of_compact_QD} plays a crucial role. However, what must be used is the case when ${\bf q} = e^{\pi {\rm i} \hbar}$, in particular $|{\bf q}|=1$. Then the above infinite product in eq.\eqref{eq:Psi_q} does not converge, and such expression can only be taken heuristically. This does not make a serious problem for the functor $\eta^\hbar$, for the functor can be described with well-defined formulas not involving the compact quantum dilogarithm. But this does become an issue for a representation $\pi^\hbar$ of quantum cluster variety, for the images of morphisms should directly involve such functions.

\vs

A remedy to this problem is the `non-compact' quantum dilogarithm; among its many guises, here I use a version used in \cite{FG09}, which is defined as a certain contour integral, as shall be seen soon. Another viewpoint is that the non-compact quantum dilogarithm can be thought of as a suitable limit of the compact quantum dilogarithm as $|\mathbf{q}|\to 1$. A na\"ive limit does not exist; only certain ratio of compact quantum dilogarithm does. Note first that if we write $\mathbf{q}=e^{\pi{\rm i}h}$ for $h\in \mathbb{C}$ (not $\hbar \in \mathbb{R}$!), then $|\mathbf{q}|<1$ if and only if $\mathrm{Im}(h)>0$, and $|\mathbf{q}|=1$ means $\mathrm{Im}(h)=0$. As $z\mapsto -1/z$ is an example of a $\mathrm{PSL}(2,\mathbb{R})$ M\"obius transformation of the upper half-plane, we see that $\mathrm{Im}(h)>0$ implies $\mathrm{Im}(-1/h)>0$. So, if we define $\mathbf{q}^\vee := e^{\pi {\rm i}/h}$ for $\mathrm{Im}(h)>0$, we have $| 1/\mathbf{q}^\vee |<1$. When we form the ratio $\Psi^\mathbf{q}(e^z)/\Psi^{1/\mathbf{q}^\vee}(e^{z/h})$, the poles and zeros make appropriate cancellation, so that the limit as $\mathrm{Im}(h) \to 0$ exists; this limit will yield the non-compact quantum dilogarithm $\Phi^\hbar(z)$ with $\hbar\in \mathbb{R}$, which is what is actually used in our construction. 

\vs

Now I recall the definition and some basic properties of the non-compact quantum dilogarithm $\Phi^h$, in a somewhat more general form than we need.
\begin{lemma}
For any $h\in \mathbb{C}$ with $\mathrm{Im}(h)\ge 0$ and $\mathrm{Re}(h)>0$, let
$$
\mathbf{q}:=e^{\pi{\rm i} h}, \qquad \mathbf{q}^\vee := e^{\pi{\rm i}/h}.
$$
\begin{enumerate}
\item[\rm (1)] 
The integral in the expression
$$
\Phi^h(z) = \exp\left( - \frac{1}{4} \int_\Omega \frac{e^{-{\rm i}pz}}{\sinh(\pi p) \sinh(\pi h p)} \frac{dp}{p} \right),
$$
where $\Omega$ is a contour on the real line that avoids the origin via a small half-circle above the origin, absolutely converges for $z$ in the strip $|\mathrm{Im}(z)| < \pi(1 + \mathrm{Re}(h))$, and yields a non-vanishing complex analytic function $\Phi^h(z)$ on this strip. Each of the functional equations
\begin{align}
\label{eq:Phi_h_difference_equations}
\left\{ {\renewcommand{\arraystretch}{1.3} \begin{array}{lcl}
\Phi^h(z+2\pi{\rm i} h) & = & (1+\mathbf{q}e^z) \, \Phi^h(z), \\
\Phi^h(z+2\pi{\rm i}) & = & (1+\mathbf{q}^\vee e^{z/h}) \, \Phi^h(z),
\end{array}} \right.
\end{align}
called the {\em difference equations}, holds when the two arguments of $\Phi^h$ are in the strip. These functional equations let us analytically continue $\Phi^h$ to a meromorphic function on the whole plane $\mathbb{C}$, with
\begin{align*}
\mbox{the set of zeros} & = \{ \, (2\ell+1) \pi {\rm i} +  (2m+1) \pi {\rm i} h \,\, | \,\, \ell,m\in \mathbb{Z}_{\ge 0} \, \}, \quad \mbox{and} \\
\mbox{the set of poles} & =  \{ \, - (2\ell+1) \pi {\rm i} - (2m+1) \pi {\rm i} h \,\, | \,\, \ell,m\in \mathbb{Z}_{\ge 0} \, \}.
\end{align*}

\item[\rm (2)] \emph{(relationship between compact and non-compact)} When $\mathrm{Im}(h)>0$, one has the equality
\begin{align}
\label{eq:Phi_h_as_ratio}
\frac{\Psi^\mathbf{q}(e^z)}{\Psi^{1/\mathbf{q}^\vee}(e^{z/h})} = \Phi^h(z),
\end{align}
for any $z$ that is not a pole of $\Phi^h$.

\item[\rm (3)] When $h = \hbar \in \mathbb{R}_{>0}$, every pole and zero of $\Phi^\hbar$ is simple if and only if $\hbar \notin \mathbb{Q}$.

\item[\rm (4)] \emph{($\hbar \leftrightarrow 1/\hbar$ duality)} When $h = \hbar \in \mathbb{R}_{>0}$, one has
\begin{align}
\label{eq:QD_self-duality}
  \Phi^{1/\hbar}(z/\hbar) = \Phi^\hbar(z).
\end{align}

\item[\rm (5)] \emph{(unitarity)} When $h = \hbar \in \mathbb{R}_{>0}$, one has
$$
|\Phi^\hbar(z)|=1, \qquad \forall z\in \mathbb{R}.
$$

\item[\rm (6)] \emph{(involutivity)} One has
\begin{align}
\label{eq:QD_identity_quadratic}
\Phi^h(z) \Phi^h(-z) = c_h\, e^{ z^2/(4\pi {\rm i}h) },
\end{align}
where
\begin{align}
\label{eq:c_hbar}
c_h := e^{- \frac{\pi {\rm i}}{12} (h + h^{-1})} \in \mathbb{C}^\times.
\end{align}
In particular, $|c_\hbar|=1$ for $\hbar = h \in \mathbb{R}_{>0}$. \qed
\end{enumerate}
\end{lemma}
From now on, by the \ul{\em non-compact quantum dilogarithm} we mean the meromorphic function $\Phi^\hbar(z)$, with $\hbar \in \mathbb{R}_{>0}$. The integral expression for $\Phi^\hbar(z)$ goes back to \cite{B01}.

\vs

Besides the difference equations eq.\eqref{eq:Phi_h_difference_equations}, the quantum dilogarithm functions satisfy a crucial equation called the pentagon equation. For the compact quantum dilogarithm $\Psi^{\bf q}$, it reads:

\begin{proposition}[the pentagon equation for the compact quantum dilogarithm; \cite{FK94}]
\label{prop:pentagon_equation_compact}
Let $X,Y$, and ${\bf q}$ be formal variables satisfying $XY = \mathbf{q}^2YX$, where ${\bf q}$ commutes with $X$ and $Y$. Then one has
\begin{align}
\label{eq:compact_QD_pentagon}
  \Psi^\mathbf{q}(Y)^{-1} \, \Psi^\mathbf{q}(X)^{-1} = \Psi^\mathbf{q}(X)^{-1} \, \Psi^\mathbf{q}(\mathbf{q}^{-1}XY)^{-1} \, \Psi^\mathbf{q}(Y)^{-1},
\end{align}
where each factor and each side is regarded as a formal power series in the variables $X$ and $Y$ with coefficients in $\mathbb{Z}_{\ge 0}[\mathbf{q}] \subset \mathbb{Z}[\mathbf{q}]$. \qed
\end{proposition}
The pentagon equation for the non-compact quantum dilogarithm $\Phi^\hbar$ is best understood as the identity of unitary operators on a Hilbert space; this shall be dealt with in later sections, independent on Prop.\ref{prop:pentagon_equation_compact}.

\subsection{Construction of intertwiners}
\label{subsec:formulas_for_intertwiners}

I finally present Fock-Goncharov's representation of quantum cluster variety, namely a functor $\pi^\hbar$ as in eq.\eqref{eq:K_functor} or \eqref{eq:K_functor2}. In the object level, we already described the answer, namely, for each $\mathcal{D}$-seed $\Gamma$, let
$$
\pi^\hbar(\Gamma) := \mathscr{H}_\Gamma,
$$
defined in eq.\eqref{eq:H_Gamma}. Now in the morphism level, let us describe the image of elementary morphisms. First, consider the morphism $\mu_k$, i.e. a mutation (Def.\ref{def:mutation}) from $\Gamma$ to $\mu_k(\Gamma) = \Gamma'$. We assign the unitary map
$$
\pi^\hbar(\mu_k) := {\bf K}^\hbar_k =  \mathbf{K}^\hbar_{\Gamma\mut{k}\Gamma'} ~ : ~ \mathscr{H}_{\Gamma'} \to \mathscr{H}_\Gamma, \quad \mbox{when} \quad \Gamma' = \mu_k(\Gamma),
$$
defined as follows. 
\begin{definition}[the intertwiner for a mutation; {\cite[Def.5.1]{FG09}}]
\label{def:bf_K}
Let $\Gamma = (\varepsilon,d,\{B_i,X_i\}_{i=1}^n)$ be a $\mathcal{D}$-seed, and let $\mu_k(\Gamma) = \Gamma' = (\varepsilon',d',\{B_i',X_i'\}_{i=1}^n)$ for some $k\in \{1,\ldots,n\}$. Let $\hbar \in \mathbb{R}_{>0}$. 

\vs

Define the unitary  map $\mathbf{K}^\hbar_{\Gamma\mut{k}\Gamma'} : \mathscr{H}_{\Gamma'} \to \mathscr{H}_\Gamma$ as
\begin{align}
\label{eq:bf_K_definition}
\mathbf{K}^\hbar_{\Gamma\mut{k}\Gamma'} := \mathbf{K}^{\sharp \hbar}_{\Gamma\mut{k}\Gamma'} \circ \mathbf{K}'_{\Gamma\mut{k}\Gamma'},
\end{align}
where the \ul{\em automorphism part}
\begin{align}
\label{eq:sharp_K_definition}
\mathbf{K}^{\sharp \hbar}_{\Gamma\mut{k}\Gamma'} := \Phi^{\hbar_k}(\mathbf{x}^\hbar_{\Gamma;k}) \left( \Phi^{\hbar_k}(\til{\mathbf{x}}^\hbar_{\Gamma;k}) \right)^{-1}  ~ : ~ \mathscr{H}_\Gamma \to \mathscr{H}_\Gamma
\end{align}
is the composition of two unitary operators obtained by applying the functional calculus to the self-adjoint operators $\mathbf{x}^\hbar_{\Gamma;k}$ and $\til{\mathbf{x}}^\hbar_{\Gamma;k}$ (Def.\ref{def:old_representation}) respectively for the functions $z\mapsto \Phi^{\hbar_k}(z)$ and $z\mapsto (\Phi^{\hbar_k}(z))^{-1}$, while the \ul{\em linear part} is the unitary map
\begin{align}
\label{eq:bf_K_prime}
  \mathbf{K}'_{\Gamma\mut{k}\Gamma'} := \mathbf{S}_{\Gamma; \mathbf{c}_{\Gamma\mut{k}\Gamma'}} \circ \mathbf{I}_{\Gamma,\Gamma'}: \mathscr{H}_{\Gamma'} \to \mathscr{H}_\Gamma,
\end{align}
where $\mathbf{I}_{\Gamma,\Gamma'} : \mathscr{H}_{\Gamma'}\to \mathscr{H}_\Gamma$ is induced by identifying each $a_i'$ with $a_i$, i.e.
\begin{align}
\label{eq:bf_I}
(\mathbf{I}_{\Gamma,\Gamma'} \, f)(a_1,\ldots,a_n) := f(a_1',\ldots,a_n'), \qquad \forall f \in \mathscr{H}_{\Gamma'} = L^2(\mathbb{R}^n, da_1'\, \cdots da_n'),
\end{align}
and $\mathbf{S}_{\Gamma;\mathbf{c}_{\Gamma\mut{k}\Gamma'}} : \mathscr{H}_\Gamma \to \mathscr{H}_\Gamma$ is the special affine shift operator defined as in eq.\eqref{eq:S_bf_c} and eq.\eqref{eq:def_of_S_c_t} for the matrix
$$
\mathbf{c}_{\Gamma\mut{k}\Gamma'} = (c_{ij})_{i,j\in \{1,\ldots,n\}} \in \mathrm{SL}_\pm(n,\mathbb{R})
$$
given by
\begin{align}
\label{eq:c_Gamma_Gamma_prime_entries}
c_{ii} = 1, \quad \forall i \neq k, \qquad c_{kk} = -1, \qquad c_{ik} = [-\varepsilon_{ki}]_+, \quad \forall i\neq k, \qquad c_{ij}=0 \quad\mbox{otherwise}.
\end{align}
\end{definition}

\begin{remark}
As written in eq.(73) of \cite{FG09}, one can view $\mathbf{K}'_{\Gamma\mut{k}\Gamma'} : \mathscr{H}_{\Gamma'} \to \mathscr{H}_\Gamma$ as induced by the following map between two $\mathbb{R}^n$'s:
\begin{align}
\label{eq:what_bf_K_prime_does}
a_i' \mapsto \left\{
\begin{array}{ll}
a_i & \mbox{if $i\neq k$}, \\
-a_k + \sum_{j=1}^n [-\varepsilon_{kj}]_+ \, a_j  & \mbox{if $i=k$.}
\end{array}
\right.
\end{align}
\end{remark}

\vs

For the morphism $P_\sigma$, i.e. a seed automorphism (Def.\ref{def:seed_automorphism}) from $\Gamma$ to $P_\sigma(\Gamma)=\Gamma'$ where $\sigma$ is a permutation of $\{1,\ldots,n\}$, we assign the permutation operator in Def.\ref{def:P_sigma}. More precisely, let
\begin{align}
\label{eq:bf_K_P_sigma}
\pi^\hbar(P_\sigma) := \mathbf{P}_{\Gamma;\sigma} \circ {\bf I}_{\Gamma,\Gamma'} ~ : ~ \mathscr{H}_{\Gamma'} \to \mathscr{H}_\Gamma,
\end{align}
where ${\bf I}_{\Gamma,\Gamma'} : \mathscr{H}_{\Gamma'} \to \mathscr{H}_\Gamma$ is as defined in eq.\eqref{eq:bf_I}, and ${\bf P}_{\Gamma;\sigma}={\bf P}_\sigma : \mathscr{H}_\Gamma \to \mathscr{H}_\Gamma$ as in Def.\ref{def:P_sigma}.

\vs

Now, consider any morphism ${\bf m}$ in the (saturated) cluster modular groupoid. So ${\bf m}$ is a sequence of mutations and seed automorphisms. Reverse this sequence; each member of the sequence is assigned a unitary operator as described so far in the present subsection. Take the composition of these operators in this (reversed) order. The resulting operator is what we assign as $\pi^\hbar({\bf m})$. This completes the construction of Fock-Goncharov representation $\pi^\hbar$ of quantum cluster $\mathcal{D}$-variety. For $\pi^\hbar$ constructed this way to make sense as a well-defined projective functor, one must check that the operators $\pi^\hbar(\mu_k)$ and $\pi^\hbar(P_\sigma)$ which we assigned to the elementary morphisms satisfy the operator identities corresponding to the relations in Lemmas \ref{lem:involution_identity}, \ref{lem:permutation_identities} and \ref{lem:h_plus_2-gon_relations}, up to constants. This is the topic of the following section.

\section{Operator identities}

\subsection{Fock-Goncharov identities for the intertwiners, and the main theorem} A major result of \cite{FG09} is that the intertwining operators $\pi^\hbar(\mu_k) = {\bf K}^\hbar_k$ and $\pi^\hbar(P_\sigma)$ constructed in \S\ref{subsec:formulas_for_intertwiners} indeed provides a representation of quantum cluster $\mathcal{D}$-variety. As mentioned earlier, one aspect of this statement is that they satisfy the operator identities corresponding to the trivial cluster transformations:
\begin{proposition}[{\cite[Thm.5.4]{FG09}}]
\label{thm:FG_constant}
For any trivial morphism $\mathbf{m}$ in the saturated cluster modular groupoid $\wh{\mcal{G}}^\mcal{D}_{|\Gamma|}$ from a $\mathcal{D}$-seed $\Gamma$ to itself written as a sequence of mutations and seed automorphisms, if $\pi^\hbar(\mathbf{m})$ denotes the composition of the reversed sequence of intertwiners corresponding to the mutations and seed automorphisms (constructed in \S\ref{subsec:formulas_for_intertwiners})  constituting the sequence $\mathbf{m}$, then
$$
\pi^\hbar(\mathbf{m}) = c_{\mathbf{m}} \cdot \mathrm{Id}_{\mathscr{H}_\Gamma},
$$
for some constant $c_{\mathbf{m}} \in \mathrm{U}(1) \subset \mathbb{C}^\times$. \qed
\end{proposition}

\begin{corollary}
The intertwiners $\pi^\hbar(\mu_k) = {\bf K}^\hbar_k$ and $\pi^\hbar(P_\sigma)$ constructed in \S\ref{subsec:formulas_for_intertwiners} for the elementary cluster transformations induce a well-defined projective functor $\pi^\hbar$ as formulated in eq.\eqref{eq:K_functor2}. \qed
\end{corollary}
The constant $c_{\mathbf{m}}$ is denoted by $\lambda$ in \cite[Thm.5.4--5.5]{FG09}, and is not precisely determined as an explicit number depending on $\hbar$. A major goal of the present paper is to show by computation that this constant $c_{\mathbf{m}}$ is always $1$. For this, I first re-write the result of the above theorem more explicitly, written for each of the generating relations of types (S1), (S2-1), (S2-2), (S2-3), and (S2-4) in \S\ref{sec:introduction}, or as in Lemmas \ref{lem:involution_identity}, \ref{lem:permutation_identities} and \ref{lem:h_plus_2-gon_relations}.

\vs

Most basic is the involution identity:
\begin{proposition}[rank $1$ identity; $A_1$ identity; `twice-flip is identity'; \cite{FG09}]
\label{prop:the_rank_1_identity_for_intertwiners}
Let $\Gamma=(\varepsilon,d,*)$ be a $\mathcal{D}$-seed, $k \in \{1,\ldots,n\}$, $\Gamma' = \mu_k(\Gamma) = (\varepsilon',d',*')$. Then $\Gamma = \mu_k(\Gamma')$, and there exists a constant $c_{A_1} \in \mathrm{U}(1)\subset \mathbb{C}^\times$ such that
\begin{align}
\label{eq:twice_flip_to_prove}
\mathbf{K}^\hbar_{\Gamma\mut{k}\Gamma'} \, \mathbf{K}^\hbar_{\Gamma'\mut{k}\Gamma}  = c_{A_1} \cdot \mathrm{Id}_{\mathscr{H}_\Gamma}.
\end{align}
\end{proposition}
Appropriately understood, this can be written compactly as
$$
{\bf K}^\hbar_k \circ {\bf K}^\hbar_k = c_{A_1} \cdot {\rm Id}, \quad\mbox{or}\quad \pi^\hbar(\mu_k) \circ \pi^\hbar(\mu_k) = c_{A_1} \cdot {\rm Id}.
$$

\vs

The following form of the permutation identities, corresponding to Lem.\ref{lem:permutation_identities}, are straightforward:
\begin{proposition}[permutation identities for the intertwiners]
\label{prop:permutation_identities_for_the_intertwiners}
One has
$$
\pi^\hbar(P_\gamma) \, \pi^\hbar(P_\sigma) = \pi^\hbar(P_{\sigma\circ \gamma}), \qquad
\pi^\hbar(P_{\sigma^{-1}}) \, {\bf K}^\hbar_k \,
\pi^\hbar(P_\sigma) = {\bf K}^\hbar_{\sigma(k)}, \qquad
\pi^\hbar(P_{\rm Id}) = {\rm Id}.
$$
\end{proposition}
Written in the style of Prop.\ref{thm:FG_constant}, these can be written as follows, with the help of Prop.\ref{prop:the_rank_1_identity_for_intertwiners}
$$
\pi^\hbar(P_\gamma) \, \pi^\hbar(P_\sigma)  \,\pi^\hbar(P_{\sigma\circ \gamma}^{-1}) = {\rm Id}, \qquad
\pi^\hbar(P_{\sigma^{-1}}) \, {\bf K}^\hbar_k \,
\pi^\hbar(P_\sigma) \,
{\bf K}^\hbar_{\sigma(k)} = c_{A_1} \cdot {\rm Id}, \qquad
\pi^\hbar(P_{\rm Id}) = {\rm Id}.
$$

\vs

The remaining ones, i.e. the $(h+2)$-gon relations, may be written in a compact form as follows. Let $i,j\in \{1,\ldots,n\}$ be two distinct indices with $\varepsilon_{ij} = - p \varepsilon_{ji} = \pm p$ or $\varepsilon_{ji} = - p \varepsilon_{ij} = \pm p$, with $p \in \{0,1,2,3\}$. Let $h=2,3,4,6$ for $p=0,1,2,3$ respectively. Then, \cite[Thm.5.5]{FG09} says that
$$
({\bf P}_{(i\,j)} \circ {\bf K}^\hbar_i)^{h+2} = \lambda \cdot {\rm Id}
$$
holds for some $\lambda\in {\rm U}(1)\subset \mathbb{C}^\times$ depending on $p$. To understand these identities more precisely, one has to be careful about the underlying seeds. Before writing down all the underlying seeds, we will first expand the left-hand-side, and use Prop.\ref{prop:permutation_identities_for_the_intertwiners} to simplify.

\begin{proposition}[rank 2 identities for the elementary intertwiners; {\cite[Thm.5.4--5.5]{FG09}}]
\label{prop:rank_2_identities_for_intertwiners}
One has

\begin{enumerate}
\item[\rm (2-1)] ($A_1\times A_1$ identity; `commuting identity') Let $\Gamma^{(0)} = (\varepsilon^{(0)}, d^{(0)}, *^{(0)})$ be a $\mathcal{D}$-seed, and assume that some two distinct indices $i,j\in \{1,\ldots,n\}$ satisfy
\begin{align}
\label{eq:A1_times_A1_condition}
\varepsilon^{(0)}_{ij} = 0 = \varepsilon^{(0)}_{ji}.
\end{align}
Let
$$
\Gamma^{(1)} := \mu_i (\Gamma^{(0)}), \quad \Gamma^{(2)} := \mu_j(\Gamma^{(1)}), \quad \Gamma^{(3)}:=\mu_i(\Gamma^{(2)}).
$$
Then $\Gamma^{(0)} = \mu_j(\Gamma^{(3)})$, and there exists a constant $c_{A_1\times A_1} \in \mathrm{U}(1) \subset\mathbb{C}^\times$ such that
\begin{align}
\label{eq:A1_times_A1_to_prove}
\mathbf{K}^\hbar_{\Gamma^{(0)}\mut{i}\Gamma^{(1)}} \, \mathbf{K}^\hbar_{\Gamma^{(1)}\mut{j}\Gamma^{(2)}} \, \mathbf{K}^\hbar_{\Gamma^{(2)}\mut{i}\Gamma^{(3)}} \, \mathbf{K}^\hbar_{\Gamma^{(3)}\mut{j}\Gamma^{(0)}} = c_{A_1 \times A_1} \cdot {\rm Id}_{\mathscr{H}_{\Gamma^{(0)}}}
\end{align}

\vs

\item[\rm (2-2)] ($A_2$ identity; `pentagon identity') Let $\Gamma^{(0)} = (\varepsilon^{(0)}, d^{(0)}, *^{(0)})$ be a $\mathcal{D}$-seed, and assume that some $i,j\in \{1,\ldots,n\}$ satisfy
$$
\varepsilon^{(0)}_{ij} = - \, \varepsilon^{(0)}_{ji} \in \{1,-1\}.
$$
Let
$$
\Gamma^{(1)} := \mu_i (\Gamma^{(0)}), \quad \Gamma^{(2)} := \mu_j(\Gamma^{(1)}),\quad \Gamma^{(3)}:=\mu_i(\Gamma^{(2)}), \quad \Gamma^{(4)}:=\mu_j(\Gamma^{(3)}), \quad \Gamma^{(5)}:=\mu_i(\Gamma^{(4)}).
$$
Then $\Gamma^{(0)} = P_{(i\, j)}(\Gamma^{(5)})$, and there exists a constant $c_{A_2} \in \mathrm{U}(1) \subset \mathbb{C}^\times$ such that 
\begin{align}
\label{eq:A2_to_prove}
\mathbf{K}^\hbar_{\Gamma^{(0)}\mut{i}\Gamma^{(1)}} \, \mathbf{K}^\hbar_{\Gamma^{(1)}\mut{j}\Gamma^{(2)}} \, \mathbf{K}^\hbar_{\Gamma^{(2)}\mut{i}\Gamma^{(3)}} \, \mathbf{K}^\hbar_{\Gamma^{(3)}\mut{j}\Gamma^{(4)}} \, \mathbf{K}^\hbar_{\Gamma^{(4)}\mut{i}\Gamma^{(5)}} \, \mathbf{P}_{\Gamma^{(5)}; (i\, j)} \, \,{\bf I}_{\Gamma^{(5)},\Gamma^{(0)}} = c_{A_2} \cdot {\rm Id}_{\mathscr{H}_{\Gamma^{(0)}}}.
\end{align}

\vs

\item[\rm (2-3)] ($B_2$ identity; `hexagon identity') Let $\Gamma^{(0)} = (\varepsilon^{(0)}, d^{(0)}, *^{(0)})$ be a $\mathcal{D}$-seed, and assume that some $i,j\in \{1,\ldots,n\}$ satisfy
$$
\varepsilon^{(0)}_{ij} = - 2 \, \varepsilon^{(0)}_{ji} \in \{2,-2\}, \qquad \mbox{or} \qquad
\varepsilon^{(0)}_{ji} = - 2 \, \varepsilon^{(0)}_{ij} \in \{2,-2\}.
$$
Let
\begin{align*}
& \Gamma^{(1)} := \mu_i (\Gamma^{(0)}), \quad \Gamma^{(2)} := \mu_j(\Gamma^{(1)}),\quad \Gamma^{(3)}:=\mu_i(\Gamma^{(2)}), \quad
\Gamma^{(4)}:=\mu_j(\Gamma^{(3)}), \quad \Gamma^{(5)}:=\mu_i(\Gamma^{(4)}).
\end{align*}
Then $\Gamma^{(0)} = \mu_j(\Gamma^{(5)})$, and there exists a constant $c_{B_2} \in \mathrm{U}(1) \subset \mathbb{C}^\times$ such that 
\begin{align}
\label{eq:B2_to_prove}
\mathbf{K}^\hbar_{\Gamma^{(0)}\mut{i}\Gamma^{(1)}} \, \mathbf{K}^\hbar_{\Gamma^{(1)}\mut{j}\Gamma^{(2)}} \, \mathbf{K}^\hbar_{\Gamma^{(2)}\mut{i}\Gamma^{(3)}} \, \mathbf{K}^\hbar_{\Gamma^{(3)}\mut{j}\Gamma^{(4)}} \, \mathbf{K}^\hbar_{\Gamma^{(4)}\mut{i}\Gamma^{(5)}} \, \mathbf{K}^\hbar_{\Gamma^{(5)}\mut{j}\Gamma^{(0)}} = c_{B_2} \cdot {\rm Id}_{\mathscr{H}_{\Gamma^{(0)}}}.
\end{align}

\vs

\item[\rm (2-4)] ($G_2$ identity; `octagon identity') Let $\Gamma^{(0)} = (\varepsilon^{(0)}, d^{(0)}, *^{(0)})$ be a $\mathcal{D}$-seed, and assume that some $i,j\in \{1,\ldots,n\}$ satisfy
$$
\varepsilon^{(0)}_{ij} = - 3 \, \varepsilon^{(0)}_{ji} \in \{3,-3\}, \qquad\mbox{or}\qquad
\varepsilon^{(0)}_{ji} = - 3 \, \varepsilon^{(0)}_{ij} \in \{3,-3\}.
$$
Let
\begin{align*}
& \Gamma^{(1)} := \mu_i (\Gamma^{(0)}), \quad \Gamma^{(2)} := \mu_j(\Gamma^{(1)}),\quad \Gamma^{(3)}:=\mu_i(\Gamma^{(2)}), \quad \Gamma^{(4)} := \mu_j(\Gamma^{(3)}), \\
& \Gamma^{(5)}:=\mu_i(\Gamma^{(4)}), \quad \Gamma^{(6)}:=\mu_j(\Gamma^{(5)}), \quad \Gamma^{(7)} := \mu_i(\Gamma^{(6)}).
\end{align*}
Then $\Gamma^{(0)} = \mu_j(\Gamma^{(7)})$, and there exists a constant $c_{G_2} \in \mathrm{U}(1) \subset \mathbb{C}^\times$ such that 
\begin{align}
\label{eq:G2_to_prove}
{\renewcommand{\arraystretch}{1.4}\begin{array}{l}
\mathbf{K}^\hbar_{\Gamma^{(0)}\mut{i}\Gamma^{(1)}} \, \mathbf{K}^\hbar_{\Gamma^{(1)}\mut{j}\Gamma^{(2)}} \, \mathbf{K}^\hbar_{\Gamma^{(2)}\mut{i}\Gamma^{(3)}} \, \mathbf{K}^\hbar_{\Gamma^{(3)}\mut{j}\Gamma^{(4)}} \, \mathbf{K}^\hbar_{\Gamma^{(4)}\mut{i}\Gamma^{(5)}} \, \mathbf{K}^\hbar_{\Gamma^{(5)}\mut{j}\Gamma^{(6)}} \, \mathbf{K}^\hbar_{\Gamma^{(6)}\mut{i}\Gamma^{(7)}} \, \mathbf{K}^\hbar_{\Gamma^{(7)}\mut{j}\Gamma^{(0)}} 
= c_{G_2} \cdot  {\rm Id}_{\mathscr{H}_{\Gamma^{(0)}}}.
\end{array}}
\end{align}
\end{enumerate}

\end{proposition}

\vs

The main theorem of the present paper is that the constants appearing in the above operator identities are all $1$.

\begin{theorem}[main theorem]
\label{thm:main}
Above constants $c_{A_1}$, $c_{A_1\times A_1}$, $c_{A_2}$, $c_{B_2}$, and $c_{G_2}$ are all $1$. Hence the constants $c_{\bf m}$ in Prop.\ref{thm:FG_constant} are all $1$.
\end{theorem}

In fact, a priori, it is not clear from the arguments of \cite{FG09} whether these constants are uniform in the sense as written in Prop.\ref{prop:rank_2_identities_for_intertwiners}; for example, for the case (2-2), it is not obvious that the constant for the case $\varepsilon_{ij} = 1$ is same as that for $\varepsilon_{ij}=-1$. In the end, this is not a problem, as a consequence of Thm.\ref{thm:main}.

\vs

The way how the above operator identities (up to constants) are proved in \cite{FG09} is not via direct operator computation, but via an indirect argument. They showed that the left-hand-side operator commutes with `enough' many operators, so that it must be a scalar operator; the situation (but not the proof) is somewhat like the classical Schur's lemma. Such an indirect argument does not pin down the scalars, so we must perform direct calculation of the operator identities.

\subsection{Signed decomposition of the mutation intertwiner}

The identities in Prop.\ref{prop:the_rank_1_identity_for_intertwiners} and Prop.\ref{prop:rank_2_identities_for_intertwiners}.(2-1)--(2-2) can be proved directly, using basic analysis and a known pentagon identity for the non-compact quantum dilogarithm $\Phi^\hbar$.  Meanwhile, direct computation of the remaining identities (2-3) and (2-4) would require the hexagon and octagon identities for $\Phi^\hbar$, which have not been rigorously established in the literature yet. Hence we will take a somewhat indirect way. First, we shall prove Prop.\ref{prop:the_rank_1_identity_for_intertwiners} by direct computation, with the trivial constant $c_{A_1}=1$. For the other identities, we will assume the statements of Prop.\ref{prop:rank_2_identities_for_intertwiners}, then extract necessary operator identities of $\Phi^\hbar$, which we use in turn to compute the scalars appearing in Prop.\ref{prop:rank_2_identities_for_intertwiners}.

\vs

To elaborate this idea, recall that each factor ${\bf K}^\hbar_k$ of the operator identities of Prop.\ref{prop:rank_2_identities_for_intertwiners} is composition of the automorphism part ${\bf K}^{\hbar \sharp}_k$ and the linear part ${\bf K}'_k$. The automorphism part ${\bf K}^{\hbar \sharp}_k$ is composition of two quantum dilogarithm factors, one for ${\bf x}_k$ and one for $\til{\bf x}_k$, and the linear part ${\bf K}'_k$ is essentially a special affine shift operator for some linear transformation on $\mathbb{R}^n$. We will move around all the linear parts to group them together, so that what is left is some operator identity involving only the quantum dilogarithm factors. We will show that the quantum dilogarithm factors for ${\bf x}_k$ and those for the tilde operators $\til{\bf x}_k$ can be `separated', yielding two operator identities, each involving half number of quantum dilogarithm factors. Finally, we will show that the constant coming from one of these two separated identities is canceled by the constant for the other identity. 

\vs

We begin by dealing with the linear parts. In general, if we move around all linear parts ${\bf K}'_k$ to group them together, they might not all cancel one another. But there is a second way to decompose the intertwiner ${\bf K}^\hbar_k$ into an automorphism part and a linear part, so that if we appropriately choose out of these two possible decompositions for each mutation, the linear parts ${\bf K}'_k$ cancel one another when grouped together. This idea of considering two different decompositions has been used in the cluster algebra literature e.g. in \cite{KN} \cite{Keller}; here in the present paper, I incorporate this idea into the Fock-Goncharov representation of quantum cluster $\mathcal{D}$-variety. 

\begin{definition}[signed version of the two parts of mutation intertwiner]
\label{def:signed_two_parts}
Let $\mu_k$ stand for a mutation of a $\mathcal{D}$-seed $\Gamma$ into $\Gamma'$, i.e. $\Gamma' = \mu_k(\Gamma)$. Let $\epsilon\in \{+,-\}$ be a sign, and $\hbar\in\mathbb{R}_{>0}$. Let ${\bf x}^\hbar_{\Gamma;k}$ and let $\til{\bf x}^\hbar_{\Gamma;k}$ be the self-adjoint operators on the Hilbert space $\mathscr{H}_\Gamma$ in eq.\eqref{eq:H_Gamma}, as defined in Def.\ref{def:old_representation}. Define the signed version of the automorphism part operator as
\begin{align}
\label{eq:sharp_K_definition_epsilon}
\mathbf{K}^{\sharp \hbar (\epsilon)}_{\Gamma\mut{k}\Gamma'} := \Phi^{\hbar_k}(\epsilon \, \mathbf{x}^\hbar_{\Gamma;k})^{\epsilon} \left( \Phi^{\hbar_k}(\epsilon\, \til{\mathbf{x}}^\hbar_{\Gamma;k}) \right)^{-\epsilon}  ~ : ~ \mathscr{H}_\Gamma \to \mathscr{H}_\Gamma,
\end{align}
where $\hbar_k = \hbar/d_k$, and $\Phi^{\hbar_k}$ is the non-compact quantum dilogarithm function defined in \S\ref{subsec:non-compact_QD}. Define the signed version of the linear part operator as
\begin{align}
\label{eq:bf_K_prime_epsilon}
{\mathbf{K}'}^{(\epsilon)}_{\Gamma\mut{k}\Gamma'} := \mathbf{S}_{\Gamma;\mathbf{c}^{(\epsilon)}_{\Gamma\mut{k}\Gamma'}} \circ \mathbf{I}_{\Gamma,\Gamma'}: \mathscr{H}_{\Gamma'} \to \mathscr{H}_\Gamma,
\end{align}
where ${\bf I}_{\Gamma,\Gamma'}: \mathscr{H}_{\Gamma'} \to \mathscr{H}_\Gamma$ is as in eq.\eqref{eq:bf_I}, and $\mathbf{S}_{\Gamma;\mathbf{c}^{(\epsilon)}_{\Gamma\mut{k}\Gamma'}} : \mathscr{H}_\Gamma \to \mathscr{H}_\Gamma$ is the special affine shift operator defined as in eq.\eqref{eq:S_bf_c} and eq.\eqref{eq:def_of_S_c_t} for the matrix
$$
\mathbf{c}^{(\epsilon)}_{\Gamma\mut{k}\Gamma'} = (c^{(\epsilon)}_{ij})_{i,j\in \{1,\ldots,n\}} \in \mathrm{SL}_\pm(n,\mathbb{R})
$$
given by
\begin{align}
\label{eq:c_Gamma_Gamma_prime_entries_epsilon}
c^{(\epsilon)}_{ii} = 1, \quad \forall i \neq k, \qquad c^{(\epsilon)}_{kk} = -1, \qquad c^{(\epsilon)}_{ik} = [-\epsilon\, \varepsilon_{ki}]_+, \quad \forall i\neq k, \qquad c^{(\epsilon)}_{ij}=0 \quad\mbox{otherwise}.
\end{align}
\end{definition}

\vs

So, Fock-Goncharov's two parts $\mathbf{K}^{\sharp \hbar}_{\Gamma\mut{k}\Gamma'}$ and ${\mathbf{K}'}_{\Gamma\mut{k}\Gamma'}$ previously defined in Def.\ref{def:bf_K} coincide just with the $(+)$-versions $\mathbf{K}^{\sharp \hbar (+)}_{\Gamma\mut{k}\Gamma'}$ and ${\mathbf{K}'}^{(+)}_{\Gamma\mut{k}\Gamma'}$, respectively. I show that each of the two signed versions provides a valid decomposition of the unitary intertwiner.

\begin{lemma}[equality of the two decompositions]
\label{lem:equality_of_the_two_decompositions}
One has
$$
{\bf K}^\hbar_{\Gamma\mut{k} \Gamma'} = \mathbf{K}^{\sharp \hbar (+)}_{\Gamma\mut{k}\Gamma'} \, {\mathbf{K}'}^{(+)}_{\Gamma\mut{k}\Gamma'} = \mathbf{K}^{\sharp \hbar (-)}_{\Gamma\mut{k}\Gamma'} \, {\mathbf{K}'}^{(-)}_{\Gamma\mut{k}\Gamma'},
$$
as equalities of unitary operators $\mathscr{H}_{\Gamma'} \to \mathscr{H}_\Gamma$.
\end{lemma}

{\it Proof.} We just have to show the latter equation $\mathbf{K}^{\sharp \hbar (+)}_{\Gamma\mut{k}\Gamma'} \, {\mathbf{K}'}^{(+)}_{\Gamma\mut{k}\Gamma'} = \mathbf{K}^{\sharp \hbar (-)}_{\Gamma\mut{k}\Gamma'} \, {\mathbf{K}'}^{(-)}_{\Gamma\mut{k}\Gamma'}$. Putting in the definitions and canceling ${\bf I}_{\Gamma,\Gamma'}$ from both sides, the equation boils down to
\begin{align*}
\Phi^{\hbar_k}({\bf x}^\hbar_k)  \,\Phi^{\hbar_k}(\til{\bf x}^\hbar_k) ^{-1} \, {\bf S}_{{\bf c}^{(+)}}
= \Phi^{\hbar_k}(-{\bf x}^\hbar_k)^{-1}  \,\Phi^{\hbar_k}(-\til{\bf x}^\hbar_k)  \, {\bf S}_{{\bf c}^{(-)}},
\end{align*}
with every factor being a unitary operator $\mathscr{H}_\Gamma\to \mathscr{H}_\Gamma$, where ${\bf c}^{(\epsilon)} = (c^{(\epsilon)}_{ij})_{i,j}$ is defined in eq.\eqref{eq:c_Gamma_Gamma_prime_entries_epsilon}. Moving some factors around, it suffices to check
\begin{align}
\label{eq:equality_of_two_decompositions_boiled_down1}
\Phi^{\hbar_k}(-\til{\bf x}^\hbar_k)^{-1} \, \Phi^{\hbar_k}(-{\bf x}^\hbar_k) \, \Phi^{\hbar_k}({\bf x}^\hbar_k)  \,\Phi^{\hbar_k}(\til{\bf x}^\hbar_k) ^{-1} 
= {\bf S}_{{\bf c}^{(-)}} \, {\bf S}_{{\bf c}^{(+)}}^{-1}.
\end{align}
The two self-adjoint operators ${\bf x}^\hbar_k$ and $\til{\bf x}_k^\hbar$ satisfy the Weyl relations for the Heisenberg relations $[{\bf x}^\hbar_k, \til{\bf x}^\hbar_k]=0$; one way of saying this situation is that these two self-adjoint operators {\em strongly commute} with each other. Hence the left-hand-side of eq.\eqref{eq:equality_of_two_decompositions_boiled_down1} can be interpreted as a result of the two-variable functional calculus for the strongly commuting self-adjoint operators ${\bf x}^\hbar_k$ and $\til{\bf x}^\hbar_k$; see e.g. \cite[\S5.5]{Schmudgen}. Therefore one can switch factors and simplify using functional relations, manipulating as if ${\bf x}^\hbar_k$ and $\til{\bf x}^\hbar_k$ are real numbers:
\begin{align*}
\Phi^{\hbar_k}(-\til{\bf x}^\hbar_k)^{-1} \, \Phi^{\hbar_k}(-{\bf x}^\hbar_k) \, \Phi^{\hbar_k}({\bf x}^\hbar_k)  \,\Phi^{\hbar_k}(\til{\bf x}^\hbar_k) ^{-1} 
& = \Phi^{\hbar_k}(-{\bf x}^\hbar_k) \, \Phi^{\hbar_k}({\bf x}^\hbar_k)  \, \Phi^{\hbar_k}(-\til{\bf x}^\hbar_k)^{-1} \, \Phi^{\hbar_k}(\til{\bf x}^\hbar_k) ^{-1}  \\
& = c_{\hbar_k} \, e^{ ({\bf x}^\hbar_k)^2/(4\pi {\rm i} \hbar_k)} \, \left( c_{\hbar_k} \, e^{ (\til{\bf x}^\hbar_k)^2/(4\pi {\rm i} \hbar_k)} \right)^{-1} \qquad (\because \mbox{eq.\eqref{eq:QD_identity_quadratic}}) \\
& = \exp\left( \frac{ ({\bf x}^\hbar_k)^2 - (\til{\bf x}^\hbar_k)^2}{4\pi{\rm i} \hbar_k} \right) \\
& = \exp \left( \frac{({\bf x}^\hbar_k - \til{\bf x}^\hbar_k)({\bf x}^\hbar_k + \til{\bf x}^\hbar_k)}{4\pi{\rm i} \hbar_k} \right).
\end{align*}
The last expression can be viewed as a result of the two-variable functional calculus for the strongly commuting self-adjoint operators ${\bf x}^\hbar_k - \til{\bf x}^\hbar_k$ and ${\bf x}^\hbar_k + \til{\bf x}^\hbar_k$, which are $- 2 \sum_{j=1}^n \varepsilon_{kj} {\bf q}^\hbar_j$ and $d_k^{-1} {\bf p}^\hbar_k$ respectively, in view of equations \eqref{eq:old_representation} and \eqref{eq:old_representation_tilde}. Therefore, the left-hand-side of eq.\eqref{eq:equality_of_two_decompositions_boiled_down1} equals $\exp\left( (-\sum_{j=1}^n \varepsilon_{kj} {\bf q}^\hbar_j) \frac{ {\bf p}^\hbar_k}{2\pi {\rm i} \hbar} \right)$. Notice ${\bf q}^\hbar_j=a_j$ and $\frac{ {\bf p}^\hbar_k}{2\pi {\rm i} \hbar} = \frac{\partial}{\partial a_k}$, so, like in \eqref{eq:bf_S_as_shift_operator}--\eqref{eq:bf_S_as_Taylor_series}, one can expect that this left-hand-side operator acts as
\begin{align}
\label{eq:equality_of_two_decompositions_boiled_down2}
\textstyle \left( \exp\left( (-\sum_{j=1}^n \varepsilon_{kj} {\bf q}^\hbar_j) \frac{ {\bf p}^\hbar_k}{2\pi {\rm i} \hbar} \right) f \right)(a_1,\ldots,a_k,\ldots,a_n) = f(a_1,\ldots,a_k-\sum_{j\neq k} \varepsilon_{kj} a_j,\ldots,a_n)
\end{align}
for $f\in \mathscr{H}_\Gamma = L^2(\mathbb{R}^n, da_1\cdots da_n)$; I used $\varepsilon_{kk}=0$. Now I give an argument for how to prove eq.\eqref{eq:equality_of_two_decompositions_boiled_down2} rigorously. First, using eq.\eqref{eq:conjugation_of_bf_S_on_bf_p}--\eqref{eq:conjugation_of_bf_S_on_bf_q}, one can easily show by a simple linear algebra that there exists a matrix ${\bf c}$ so that ${\bf S}_{{\bf c}} \, (-\sum_{j=1}^n \varepsilon_{kj} {\bf q}^\hbar_j) \, {\bf S}_{{\bf c}}^{-1} = {\bf q}^\hbar_1$ and ${\bf S}_{{\bf c}} \, {\bf p}^\hbar_k \, {\bf S}_{{\bf c}}^{-1} = c \, {\bf p}^\hbar_2$ for some $c\neq 0$; namely, ${\bf c}$ can be any matrix in ${\rm SL}_\pm (n,\mathbb{R})$ such that the inverse ${\bf c}^{-1} = (c^{ij})_{i,j}$ satisfies $c^{j1} = -\varepsilon_{kj}$ for all $j$ and $c^{k2}\neq 0$. One can observe that now it suffices to show
\begin{align}
\label{eq:equality_of_two_decompositions_boiled_down3}
\textstyle
(\exp( {\bf q}^\hbar_1 \, \frac{ c \, {\bf p}^\hbar_2 }{2\pi {\rm i} \hbar} ) f)(a_1,a_2,\ldots,a_n) = f(a_1,a_2+c\,a_1,a_3,\ldots,a_n),
\end{align}
which is essentially a problem about $L^2(\mathbb{R}^2, da_1 \, da_2)$. To prove this, consider the Fourier transform $\mathcal{F}_i$ in the $i$-th variable $a_i$, defined as
\begin{align}
\label{eq:F_i}
(\mathcal{F}_i f)(a_1,\ldots,a_i,\ldots,a_n) = \int_\mathbb{R} e^{-2\pi{\rm i} a_i b_i} f(a_1,\ldots,b_i,\ldots,a_n) \, db_i
\end{align}
for $f\in L^1\cap L^2$. It is well known that $\mathcal{F}_i$ is a unitary operator on $\mathscr{H}_\Gamma$, and satisfies
\begin{align}
\label{eq:F_i_conjugation}
\textstyle \mathcal{F}_i {\bf p}^\hbar_i \mathcal{F}_i^{-1} = - (2\pi)^2 \hbar \, {\bf q}^\hbar_i, \qquad
\mathcal{F}_i \,  {\bf q}^\hbar_i \mathcal{F}_i^{-1} =\frac{ {\bf p}^\hbar_i }{(2\pi)^2 \hbar},  \qquad
\mathcal{F}_i {\bf p}^\hbar_j \mathcal{F}_i^{-1}={\bf p}^\hbar_j, \qquad
\mathcal{F}_i \, {\bf q}^\hbar_j \mathcal{F}_i^{-1} = {\bf q}^\hbar_j
\end{align}
hold for all $i,j$ with $i\neq j$. So $\mathcal{F}_2 \, \exp( {\bf q}^\hbar_1 \, \frac{ c \, {\bf p}^\hbar_2 }{2\pi {\rm i} \hbar} ) \, \mathcal{F}_2^{-1} = \exp(2\pi {\rm i} c \, {\bf q}^\hbar_1 {\bf q}^\hbar_2)$, which acts as multiplication by $e^{2\pi {\rm i} c a_1 a_2}$. From here, it is then a straightforward exercise about Fourier transform to prove eq.\eqref{eq:equality_of_two_decompositions_boiled_down3}.

\vs

So we justified eq.\eqref{eq:equality_of_two_decompositions_boiled_down2}, which is in fact saying that the operator $\exp\left( (-\sum_{j=1}^n \varepsilon_{kj} {\bf q}^\hbar_j) \frac{ {\bf p}^\hbar_k}{2\pi {\rm i} \hbar} \right)$, hence so is the left-hand-side of eq.\eqref{eq:equality_of_two_decompositions_boiled_down1}, is a special affine shift operator:
$$
\textstyle 
\Phi^{\hbar_k}(-\til{\bf x}^\hbar_k)^{-1} \, \Phi^{\hbar_k}(-{\bf x}^\hbar_k) \, \Phi^{\hbar_k}({\bf x}^\hbar_k)  \,\Phi^{\hbar_k}(\til{\bf x}^\hbar_k) ^{-1} = \exp\left( (-\sum_{j=1}^n \varepsilon_{kj} {\bf q}^\hbar_j) \frac{ {\bf p}^\hbar_k}{2\pi {\rm i} \hbar} \right) = {\bf S}_{\wh{\bf c}},
$$
where
$$
\wh{\bf c}=(\wh{c}_{ij})_{i,j},\qquad
\wh{c}_{ii} = 1,\quad \forall i=1,\ldots,n, \quad
\wh{c}_{ik} = -\varepsilon_{ki}, \quad \forall i\neq k, \quad
\wh{c}_{ij}=0 \quad \mbox{otherwise}.
$$
Now, in view of the multiplicativity of the special affine shift operators (Lem.\ref{lem:special_affine_shift_operators_multiplicative}), what remains to show in order to prove eq.\eqref{eq:equality_of_two_decompositions_boiled_down1} is just the matrix identity
$$
\wh{\bf c} = {\bf c}^{(-)} \, ({\bf c}^{(+)})^{-1}.
$$
Let's show $\wh{\bf c} \, {\bf c}^{(+)} = {\bf c}^{(-)}$.  It is easy to check that each side coincides with the identity matrix except at the $k$-th column. The $k$-th entry of the $k$-th column is $-1$ for each side. The $i$-th entry of the $k$-th column, for $i\neq k$, is $\varepsilon_{ki} + [-\varepsilon_{ki}]_+$ for the left-hand-side $\wh{\bf c} \, {\bf c}^{(+)}$, and is $[\varepsilon_{ki}]_+$ for the right-hand-side ${\bf c}^{(-)}$. Since $a+[-a]_+ = [a]_+$ holds for any real number $a$, we get the desired equality. \qed

\vs

The above Lem.\ref{lem:equality_of_the_two_decompositions} is the intertwiner operator version of Prop.3.1 and Prop.4.2 of \cite{KN}; for our concern, an emphasis must be on the fact that the equality $\mathbf{K}^{\sharp \hbar (+)}_{\Gamma\mut{k}\Gamma'} \, {\mathbf{K}'}^{(+)}_{\Gamma\mut{k}\Gamma'} = \mathbf{K}^{\sharp \hbar (-)}_{\Gamma\mut{k}\Gamma'} \, {\mathbf{K}'}^{(-)}_{\Gamma\mut{k}\Gamma'}$ holds on the nose, not just up to constant. In the meantime, as mentioned in \cite{KN}, the difference between the two possible choices of signed decompositions is not just a technical matter, but represents a deeper phenomenon studied e.g. in \cite{Keller}.

\vs

For later use, I record the intertwining equations satisfied by the linear part.
\begin{lemma}[some intertwining equations for the linear part; {\cite[Thm.5.6]{FG09}}]
\label{lem:some_intertwining_equations_for_the_linear_part}
Let $\Gamma,\Gamma=\mu_k(\Gamma')$, and ${{\bf K}'}^{(\epsilon)}_{\Gamma\mut{k}\Gamma'} : \mathscr{H}_{\Gamma'} \to \mathscr{H}_\Gamma$ be as in Def.\ref{def:signed_two_parts}, for a sign $\epsilon\in\{+,-\}$. Then the following equalities of self-adjoint operators hold:
\begin{align}
\label{eq:bf_K_prime_conjugation_on_bf_x_i}
{\mathbf{K}'}^{(\epsilon)}_{\Gamma\mut{k}\Gamma'} \,\, {\bf x}^\hbar_{\Gamma';i} \,\, ({\mathbf{K}'}^{(\epsilon)}_{\Gamma\mut{k}\Gamma'})^{-1} & = 
\left\{
{\renewcommand{\arraystretch}{1.2} \begin{array}{ll}
{\bf x}^\hbar_{\Gamma;i} + [\epsilon\, \varepsilon_{ik}]_+ \, {\bf x}^\hbar_{\Gamma;k}, & \mbox{if $i\neq k$}, \\
-{\bf x}^\hbar_{\Gamma;k} & \mbox{if $i=k$,}
\end{array}} \right.
\\
\label{eq:bf_K_prime_conjugation_on_til_bf_x_i}
{\mathbf{K}'}^{(\epsilon)}_{\Gamma\mut{k}\Gamma'} \,\, \til{\bf x}^\hbar_{\Gamma';i} \,\, ({\mathbf{K}'}^{(\epsilon)}_{\Gamma\mut{k}\Gamma'})^{-1} & = 
\left\{
{\renewcommand{\arraystretch}{1.2} \begin{array}{ll}
\til{\bf x}^\hbar_{\Gamma;i} + [\epsilon\, \varepsilon_{ik}]_+ \, \til{\bf x}^\hbar_{\Gamma;k}, & \mbox{if $i\neq k$}, \\
-\til{\bf x}^\hbar_{\Gamma;k} & \mbox{if $i=k$.}
\end{array}} \right.
\end{align}
\end{lemma}

{\it Proof.} Note that the matrix ${\bf c}^{(\epsilon)}_{\Gamma\mut{k}\Gamma'}$ defined in eq.\eqref{eq:c_Gamma_Gamma_prime_entries_epsilon} equals its own inverse. Hence from Lem.\ref{lem:bf_S_conjugation_on_bf_p_and_bf_q} one obtains:
\begin{align*}
& {\mathbf{K}'}^{(\epsilon)}_{\Gamma\mut{k}\Gamma'} \, \mathbf{q}^\hbar_{\Gamma';i} \, ({\mathbf{K}'}^{(\epsilon)}_{\Gamma\mut{k}\Gamma'})^{-1}
= (\mathbf{S}_{\Gamma;\mathbf{c}^{(\epsilon)}_{\Gamma\mut{k}\Gamma'}} \circ \mathbf{I}_{\Gamma,\Gamma'}) \left( \mathbf{q}^\hbar_{\Gamma';i} \right) (\mathbf{I}_{\Gamma,\Gamma'}^{-1} \circ \mathbf{S}_{\Gamma;\mathbf{c}^{(\epsilon)}_{\Gamma\mut{k}\Gamma'}}^{-1} ) \\
& = \mathbf{S}_{\Gamma;\mathbf{c}^{(\epsilon)}_{\Gamma\mut{k}\Gamma'}} \left( \mathbf{q}^\hbar_{\Gamma;i} \right) \mathbf{S}_{\Gamma;\mathbf{c}^{(\epsilon)}_{\Gamma\mut{k}\Gamma'}}^{-1} = \left\{
  \begin{array}{ll}
    \wh{\mathbf{q}}^\hbar_{\Gamma;i}, & i\neq k, \\
    - \wh{\mathbf{q}}^\hbar_{\Gamma;k} + \sum_{j=1}^n [-\epsilon\, \varepsilon_{kj}]_+ \, \wh{\mathbf{q}}^\hbar_{\Gamma;j}, & i=k,
  \end{array}
\right.
\end{align*}
and similarly
\begin{align*}
& {\mathbf{K}'}^{(\epsilon)}_{\Gamma\mut{k}\Gamma'} \, \mathbf{p}^\hbar_{\Gamma';i} \, ({\mathbf{K}'}^{(\epsilon)}_{\Gamma\mut{k}\Gamma'})^{-1}
= \mathbf{S}_{\Gamma;\mathbf{c}^{(\epsilon)}_{\Gamma\mut{k}\Gamma'}} \left( \mathbf{p}^\hbar_{\Gamma;i} \right) \mathbf{S}_{\Gamma;\mathbf{c}^{(\epsilon)}_{\Gamma\mut{k}\Gamma'}}^{-1} = \left\{
\begin{array}{ll}
 {\bf p}^\hbar_{\Gamma;i} + [-\epsilon \, \varepsilon_{ki}]_+ {\bf p}^\hbar_{\Gamma;k},  & i\neq k, \\
-{\bf p}^\hbar_{\Gamma;k}, & i=k,
\end{array}
\right.
\end{align*}
Observe now
\begin{align*}
\textstyle {\mathbf{K}'}^{(\epsilon)}_{\Gamma\mut{k}\Gamma'} \,\, \left( \frac{d_k^{-1}}{2} \mathbf{p}^\hbar_{\Gamma';k} \mp \sum_{j=1}^n \varepsilon'_{kj} \mathbf{q}^\hbar_{\Gamma';j}\right) ({\mathbf{K}'}^{(\epsilon)}_{\Gamma\mut{k}\Gamma'})^{-1}  
\stackrel{\vee}{=}\frac{d_k^{-1}}{2} ( - \mathbf{p}^\hbar_{\Gamma;k}) \mp \sum_{j\neq k} (-\varepsilon_{kj}) \mathbf{q}^\hbar_{\Gamma;j} 
= - \left( \frac{d_k^{-1}}{2} \mathbf{p}^\hbar_{\Gamma;k} \mp \sum_{j=1}^n \varepsilon_{kj} \mathbf{q}^\hbar_{\Gamma;j} \right),
\end{align*}
where for the checked equality $\vee$, I used the fact that $\varepsilon'_{kj} = - \varepsilon_{kj}$ and $\varepsilon_{kk}=0$. So we just proved eq.\eqref{eq:bf_K_prime_conjugation_on_bf_x_i}--\eqref{eq:bf_K_prime_conjugation_on_til_bf_x_i} for the case $i=k$. On the other hand, for each $i\neq k$  observe
\begin{align}
\label{eq:bf_K_prime_conjugation_eq1}
{\renewcommand{\arraystretch}{1.7}  \begin{array}{l}
{\mathbf{K}'}^{(\epsilon)}_{\Gamma\mut{k}\Gamma'} \,\, \left( \frac{d_i^{-1}}{2} \mathbf{p}^\hbar_{\Gamma';i} \mp \sum_{j=1}^n \varepsilon'_{ij} \mathbf{q}^\hbar_{\Gamma';j} \right) \,\, ({\mathbf{K}'}^{(\epsilon)}_{\Gamma\mut{k}\Gamma'})^{-1}  \\
= \frac{d_i^{-1}}{2}\left(\mathbf{p}^\hbar_{\Gamma;i} + [- \epsilon\, \varepsilon_{ki}]_+  \, \mathbf{p}^\hbar_{\Gamma;k}\right) \mp  \sum_{j\neq k} \varepsilon'_{ij} \, \mathbf{q}^\hbar_{\Gamma;j}
\mp \varepsilon'_{ik} \left( - \mathbf{q}^\hbar_{\Gamma;k} + \sum_{\ell\neq k} [-\epsilon\, \varepsilon_{k\ell}]_+  \,\mathbf{q}^\hbar_{\Gamma;\ell} \right).
\end{array}}
\end{align}
Note from eq.\eqref{eq:wh_varepsilon} that $[-\epsilon\, \varepsilon_{ki}]_+ \, d_i^{-1} = [-\epsilon\, \wh{\varepsilon}_{ki} ]_+ = [\epsilon\, \wh{\varepsilon}_{ik} ]_+ = [\epsilon \, \varepsilon_{ik}]_+ \, d_k^{-1}$, and from the matrix mutation formula in eq.\eqref{eq:varepsilon_prime_formula} that $\varepsilon'_{ij} = \varepsilon_{ij} + \frac{1}{2}(|\varepsilon_{ik}|\varepsilon_{kj} + \varepsilon_{ik} |\varepsilon_{kj}|)$ for $j\neq k$, and $\varepsilon'_{ik}=-\varepsilon_{ik}$. Hence, by renaming $\ell$ by $j$, one gets
\begin{align*}
\textstyle
{\rm eq}.\eqref{eq:bf_K_prime_conjugation_eq1}
= \left( \frac{d_i^{-1}}{2}\mathbf{p}^\hbar_{\Gamma;i} \mp \sum_{j=1}^n  \varepsilon_{ij} \, \mathbf{q}^\hbar_{\Gamma;j}\right)
+ \left(\frac{d_k^{-1}}{2} [\epsilon\, \varepsilon_{ik}]_+ \, \mathbf{p}^\hbar_{\Gamma;k} \mp \sum_{j\neq k}  (\underbrace{\textstyle \frac{1}{2}(|\varepsilon_{ik}|\varepsilon_{kj} + \varepsilon_{ik}|\varepsilon_{kj}|) - \varepsilon_{ik} [-\epsilon\, \varepsilon_{kj}]_+ }) \,\mathbf{q}^\hbar_{\Gamma;j} \right)
\end{align*}
Putting $|\varepsilon_{ik}| = |\epsilon\, \varepsilon_{ik}| = 2 [\epsilon\, \varepsilon_{ik}]_+ - \epsilon\, \varepsilon_{ik}$ and $|\varepsilon_{kj}| = |\epsilon \,\varepsilon_{kj}| = 2[- \epsilon\,\varepsilon_{kj}]_+ + \epsilon\, \varepsilon_{kj}$ makes the underbraced part to
$$
\textstyle \frac{1}{2}\left(  (2[\epsilon\, \varepsilon_{ik}]_+-\epsilon\,\varepsilon_{ik})\varepsilon_{kj} + \varepsilon_{ik}(2[-\epsilon\, \varepsilon_{kj}]_+ + \epsilon\, \varepsilon_{kj}) \right) - \varepsilon_{ik} [- \epsilon\, \varepsilon_{kj}]_+ = [\epsilon\, \varepsilon_{ik}]_+ \,\varepsilon_{kj}.
$$
Thus, keeping in mind $\varepsilon_{kk}=0$, we get
$$
\textstyle {\rm eq}.\eqref{eq:bf_K_prime_conjugation_eq1} = \left(\frac{d_i^{-1}}{2} \mathbf{p}^\hbar_{\Gamma;i} \mp \sum_{j=1}^n  \varepsilon_{ij} \mathbf{q}^\hbar_{\Gamma;j}\right)
+ [\epsilon\, \varepsilon_{ik}]_+ \left( \frac{d_k^{-1}}{2} \mathbf{p}^\hbar_{\Gamma;k} \mp \sum_{j=1}^n  \varepsilon_{kj} \mathbf{q}^\hbar_{\Gamma;j}\right),
$$
So we proved eq.\eqref{eq:bf_K_prime_conjugation_on_bf_x_i}--\eqref{eq:bf_K_prime_conjugation_on_til_bf_x_i} for the case $i\neq k$. \qed

\vs

However straightforward this proof is, it is worthwhile to write it down, because of the following situation: this lemma for $\epsilon=+$ is claimed in \cite[Thm5.6]{FG09} without proof, but is wrong as is written there; the statement does not hold for the operators ${\bf x}^\hbar_i$, $\til{\bf x}^\hbar_i$, and ${\bf K}'_{\Gamma\mut{k}\Gamma'}$ as defined in \cite{FG09}. What led to this trouble is that the operators to represent ${\bf x}^\hbar_i$, ${\bf b}^\hbar_i$, $\til{\bf x}^\hbar_i$, $\til{\bf b}^\hbar_i$ were changed to new ones in \cite{FG09}, compared to previous papers of Fock and Goncharov, but ${\bf K}'_{\Gamma\mut{k}\Gamma'}$ was not changed accordingly. Instead of trying to come up with a correct linear part ${\bf K}'_{\Gamma\mut{k}\Gamma'}$, I suggest we better use Fock Goncharov's old representation for ${\bf x}^\hbar_i$, ${\bf b}^\hbar_i$, $\til{\bf x}^\hbar_i$, $\til{\bf b}^\hbar_i$ which are used e.g. in \cite{FG07} and mentioned also in eq.(78) of \cite{FG09}. In my opinion, using Fock-Goncharov's `new' operators in \cite{FG09}, together with their viewpoint of considering the variables $a_1,\ldots,a_n$ in the Hilbert space $\mathscr{H}_\Gamma = L^2(\mathbb{R}^n,da_1 \cdots da_n)$ as ${\rm log}$ version of the cluster $\mathcal{A}$-variables, is not natural; see e.g. \cite{K16} for some more details. The canonicity of the `old' operators, which are used in the present paper, will be discussed in more depth in a forthcoming paper \cite{KS}.

\vs

\subsection{Collecting the linear parts}
\label{subsec:collectin_the_linear_parts}

Suppose we have a trivial cluster transformation ${\bf m}$ of $\mathcal{D}$-seeds; namely,
$$
{\bf m} = {\bf m}_r \circ {\bf m}_{r-1} \circ \cdots \circ {\bf m}_1
$$
where each ${\bf m}_i$ is a mutation or a seed automorphism, and ${\bf m}$ connects a $\mathcal{D}$-seed $\Gamma$ to itself. i.e. it starts and ends at the same seed. Let
\begin{align}
\label{eq:seed_Gamma_i}
\Gamma^{(0)} = \Gamma, \qquad
\Gamma^{(i)} = {\bf m}_i(\Gamma^{(i-1)}), \quad \forall i =1,2,\ldots,r.
\end{align}
Then $\Gamma^{(r)} = \Gamma^{(0)}$. The intertwiner corresponding to this sequence ${\bf m}$ is
$$
\pi^\hbar({\bf m}) = \pi^\hbar({\bf m}_1) \circ \cdots \circ \pi^\hbar({\bf m}_r) ~ : ~ \mathscr{H}_{\Gamma^{(0)}} \to \mathscr{H}_{\Gamma^{(0)}}.
$$
Suppose ${\bf m}$ is a trivial morphism of the saturated cluster modular groupoid $\wh{\mathcal{G}}^\mathcal{D}_{|\Gamma|}$, so that by Prop.\ref{thm:FG_constant}, the composed operator $\pi^\hbar({\bf m})$ is a scalar operator, namely $\pi^\hbar({\bf m}) = c_{\bf m} \cdot {\rm Id}$. Our task is to compute what this scalar $c_{\bf m}\in {\rm U}(1)$ is. In the present subsection, we study how to manipulate $\pi^\hbar({\bf m})$ to make it easier to investigate.

\vs

Let
$$
\vec{\epsilon} = (\epsilon_1 \, \epsilon_2 \, \cdots \, \epsilon_r)
$$
be a sequence of signs $\epsilon_i \in \{+,-,0\}$, with the condition that
$$
\mbox{$\epsilon_i=0$ if and only if ${\bf m}_i$ is a seed automorphism}.
$$
Note that each ${\bf m}_i$ connects $\Gamma^{(i-1)}$ to $\Gamma^{(i)}$. In case when ${\bf m}_i$ is a mutation, say $\mu_{k_i}$ for some $k_i \in \{1,\ldots,n\}$, then write the intertwiner $\pi^\hbar({\bf m}_i) : \mathscr{H}_{\Gamma^{(i)}} \to \mathscr{H}_{\Gamma^{(i-1)}}$ for ${\bf m}_i$ in terms of the signed decomposition with sign $\epsilon_i$:
\begin{align*}
\pi^\hbar({\bf m}_i)
& = {\bf K}^\hbar_{{\bf m}_i}
= {\bf K}^\hbar_{\Gamma^{(i-1)}\mut{k_i}\Gamma^{(i)}}
= {\bf K}^{\hbar\sharp(\epsilon_i)}_{\Gamma^{(i-1)}\mut{k_i}\Gamma^{(i)}}
{{\bf K}'}^{(\epsilon_i)}_{\Gamma^{(i-1)}\mut{k_i}\Gamma^{(i)}} \\
& = \Phi^{\hbar_{k_i}}(\epsilon_i \, {\bf x}^\hbar_{\Gamma^{(i-1)};k_i})^{\epsilon_i} \, \Phi^{\hbar_{k_i}}(\epsilon_i \, \til{\bf x}^\hbar_{\Gamma^{(i-1)};k_i})^{-\epsilon_i}
\, {\bf S}_{\Gamma^{(i-1)};{\bf c}^{(\epsilon_i)}_{{\bf m}_i}} \, {\bf I}_{\Gamma^{(i-1)},\Gamma^{(i)}} \quad (\mbox{when ${\bf m}_i$ is the mutation $\mu_{k_i}$}),
\end{align*}
where
\begin{align}
\label{eq:bf_c_m_i_epsilon_i}
{\bf c}^{(\epsilon_i)}_{{\bf m}_i} = {\bf c}^{(\epsilon_i)}_{\Gamma^{(i-1)}\mut{k_i}\Gamma^{(i)}} \qquad (\mbox{when ${\bf m}_i$ is the mutation $\mu_{k_i}$}),
\end{align}
is as defined in eq.\eqref{eq:c_Gamma_Gamma_prime_entries_epsilon}. So such $\pi^\hbar({\bf m}_i)$ is composition of four factors. Now, in case when ${\bf m}_i$ is a seed automorphism, say $P_{\sigma_i}$ for some permutation $\sigma_i$ of $\{1,\ldots,n\}$, then the corresponding intertwiner $\pi^\hbar({\bf m}_i) : \mathscr{H}_{\Gamma^{(i)}}\to\mathscr{H}_{\Gamma^{(i-1)}}$ is as in eq.\eqref{eq:bf_K_P_sigma}:
$$
\pi^\hbar({\bf m}_i) = {\bf P}_{\Gamma^{(i-1)};\sigma_i} \circ {\bf I}_{\Gamma^{(i-1)},\Gamma^{(i)}},
$$
which is composition of two factors. For the sake of a uniform notation, denote the permutation matrix for $\sigma_i$ as
\begin{align}
\label{eq:bf_c_sigma_i}
{\bf c}^{(\epsilon_i)}_{{\bf m}_i} = (\delta_{j,\sigma_i(l)})_{j,l} \qquad (\mbox{when ${\bf m}_i$ is the seed automorphism $\sigma_i$})
\end{align}
so that
$$
\pi^\hbar({\bf m}_i) = {\bf S}_{\Gamma^{(i-1)};{\bf c}^{(\epsilon_i)}_{{\bf m}_i}} \, {\bf I}_{\Gamma^{(i-1)},\Gamma^{(i)}} \qquad (\mbox{when ${\bf m}_i$ is the seed automorphism $\sigma_i$}).
$$

\vs

Now take $\pi^\hbar({\bf m}_1) \circ \cdots \circ \pi^\hbar({\bf m}_r)$, write each $\pi^\hbar({\bf m}_i)$ as composition of four or two factors, as done above. Move all ${\bf I}_{\Gamma^{(i-1)}, \Gamma^{(i)}}$ to the right, so that they group together as
$$
{\bf I}_{\Gamma^{(0)},\Gamma^{(1)}} \,
{\bf I}_{\Gamma^{(1)},\Gamma^{(2)}} \,
\cdots \,
{\bf I}_{\Gamma^{(r-1)},\Gamma^{(r)}}
$$
at the rightmost place. In view of eq.\eqref{eq:bf_I}, one can see that this composition of all ${\bf I}_{\Gamma^{(i-1)},\Gamma^{(i)}}$'s is the identity map $\mathscr{H}_{\Gamma^{(0)}} \to \mathscr{H}_{\Gamma^{(0)}}$, so we can now delete them. As a consequence of moving all the ${\bf I}$'s to the right as we did, the operators ${\bf x}^\hbar_{\Gamma^{(i-1)};k_i}$, $\til{\bf x}^\hbar_{\Gamma^{(i-1)};k_i}$, and ${\bf S}_{\Gamma^{(i-1)};{\bf c}^{(\epsilon_i)}_{{\bf m}_i}}$ get conjugated by bunch of ${\bf I}$'s to become now operators on $\mathscr{H}_{\Gamma^{(0)}} = \mathscr{H}_\Gamma$ to itself, instead of on $\mathscr{H}_{\Gamma^{(i)}}$. Namely, the variables $a_1^{(i)},\ldots,a_n^{(i)}$ in $\mathscr{H}_{\Gamma^{(i)}} = L^2(\mathbb{R}^n, da_1^{(i)} \cdots da_n^{(i)})$ get re-labeled as $a_1,\ldots,a_n$ in $\mathscr{H}_\Gamma = L^2(\mathbb{R}^n, da_1 \cdots da_n)$, and the operator ${\bf S}_{\Gamma^{(i-1)};{\bf c}^{(\epsilon_i)}_{{\bf m}_i}}$ becomes ${\bf S}_{\Gamma^{(0)};{\bf c}^{(\epsilon_i)}_{{\bf m}_i}}$.

\vs

Now we move all the factors ${\bf S}_{\Gamma^{(0)};{\bf c}^{(\epsilon_i)}_{{\bf m}_i}}={\bf S}_{{\bf c}^{(\epsilon_i)}_{{\bf m}_i}}$ to the right, so that they group together as
\begin{align}
\label{eq:bf_S_only}
{\bf S}_{{\bf c}^{(\epsilon_1)}_{{\bf m}_1}} \, {\bf S}_{{\bf c}^{(\epsilon_2)}_{{\bf m}_2}} \, \cdots \, {\bf S}_{{\bf c}^{(\epsilon_r)}_{{\bf m}_r}}
\end{align}
at the rightmost place; by the multiplicativity in Lem.\ref{lem:special_affine_shift_operators_multiplicative}, this part equals
$$
{\bf S}_{{\bf c}^{(\epsilon_1)}_{{\bf m}_1} \, {\bf c}^{(\epsilon_2)}_{{\bf m}_2} \, \cdots \, {\bf c}^{(\epsilon_r)}_{{\bf m}_r}},
$$
which equals the identity operator if and only if the matrix ${\bf c}^{(\epsilon_1)}_{{\bf m}_1} \, {\bf c}^{(\epsilon_2)}_{{\bf m}_2} \, \cdots \, {\bf c}^{(\epsilon_r)}_{{\bf m}_r}$ equals the identity matrix, by the injectivity in Lem.\ref{lem:special_affine_shift_operators_multiplicative}. The matrix ${\bf c}^{(\epsilon_1)}_{{\bf m}_1} \, {\bf c}^{(\epsilon_2)}_{{\bf m}_2} \, \cdots \, {\bf c}^{(\epsilon_r)}_{{\bf m}_r}$ is not in general the identity matrix, even though ${\bf m}$ is a trivial cluster transformation. Fortunately, it is known that, this is the identity for a suitably chosen sign sequence $\vec{\epsilon}$; this special sign sequence is what we will soon recall in the next subsection. 

\vs

Suppose for now that we chose such a sign sequence $\vec{\epsilon}$. Then the product ${\bf S}_{{\bf c}^{(\epsilon_1)}_{{\bf m}_1}} \, {\bf S}_{{\bf c}^{(\epsilon_2)}_{{\bf m}_2}} \, \cdots \, {\bf S}_{{\bf c}^{(\epsilon_r)}_{{\bf m}_r}}$ of special affine shift operators is identity, hence can be deleted. We are then left with only the factors written in terms of quantum dilogarithm functions. However, the self-adjoint-operator arguments $\epsilon_i {\bf x}^\hbar_{\Gamma^{(i-1)};k_i}$ and $\epsilon_i \til{\bf x}^\hbar_{\Gamma^{(i-1)};k_i}$ are now altered while we moved ${\bf S}_{{\bf c}^{\epsilon_i}_{{\bf m}_i}}$ to the right. Let us denote the resulting situation as
\begin{align}
\label{eq:quantum_dilog_factors_only}
\Phi^{\hbar_{k_1}}(\epsilon_1\, {\bf x}^\hbar_{{\bf m}_1})^{\epsilon_1}
\, \Phi^{\hbar_{k_1}}(\epsilon_1\, \til{\bf x}^\hbar_{{\bf m}_1})^{-\epsilon_1} \,
\Phi^{\hbar_{k_2}}(\epsilon_2\, {\bf x}^\hbar_{{\bf m}_2})^{\epsilon_2}
\, \Phi^{\hbar_{k_2}}(\epsilon_2\, \til{\bf x}^\hbar_{{\bf m}_2})^{-\epsilon_2} \, \cdots \,
\Phi^{\hbar_{k_r}}(\epsilon_r\, {\bf x}^\hbar_{{\bf m}_r})^{\epsilon_r}
\, \Phi^{\hbar_{k_r}}(\epsilon_r\, \til{\bf x}^\hbar_{{\bf m}_r})^{-\epsilon_r} = c_{\bf m} \cdot {\rm Id},
\end{align}
where the self-adjoint operators ${\bf x}^\hbar_{{\bf m}_i}$ and $\til{\bf x}^\hbar_{{\bf m}_i}$ on $\mathscr{H}_\Gamma$ are defined as ${\bf x}^\hbar_{{\bf m}_1} = {\bf x}^\hbar_{\Gamma^{(0)};k_1}$ and $\til{\bf x}^\hbar_{{\bf m}_1} = \til{\bf x}^\hbar_{\Gamma^{(0)};k_1}$ and
\begin{align}
\nonumber
{\bf x}^\hbar_{ {\bf m}_i } & := 
\left( {\bf S}_{{\bf c}^{(\epsilon_1)}_{{\bf m}_1}} \,  \cdots \,
{\bf S}_{{\bf c}^{(\epsilon_{i-1})}_{{\bf m}_{i-1}}} \right) \, 
{\bf x}^\hbar_{\Gamma^{(i-1)};k_i} \, \left( {\bf S}_{{\bf c}^{(\epsilon_1)}_{{\bf m}_1}} \,  \cdots \,
{\bf S}_{{\bf c}^{(\epsilon_{i-1})}_{{\bf m}_{i-1}}} \right)^{-1}, \\
\nonumber
\til{\bf x}^\hbar_{ {\bf m}_i } & := 
\left( {\bf S}_{{\bf c}^{(\epsilon_1)}_{{\bf m}_1}} \,  \cdots \,
{\bf S}_{{\bf c}^{(\epsilon_{i-1})}_{{\bf m}_{i-1}}} \right) \, 
\til{\bf x}^\hbar_{\Gamma^{(i-1)};k_i} \, \left( {\bf S}_{{\bf c}^{(\epsilon_1)}_{{\bf m}_1}} \,  \cdots \,
{\bf S}_{{\bf c}^{(\epsilon_{i-1})}_{{\bf m}_{i-1}}} \right)^{-1}.
\end{align}
for $i\ge 2$. For $i$ for which ${\bf m}_i$ is a seed automorphism, we may set $k_i$ to be any number in $\{1,\ldots,n\}$; since $\epsilon_i=0$, anyways the factor $\Phi^{\hbar_{k_i}}(\epsilon_i\, {\bf x}^\hbar_{{\bf m}_i})^{\epsilon_i}
\, \Phi^{\hbar_{k_i}}(\epsilon_i\, \til{\bf x}^\hbar_{{\bf m}_i})^{-\epsilon_i}$ equals the identity, hence can be thought of as being absent.

\vs

To summarize, for a $\mathcal{D}$-trivial cluster transformation ${\bf m}$ of length $r$, with $r'$ of them being mutations, the operator identity coming from representation of quantum cluster $\mathcal{D}$-varity can be  re-written in terms of an operator identity written with $2r'$ factors involving quantum dilogarithms, with $r'$ of them for the operators ${\bf x}^\hbar_{{\bf m}_i}$, and the remaining $r'$ of them for the tilde version operators $\til{\bf x}^\hbar_{{\bf m}_i}$. Such is how we will investigate the operator identities for trivial cluster transformations.

\vs

 For our purposes, one may assume that ${\bf m}_1,\ldots,{\bf m}_{r-1}$ are all mutations. Taking into consider the conjugation by ${\bf I}$'s,  more precise definition of the argument operators, for $i\ge 2$, is:
\begin{align}
\label{eq:bf_x_m_i_1}
{\bf x}^\hbar_{{\bf m}_i} & := 
{{\bf K}'}^{(\epsilon_1)}_{\Gamma^{(0)}\mut{k_1} \Gamma^{(1)} } \, \cdots \,
{{\bf K}'}^{(\epsilon_{i-1})}_{\Gamma^{(i-2)}\mut{k_{i-1}} \Gamma^{(i-1)} } \,
{\bf x}^\hbar_{\Gamma^{(i-1)};k_i} \,
({{\bf K}'}^{(\epsilon_{i-1})}_{\Gamma^{(i-2)}\mut{k_{i-1}} \Gamma^{(i-1)} })^{-1} \, \cdots \, 
({{\bf K}'}^{(\epsilon_1)}_{\Gamma^{(0)}\mut{k_1} \Gamma^{(1)} })^{-1} \\
\label{eq:til_bf_x_m_i_1}
\til{\bf x}^\hbar_{{\bf m}_i} & := 
{{\bf K}'}^{(\epsilon_1)}_{\Gamma^{(0)}\mut{k_1} \Gamma^{(1)} } \, \cdots \,
{{\bf K}'}^{(\epsilon_{i-1})}_{\Gamma^{(i-2)}\mut{k_{i-1}} \Gamma^{(i-1)} } \,
\til{\bf x}^\hbar_{\Gamma^{(i-1)};k_i} \,
({{\bf K}'}^{(\epsilon_{i-1})}_{\Gamma^{(i-2)}\mut{k_{i-1}} \Gamma^{(i-1)} })^{-1} \, \cdots \, 
({{\bf K}'}^{(\epsilon_1)}_{\Gamma^{(0)}\mut{k_1} \Gamma^{(1)} })^{-1}
\end{align}

\subsection{Tropical signs} 

Now I recall from \cite{KN} how to choose the sign sequence $\vec{\epsilon}$, so that the product ${\bf c}^{(\epsilon_1)}_{{\bf m}_1} \, {\bf c}^{(\epsilon_2)}_{{\bf m}_2} \, \cdots \, {\bf c}^{(\epsilon_r)}_{{\bf m}_r}$ is the identity matrix. First, recall that a {\em semifield} refers to a set $\mathbb{P}$ equipped with two binary operations, namely addition $+$ and multiplication $\cdot$, both being associative and commutative, so that $(\mathbb{P}, \cdot)$ is an abelian group, and the distributive law $(a+ b)\cdot c = (a\cdot c) + (b\cdot c)$ holds. I first recall two basic well-known semifields to use.
\begin{lemma}[see e.g. {\cite[\S2]{KN}}]
\begin{enumerate}
\item[\rm (1)] (the universal semifield of $n$ variables) Consider the field $\mathbb{Q}(X_1,\ldots,X_n)$ of all rational functions in the algebraically independent symbols $X_1,\ldots,X_n$. Let $\mathbb{P}_{\rm univ}(X_1,\ldots,X_n)$ be the set of all non-zero elements of $\mathbb{Q}(X_1,\ldots,X_n)$ that can be written as rational functions in $X_1,\ldots,X_n$ not involving subtractions (i.e. the `minus' operation). Then $\mathbb{P}_{\rm univ}(X_1,\ldots,X_n)$ becomes a semifield, when the addition and multiplication are inherited from those of $\mathbb{Q}(X_1,\ldots,X_n)$.

\vs

\item[\rm (2)] (the tropical semifield of $n$ variables) Let $y_1,\ldots,y_n$ be formal symbols, and let $\mathbb{P}_{\rm trop}(y_1,\ldots,y_n)$ be the set of all elements of the form $\prod_{i=1}^n y_i^{a_i}$ with $a_1,\ldots,a_n \in \mathbb{Z}$. Define the binary operations $+$ and $\cdot$ as
$$
\textstyle 
(\prod_{i=1}^n y_i^{a_i}) + (\prod_{i=1}^n y_i^{b_i}) := \prod_{i=1}^n y_i^{\min(a_i,b_i)},
\qquad
(\prod_{i=1}^n y_i^{a_i}) \cdot (\prod_{i=1}^n y_i^{b_i})
= \prod_{i=1}^n y_i^{a_i+b_i}.
$$
Then $\mathbb{P}_{\rm trop}(y_1,\ldots,y_n)$ is a semifield. 

\vs

\item[\rm (3)] (the tropical evaluation) There is a unique semifield homomorphism 
$$
{\bf T} : \mathbb{P}_{\rm univ}(X_1,\ldots,X_n) \to \mathbb{P}_{\rm trop}(y_1,\ldots,y_n)
$$
that satisfies
$$
{\bf T}(X_i) = y_i, \quad \forall i =1,\ldots,n, \qquad
{\bf T}(\alpha) = 1, \quad \forall \alpha \in \mathbb{Q}_{>0}. \qed
$$
\end{enumerate}
\end{lemma}

\begin{definition}
Let $\prod_{i=1}^n y_i^{a_i}$ be an element of the semifield $\mathbb{P}_{\rm trop}(y_1,\ldots,y_n)$, such that at least one $a_i$ is nonzero. Such an element $\prod_{i=1}^n y_i^{a_i}$ is said to be \ul{\em positive} if $a_i\ge 0$ holds for every $i=1,\ldots,n$, and \ul{\em negative} if $a_i\le 0$ holds for every $i=1,\ldots,n$.
\end{definition}

\begin{proposition}[sign-coherence; \cite{FoZ07} \cite{GHKK}]
\label{prop:sign}
Let $\Gamma = (\varepsilon, d, \{X_i\}_{i=1}^n)$ be an initial $\mathcal{X}$-seed. Let $\Gamma' \in |\Gamma|$, i.e. let $\Gamma' = (\varepsilon',d',\{X_i'\}_{i=1}^n)$ be an $\mathcal{X}$-seed equivalent to $\Gamma$, i.e. $\Gamma$ is connected to $\Gamma'$ by a cluster transformation. 

\vs

For each $i=1,\ldots,n$, note $X_i' \in \mathbb{P}_{\rm univ}(X_1,\ldots,X_n)$, so ${\bf T}(X_i') \in \mathbb{P}_{\rm trop}(y_1,\ldots,y_n)$, where $y_1,\ldots,y_n$ are formal symbols. Then ${\bf T}(X_i')$ is either positive, or negative.
\end{proposition}
The statement of Prop.\ref{prop:sign} is referred to as the `sign-coherence of ${\bf c}$-vectors' in the cluster algebra literature; it was first formulated as a conjecture in \cite{FoZ07}, proved for the cases of skew-symmetric exchange matrices $\varepsilon$ in \cite{DWZ} and for general cases in \cite{GHKK}. Thanks to Prop.\ref{prop:sign}, the following definition is well-defined.
\begin{definition}
Let $\Gamma,\Gamma'$, $X_i'$, ${\bf T}$, and $\mathbb{P}_{\rm trop}(y_1,\ldots,y_n)$ be as in Prop.\ref{prop:sign}. Define the \ul{\em sign} $\epsilon_\Gamma(X_i') \in \{+,-\}$ of the cluster variable $X_i'$ of $\Gamma'$ with respect to the initial seed $\Gamma$ as
$$
\epsilon_\Gamma(X_i') = \left\{
\begin{array}{ll}
+, & \mbox{if ${\bf T}(X_i')$ is a positive element of $\mathbb{P}_{\rm trop}(y_1,\ldots,y_n)$}, \\
-, & \mbox{if ${\bf T}(X_i')$ is a negative element of $\mathbb{P}_{\rm trop}(y_1,\ldots,y_n)$}.
\end{array}
\right.
$$
\end{definition}
Notice that, for a seed $\Gamma'$, it may be that the $n$ signs $\epsilon_\Gamma(X_1'),\ldots,\epsilon_\Gamma(X_n')$ are not all the same.

\begin{definition}[\cite{KN}]
\label{def:tropical_sign-sequence}
Suppose that ${\bf m} = {\bf m}_r \circ {\bf m}_{r-1} \circ \cdots \circ {\bf m}_1$ is a trivial cluster transformation of $\mathcal{X}$-seeds, being applied to the initial $\mathcal{X}$-seed $\Gamma = (\varepsilon,d,\{X_i\}_{i=1}^n)$. Use the notation for $\Gamma^{(i)}$ as in eq.\eqref{eq:seed_Gamma_i}, and write $\Gamma^{(i)} = (\varepsilon^{(i)},d^{(i)},\{X^{(i)}_j\}_{j=1}^n)$, for $i=0,1,\ldots,r-1$. For each $i$ for which ${\bf m}_i$ is a mutation, say it is $\mu_{k_i}$ for some $k_i \in \{1,\ldots,n\}$. Define the \ul{\em tropical sign-sequence} for the trivial cluster transformation ${\bf m}$ as the sequence
$$
\vec{\epsilon}_{\bf m} := (\epsilon_1,\epsilon_2,\ldots,\epsilon_r), \qquad \epsilon_i = \epsilon_{{\bf m};i} := \left\{
\begin{array}{ll}
\epsilon_\Gamma(X^{(i-1)}_{k_i}) \in \{+,-\}, & \mbox{if ${\bf m}_i$ is a mutation $\mu_{k_i}$}, \\
0, & \mbox{if ${\bf m}_i$ is a seed automorphism},
\end{array} \right.
 \quad \forall i=1,\ldots,r.
$$
\end{definition}
Finally, what we really need is the following.
\begin{proposition}[e.g. in the proof of {\cite[Thm.3.5]{KN}}]
\label{prop:product_of_matrices}
Suppose the situation as in Def.\ref{def:tropical_sign-sequence}. Suppose further that each ${\bf m}_i$ is a mutation, say $\mu_{k_i}$, and that each exchange matrix $\varepsilon^{(i)}$ is skew-symmetric.  For each $i=1,\ldots,r$, define the matrix ${\bf c}^{(\epsilon_i)}_{{\bf m}_i} \in {\rm SL}_\pm(n,\mathbb{R})$ by eq.\eqref{eq:bf_c_m_i_epsilon_i}. Then, one has
\begin{align}
\label{eq:bf_c_product_is_identity}
{\bf c}^{(\epsilon_1)}_{{\bf m}_1} \, {\bf c}^{(\epsilon_2)}_{{\bf m}_2} \, \cdots \, {\bf c}^{(\epsilon_r)}_{{\bf m}_r} = {\rm Id}.
\end{align}
\end{proposition}
First, notice that whether $\varepsilon^{(i)}$'s are skew-symmetric or merely skew-symmetrizable is determined by the initial one $\varepsilon$. So Prop.\ref{prop:product_of_matrices} is dealing with the case when $\varepsilon$ is skew-symmetric. I think the skew-symmetrizable version of Prop.\ref{prop:product_of_matrices} is implicit in the work \cite{GHKK}; here, let us try writing down an explicit proof, by mimicking \cite{KN}. In the meantime, here we are generalizing in another direction too, namely we are allowing seed automorphisms too. We shall try to interpret the matrices ${\bf c}^{(\epsilon_i)}_{{\bf m}_i}$ of eq.\eqref{eq:bf_c_m_i_epsilon_i}--\eqref{eq:bf_c_sigma_i} as coming from mutations and seed automorphisms of the cluster $\mathcal{X}$-variables. 

\vs

Fix an initial $\mathcal{X}$-seed $\Gamma$. Let $\Gamma' = (\varepsilon',d',\{X_i'\}_{i=1}^n)$ be any seed equivalent to $\Gamma$, and let $\mu_k(\Gamma') = \Gamma'' = (\varepsilon'',d'',\{X_i''\}_{i=1}^n)$. It is a straightforward exercise to deduce, from the mutation formula in eq.\eqref{eq:mu_k_on_X} for the cluster variables of $\Gamma''$ in terms of those of $\Gamma'$, the following tropical mutation formula for the tropical variables (see e.g. eq.(2.10) of \cite{KN}):
$$
{\bf T}(X_i'') = \left\{
\begin{array}{ll}
{\bf T}(X_k')^{-1} & \mbox{if $i=k$}, \\
{\bf T}(X_i') \cdot {\bf T}(X_k')^{[\epsilon_\Gamma(X_k') \varepsilon'_{ik}]_+} & \mbox{if $i\neq k$,}
\end{array}
\right.
$$
where the powers $-1$ and $[\epsilon_\Gamma(X_k') \varepsilon'_{ik}]_+$ are with respect to the tropical multiplication $\cdot$ in the semifield $\mathbb{P}_{\rm trop}(y_1,\ldots,y_n)$. Above tropical formula is an equality of elements of $\mathbb{P}_{\rm trop}(y_1,\ldots,y_n)$. Denote by 
$$
{\rm pow} : \mathbb{P}_{\rm trop}(y_1,\ldots,y_n) \to \mathbb{Z}^n
$$ 
the map ${\rm pow}(\prod_{i=1}^n y_i^{a_i})=(a_i)_{i=1}^n$, reading the powers; then ${\rm pow}$ is a group isomorphism, with respect to $\odot$. For a cluster variable $X_i'$ of a seed $\Gamma'$ equivalent to $\Gamma$, The vector $({\rm pow}\circ {\bf T})(X_i') \in \mathbb{Z}^n$ is called the `${\bf c}$-vector of $X_i'$' in \cite{FoZ07}. The tropical mutation formula can be written in terms of the ${\bf c}$-vectors as
$$
({\rm pow}\circ{\bf T})(X_i'') = \left\{
\begin{array}{ll}
- ({\rm pow}\circ{\bf T})(X_k') & \mbox{if $i=k$}, \\
({\rm pow}\circ{\bf T}) (X_i') + [\epsilon_\Gamma(X_k') \varepsilon'_{ik}]_+ \, \cdot ({\rm pow}\circ{\bf T})(X_k') & \mbox{if $i\neq k$,}
\end{array}
\right.
$$
Denote by $({\rm pow}\circ {\bf T})(\Gamma')$ the $n\times n$ matrix whose $i$-th row is the row vector $({\rm pow}\circ {\bf T})(X_i')$. Then the above formula can be written as
$$
({\rm pow}\circ {\bf T})(\Gamma'') = \til{\bf c}^{(\epsilon_\Gamma(X_k'))}_{\Gamma'\mut{k}\Gamma''} \cdot ({\rm pow}\circ {\bf T})(\Gamma') 
$$
where the $n\times n$ matrix $\til{\bf c}^{(\epsilon)}_{\Gamma'\mut{k}\Gamma''} = (\til{c}^{(\epsilon)}_{ij})_{i,j}$ is given by
\begin{align}
\label{eq:bf_c_tilde}
\til{c}^{(\epsilon)}_{ii} = 1, \quad \forall i \neq k, \qquad \til{c}^{(\epsilon)}_{kk} = -1, \qquad \til{c}^{(\epsilon)}_{ik} = [\epsilon\, \varepsilon'_{ik}]_+, \quad \forall i\neq k, \qquad \til{c}^{(\epsilon)}_{ij}=0 \quad\mbox{otherwise}.
\end{align}
In particular, in case when $\varepsilon'$ is skew-symmetric, one has $\varepsilon'_{ik} = - \varepsilon'_{ki}$, and hence this matrix $\til{\bf c}^{(\epsilon)}_{\Gamma'\mut{k}\Gamma''}$ equals the matrix ${\bf c}^{(\epsilon)}_{\Gamma'\mut{k}\Gamma''} \in {\rm SL}_\pm(n,\mathbb{R})$ define in eq.\eqref{eq:c_Gamma_Gamma_prime_entries_epsilon}.

\vs

On the other hand, this time suppose $P_\sigma(\Gamma') = \Gamma''$, for some permutation $\sigma$ of $\{1,\ldots,n\}$. Then from Def.\ref{def:seed_automorphism},
$$
{\bf T}(X_{\sigma(i)}'') = {\bf T}(X_i'), \quad \forall i=1,\ldots,n,
$$
which can be written in the form of
$$
({\rm pow}\circ {\bf T})(\Gamma'') = \til{\bf c}^{(0)}_{\sigma} \cdot ({\rm pow}\circ {\bf T})(\Gamma'),
$$
where the $n\times n$ matrix $\til{\bf c}^{(0)}_\sigma = (\til{c}^{(0)}_{ij})_{i,j}$ is given by
\begin{align}
\label{eq:til_bf_c_sigma}
\til{c}^{(0)}_{ij} = \delta_{\sigma(i),\,j}, \quad \forall i,j.
\end{align}

\vs

Now, coming back to our cluster transformation ${\bf m}$, for a mutation ${\bf m}_i = \mu_{k_i}$ let $\til{\bf c}^{(\epsilon_i)}_{{\bf m}_i} = \til{\bf c}^{(\epsilon_i)}_{\Gamma^{(i-1)}\mut{k_i}\Gamma^{(i)}}$, and for a seed automorphism ${\bf m}_i = P_{\sigma_i}$ let $\til{\bf c}^{(\epsilon_i)}_{{\bf m}_i} = \til{\bf c}^{(0)}_{\sigma_i}$. Using these notations, one has
\begin{align}
\label{eq:equality_for_bf_c_tilde}
({\rm pow}\circ {\bf T})(\Gamma^{(0)}) = ({\rm pow}\circ {\bf T})(\Gamma^{(r)}) = 
\underbrace{ \til{\bf c}^{(\epsilon_r)}_{{\bf m}_r} \, \cdots \,\, \til{\bf c}^{(\epsilon_1)}_{{\bf m}_1}}
\, \cdot \, ({\rm pow}\circ {\bf T})(\Gamma^{(0)})
\end{align}
Since $({\rm pow}\circ {\bf T})(\Gamma^{(0)})$ is the identity matrix, it follows that the underbraced part in the above equation equals the identity matrix. In case when all ${\bf m}_i$'s are mutations and all $\varepsilon^{(i)}$'s are skew-symmetric, each matrix $\til{\bf c}^{(\epsilon_i)}_{{\bf m}_i}$ equals ${\bf c}^{(\epsilon_i)}_{{\bf m}_i}$, and in the meantime, it is its own inverse. Hence indeed the sought-for eq.\eqref{eq:bf_c_product_is_identity} follows.

\vs

We now should generalize this proof to the case when $\varepsilon^{(i)}$'s are not necessarily skew-symmetric, and not necessarily all ${\bf m}_i$'s are mutations. It is easy to deal with the seed automorphisms ${\bf m}_i =P_{\sigma_i}$, for one can easily see that $\til{\bf c}^{(0)}_{\sigma_i}$ defined in eq.\eqref{eq:til_bf_c_sigma} is the inverse of the matrix ${\bf c}^{(0)}_{\sigma_i}$ defined in eq.\eqref{eq:bf_c_sigma_i}. Now, for skew-symmetrizable exchange matrices, $\varepsilon'_{ik} = - \varepsilon_{ki}'$ need not hold. So, although the underbraced product of the $\til{\bf c}$ matrices in eq.\eqref{eq:equality_for_bf_c_tilde} is still the identity matrix, it is not clear whether the sought-for product ${\bf c}^{(\epsilon_1)}_{{\bf m}_1} \, {\bf c}^{(\epsilon_2)}_{{\bf m}_2} \, \cdots \, {\bf c}^{(\epsilon_r)}_{{\bf m}_r}$ is also the identity matrix. One way to reformulate the statement is as follows. For any $\mathcal{X}$-seed $\Gamma=(\varepsilon,d,\{X_i\}_{i=1}^n)$, define its \ul{\em Langlands dual} seed as the $\mathcal{X}$-seed $\Gamma^\vee = (\varepsilon^\vee, d^\vee, \{X_i^\vee\}_{i=1}^n)$, where
$$
\varepsilon^\vee_{ij} = -\varepsilon_{ji},  \quad d^\vee_i = N d_i^{-1},
$$
where $N$ is the least common multiple of $d_i$'s. Since $\Gamma^\vee$ is an example of an $\mathcal{X}$-seed, one can consider mutations and seed automorphisms. One can easily check that the transformation formulas for $\varepsilon$ and $d$ under mutations and seed automorphisms commute with the dualizing maps $\varepsilon \mapsto \varepsilon^\vee$, $d\mapsto d^\vee$. However, it is not clear how to establish a dualizing map for the cluster variables so that it commutes with the mutation formulas; so, at the moment, one cannot really refer to a map or correspondence $\Gamma \mapsto \Gamma^\vee$. Anyhow, view ${\bf m}$ as the sequence of symbols $\mu_{k_r} \circ \mu_{k_{r-1}} \circ \cdots \circ \mu_{k_1}$, which was a priori applied to an initial seed $\Gamma=\Gamma^{(0)}$. Now, consider applying this sequence of symbols to the dual seed $\Gamma^\vee$, to obtain a dual version cluster transformation ${\bf m}^\vee = {\bf m}_r^\vee \circ \cdots \circ {\bf m}^\vee_1$ starting from $\Gamma^\vee = (\Gamma^{(0)})^\vee$. Suppose now that this new cluster transformation ${\bf m}^\vee$ is also a trivial cluster transformation, i.e. it starts and ends at a same seed $\Gamma^\vee$. Denote by $\vec{\epsilon}^{\,\vee}$ the tropical sign-sequence of ${\bf m}^\vee$. Suppose further that $\vec{\epsilon}^{\,\vee}$ equals the original one $\vec{\epsilon}$ for ${\bf m}$. Then, by applying our argument above, we obtain the sought-for eq.\eqref{eq:bf_c_product_is_identity}. One can formulate a relevant conjecture as follows.
\begin{conjecture}
\label{conj:dual_triviality}
Let ${\bf m} = \mu_{k_r}\circ \cdots \circ \mu_{k_1}$ be a sequence of mutations and seed automorphisms. Suppose that ${\bf m}$ is a trivial cluster transformation when applied to some $\mathcal{X}$-seed $\Gamma$. Then, 
\begin{enumerate}
\item[\rm (1)] when this sequence is applied to the Langlands dual $\mathcal{X}$-seed $\Gamma^\vee$, it is still a trivial cluster transformation;

\item[\rm (2)] the tropical sign-sequence for ${\bf m}$ when applied to $\Gamma$ coincides with that for ${\bf m}$ when applied to $\Gamma^\vee$.
\end{enumerate}
\end{conjecture}
\begin{lemma}
\label{lem:cancellation_of_linear_parts_general}
For each sequence ${\bf m}$ for which Conjecture \ref{conj:dual_triviality} is true, the statement of Prop.\ref{prop:product_of_matrices} holds, without the stipulation that all ${\bf m}_i$'s be mutations and $\varepsilon$ be skew-symmetric.
\end{lemma}
For the purpose of the present paper, we only deal with the trivial morphisms ${\bf m}$ of the saturated cluster modular groupoid, i.e. only the ones corresponding to (S1)--(S2) in \S\ref{sec:introduction}. In the next subsection, we will show that the Conjecture \ref{conj:dual_triviality} holds for these trivial cluster transformations.

\subsection{Triviality of the linear parts for the rank 1 and rank 2 relations}
\label{subsec:triviality_of_linear_parts}

In this subsection, we check Conjecture \ref{conj:dual_triviality} which is formulated in the previous subsection, for the trivial cluster transformations ${\bf m}$ of (S1) and (S2) which appeared in \S\ref{sec:introduction}, or equivalently, those written in Lemmas \ref{lem:involution_identity} and \ref{lem:h_plus_2-gon_relations}. Part (1) of the conjecture is already known to hold for these ${\bf m}$'s, so we only need to check part (2). We follow the notations in the previous subsection, unless otherwise specified. In particular, we use $\Gamma = (\varepsilon,d,\{X_i\}_{i=1}^n)=\Gamma^{(0)}$ for the initial $\mathcal{X}$-seed, and use the tropicalization map ${\bf T}$ from $\mathbb{P}_{\rm univ}(X_1,\ldots,X_n)$ to $\mathbb{P}_{\rm trop}(y_1,\ldots,y_n)$.

\vs

\ul{\em Case (S1) : the rank 1 identity, or $A_1$ identity, or the involution identity}

\vs

Consider ${\bf m} = \mu_k \circ \mu_k$, and consider applying it to an $\mathcal{X}$-seed $\Gamma = \Gamma^{(0)}$. Indeed, for any $\mathcal{X}$-seed, this is a trivial cluster transformation (Lem.\ref{lem:involution_identity}). Let us compute the corresponding tropical sign-sequence $\vec{\epsilon} = (\epsilon_1,\epsilon_2)$. Note ${\bf T}( X_k^{(0)} ) = {\bf T}(X_k) = y_k$, so $\epsilon_1 = \epsilon_\Gamma(X_k^{(0)})=+$. Note $X_k^{(1)} = (X_k^{(0)})^{-1} = X_k^{-1}$, so ${\bf T}(X_k^{(1)}) = y_k^{-1}$, hence $\epsilon_2 = \epsilon_\Gamma(X_k^{(1)})=-$. We did not use any assumption other than ${\bf m}$ is $\mu_k \circ \mu_k$. So Conjecture \ref{conj:dual_triviality} holds for this ${\bf m}$.

\vs

\ul{\em Case (S2-1) : the $A_1\times A_1$ identity}

\vs

Consider ${\bf m} = \mu_j \circ \mu_i \circ \mu_j \circ \mu_i$, applied to an $\mathcal{X}$-seed $\Gamma = \Gamma^{(0)}$ such that $\varepsilon_{ij}=0$ with $i\neq j$. Indeed it is a trivial cluster transformation (Lem.\ref{lem:h_plus_2-gon_relations}). Let us compute the tropical sign-sequence $\vec{\epsilon}=(\epsilon_1,\epsilon_2,\epsilon_3,\varepsilon_4)$. Note ${\bf T}(X_i^{(0)}) = {\bf T}(X_i) = y_i$, so $\epsilon_1 = \epsilon_\Gamma(X_i^{(0)})=+$. Note $X_j^{(1)} = X_j^{(0)}$ by eq.\eqref{eq:mu_k_on_X} (mutation at $i$), so ${\bf T}(X_j^{(1)}) = {\bf T}(X_j) = y_j$, so $\epsilon_2=+$. Note $\varepsilon^{(1)}_{ij} = -\varepsilon^{(0)}_{ij}=0$ (eq.\eqref{eq:varepsilon_prime_formula}), so $X_i^{(2)} = X_i^{(1)}$ (mutation at $j$) while $X_i^{(1)} = (X_i^{(0)})^{-1}$ (mutation at $i$), so ${\bf T}(X_i^{(2)}) = {\bf T}(X_i^{-1}) = y_i^{-1}$, hence $\epsilon_3 = -$. Note $\varepsilon^{(2)}_{ij} = - \varepsilon^{(1)}_{ij} = 0$ (eq.\eqref{eq:varepsilon_prime_formula}), so $X_j^{(3)} = X_j^{(2)}$ (mutation at $i$), while $X_j^{(2)} = (X_j^{(1)})^{-1}$ (mutation at $j$), and $X_j^{(1)} = X_j^{(0)}$ (mutation at $i$), hence ${\bf T}(X_j^{(3)}) = {\bf T}(X_j^{-1})= y_j^{-1}$, so $\epsilon_4=-$. So $\vec{\epsilon}=(+,+,-,-)$, and in particular, Conjecture \ref{conj:dual_triviality} holds for this ${\bf m}$. 

\vs

\ul{\em Case (S2-2) : the $A_2$ identity}

\vs

Consider ${\bf m} = P_{(i\,j)} \circ \mu_i \circ \mu_j \circ \mu_i \circ \mu_j \circ \mu_i$, applied to $\Gamma = \Gamma^{(0)}$ such that $\varepsilon_{ij}=-\varepsilon_{ji}\in\{1,-1\}$. First, suppose
$$
\varepsilon_{ij} = -\varepsilon_{ji}=1.
$$
By the mutation formula in eq.\eqref{eq:varepsilon_prime_formula} for the exchange matrices, one has $\varepsilon^{(\ell)}_{ij} = - \varepsilon^{(\ell)}_{ji} = 1$ for $\ell = 0,2$ and $\varepsilon^{(\ell)}_{ij} = - \varepsilon^{(\ell)}_{ji} = -1$ for $\ell = 1,3$. Note
\begin{align*}
\Gamma^{(0)} & : {\bf T}(X^{(0)}_i) = {\bf T}(X_i) = y_i, \quad \epsilon_1=+, \\
\Gamma^{(1)} = \mu_i(\Gamma^{(0)}) & : 
{\bf T}(X_i^{(1)}) = {\bf T}(X_i^{-1}) = y_i^{-1}, \\
& \quad  {\bf T}(X_j^{(1)}) = {\bf T}(X_j(1+X_i)) = y_j, \quad \epsilon_2=+, \\
\Gamma^{(2)} = \mu_j(\Gamma^{(1)}) & : 
{\bf T}(X_i^{(2)}) = {\bf T}(X_i^{(1)}(1+X_j^{(1)})) = y_i^{-1}(1+y_j) = y_i^{-1}, \quad \epsilon_3 = -, \\
& \quad {\bf T}(X_j^{(2)}) = {\bf T}((X_j^{(1)})^{-1}) = y_j^{-1}, \\
\Gamma^{(3)} = \mu_i(\Gamma^{(2)}) & : 
{\bf T}(X_i^{(3)}) = {\bf T}((X_i^{(2)})^{-1}) = y_i, \\
& \quad {\bf T}(X_j^{(3)}) = {\bf T}(X_j^{(2)}(1+X_i^{(2)})) = y_j^{-1}(1+y_i^{-1}) = y_i^{-1}y_j^{-1}, \quad \epsilon_4 = -, \\
\Gamma^{(4)} = \mu_j(\Gamma^{(3)}) & : 
{\bf T}(X_i^{(4)}) = {\bf T}(X_i^{(3)}(1+X_j^{(3)})) = y_i(1+y_i^{-1}y_j^{-1}) = y_j^{-1}, \quad \epsilon_5=-.
\end{align*}
So the tropical sign-sequence is $\vec{\epsilon}=(+,+,-,-,-,0)$; in particular, Conjecture \ref{conj:dual_triviality} holds for this ${\bf m}$. Such computation of tropicalizations of $\mathcal{X}$-variables are straightforward, as shown above. Hence, from now on, I will only record the results of computations.

\vs

Now suppose
$$
\varepsilon_{ij} = - \varepsilon_{ji} = -1.
$$
Through similar computations as above, one can obtain
$$
{\bf T}(X_i^{(0)}) = y_i, \quad
{\bf T}(X_j^{(1)}) = y_i y_j, \quad
{\bf T}(X_i^{(2)}) = y_j, \quad
{\bf T}(X_j^{(3)}) = y_i^{-1}, \quad
{\bf T}(X_i^{(4)}) = y_j^{-1},
$$
so that the tropical sign-sequence is $\vec{\epsilon} = (+,+,+,-,-,0)$; in particular, Conjecture \ref{conj:dual_triviality} holds for this ${\bf m}$. 

\vs

\ul{\em Case (S2-3) : the $B_2$ identity}

\vs

Consider ${\bf m} = \mu_j \circ \mu_i \circ \mu_j \circ \mu_i \circ \mu_j \circ \mu_i$, applied to $\Gamma = \Gamma^{(0)}$ such that $\varepsilon_{ij}=-2\varepsilon_{ji}\in\{2,-2\}$ or $\varepsilon_{ji}=-2\varepsilon_{ij}\in\{2,-2\}$. First, suppose
$$
\varepsilon_{ij} = 2, \quad \varepsilon_{ji} = -1.
$$
By the mutation formula in eq.\eqref{eq:varepsilon_prime_formula} for the exchange matrices, one has $\varepsilon^{(\ell)}_{ij} = 2$ and $\varepsilon^{(\ell)}_{ji} = -1$ for $\ell = 0,2,4$, while $\varepsilon^{(\ell)}_{ij} = - 2$ and  $\varepsilon^{(\ell)}_{ji} = 1$ for $\ell = 1,3$. By similar computations as before, one obtains
$$
{\bf T}(X_i^{(0)}) = y_i, \quad
{\bf T}(X_j^{(1)}) = y_j, \quad
{\bf T}(X_i^{(2)}) = y_i^{-1}, \quad
{\bf T}(X_j^{(3)}) = y_i^{-1} y_j^{-1}, \quad
{\bf T}(X_i^{(4)}) = y_i^{-1} y_j^{-2}, \quad
{\bf T}(X_j^{(5)}) = y_j^{-1},
$$
so the tropical sign-sequence is $(+,+,-,-,-,-)$. In order to check Conjecture \ref{conj:dual_triviality}, assume now
$$
\varepsilon_{ij} = 1, \quad \varepsilon_{ji} = -2.
$$
By the mutation formula in eq.\eqref{eq:varepsilon_prime_formula} for the exchange matrices, one has $\varepsilon^{(\ell)}_{ij} = 1$ and $\varepsilon^{(\ell)}_{ji} = -2$ for $\ell = 0,2,4$, while $\varepsilon^{(\ell)}_{ij} = - 1$ and  $\varepsilon^{(\ell)}_{ji} = 2$ for $\ell = 1,3$. By similar computations, one gets
$$
{\bf T}(X_i^{(0)}) = y_i, \quad
{\bf T}(X_j^{(1)}) = y_j, \quad
{\bf T}(X_i^{(2)}) = y_i^{-1}, \quad
{\bf T}(X_j^{(3)}) = y_i^{-2} y_j^{-1}, \quad
{\bf T}(X_i^{(4)}) = y_i^{-1} y_j^{-1}, \quad
{\bf T}(X_j^{(5)}) = y_j^{-1},
$$
so the tropical sign-sequence is again $(+,+,-,-,-,-)$, hence indeed Conjecture \ref{conj:dual_triviality} holds for this ${\bf m}$.

\vs

Now, suppose
$$
\varepsilon_{ij} = -2, \quad \varepsilon_{ji} = 1.
$$
By the mutation formula in eq.\eqref{eq:varepsilon_prime_formula} for the exchange matrices, one has $\varepsilon^{(\ell)}_{ij} = -2$ and $\varepsilon^{(\ell)}_{ji} = 1$ for $\ell = 0,2,4$, while $\varepsilon^{(\ell)}_{ij} = 2$ and  $\varepsilon^{(\ell)}_{ji} = -1$ for $\ell = 1,3$. One can compute
$$
{\bf T}(X_i^{(0)}) = y_i, \quad
{\bf T}(X_j^{(1)}) = y_i y_j, \quad
{\bf T}(X_i^{(2)}) = y_i y_j^2, \quad
{\bf T}(X_j^{(3)}) = y_j, \quad
{\bf T}(X_i^{(4)}) = y_i^{-1}, \quad
{\bf T}(X_j^{(5)}) = y_j^{-1},
$$
so the tropical sign-sequence is $(+,+,+,+,-,-)$. In order to check Conjecture \ref{conj:dual_triviality}, assume now
$$
\varepsilon_{ij} = -1, \quad \varepsilon_{ji} = 2.
$$
By the mutation formula in eq.\eqref{eq:varepsilon_prime_formula} for the exchange matrices, one has $\varepsilon^{(\ell)}_{ij} = -1$ and $\varepsilon^{(\ell)}_{ji} = 2$ for $\ell = 0,2,4$, while $\varepsilon^{(\ell)}_{ij} = 1$ and  $\varepsilon^{(\ell)}_{ji} = -2$ for $\ell = 1,3$. One computes
$$
{\bf T}(X_i^{(0)}) = y_i, \quad
{\bf T}(X_j^{(1)}) = y_i^2 y_j, \quad
{\bf T}(X_i^{(2)}) = y_i y_j, \quad
{\bf T}(X_j^{(3)}) = y_j, \quad
{\bf T}(X_i^{(4)}) = y_i^{-1}, \quad
{\bf T}(X_j^{(5)}) = y_j^{-1},
$$
so the tropical sign-sequence is again $(+,+,+,+,-,-)$, hence indeed Conjecture \ref{conj:dual_triviality} holds for this ${\bf m}$.

\vs

\ul{\em Case (S2-4) : the $G_2$ identity}

\vs

Consider ${\bf m} = \mu_j \circ \mu_i \circ \mu_j \circ \mu_i \circ \mu_j \circ \mu_i \circ \mu_j \circ \mu_i$, applied to $\Gamma = \Gamma^{(0)}$ such that $\varepsilon_{ij}=-3\varepsilon_{ji}\in\{3,-3\}$ or $\varepsilon_{ji}=-3\varepsilon_{ij}\in\{3,-3\}$. First, suppose
$$
\varepsilon_{ij} = 3, \quad \varepsilon_{ji} = -1.
$$
By the mutation formula in eq.\eqref{eq:varepsilon_prime_formula} for the exchange matrices, one has $\varepsilon^{(\ell)}_{ij} = 3$ and $\varepsilon^{(\ell)}_{ji} = -1$ for $\ell = 0,2,4,6$, while $\varepsilon^{(\ell)}_{ij} = - 3$ and  $\varepsilon^{(\ell)}_{ji} = 1$ for $\ell = 1,3,5$. One computes
\begin{align*}
& {\bf T}(X_i^{(0)}) = y_i, \quad
{\bf T}(X_j^{(1)}) = y_j, \quad
{\bf T}(X_i^{(2)}) = y_i^{-1}, \quad
{\bf T}(X_j^{(3)}) = y_i^{-1} y_j^{-1}, \quad
{\bf T}(X_i^{(4)}) = y_i^{-2} y_j^{-3}, \\
& {\bf T}(X_j^{(5)}) = y_i^{-1} y_j^{-2}, \quad
{\bf T}(X_i^{(6)}) = y_i^{-1} y_j^{-3}, \quad
{\bf T}(X_j^{(7)}) = y_j^{-1}, \quad
\end{align*}
so the tropical sign-sequence is $(+,+,-,-,-,-,-,-)$.  In order to check Conjecture \ref{conj:dual_triviality}, assume now
$$
\varepsilon_{ij} = 1, \quad \varepsilon_{ji} = -3.
$$
By the mutation formula in eq.\eqref{eq:varepsilon_prime_formula} for the exchange matrices, one has $\varepsilon^{(\ell)}_{ij} = 1$ and $\varepsilon^{(\ell)}_{ji} = -3$ for $\ell = 0,2,4,6$, while $\varepsilon^{(\ell)}_{ij} = - 1$ and  $\varepsilon^{(\ell)}_{ji} = 3$ for $\ell = 1,3,5$. One computes
\begin{align*}
& {\bf T}(X_i^{(0)}) = y_i, \quad
{\bf T}(X_j^{(1)}) = y_j, \quad
{\bf T}(X_i^{(2)}) = y_i^{-1}, \quad
{\bf T}(X_j^{(3)}) = y_i^{-3}y_j^{-1}, \quad
{\bf T}(X_i^{(4)}) = y_i^{-2} y_j^{-1}, \\
& {\bf T}(X_j^{(5)}) = y_i^{-3} y_j^{-2}, \quad
{\bf T}(X_i^{(6)}) = y_i^{-1} y_j^{-1}, \quad
{\bf T}(X_j^{(7)}) = y_j^{-1}, \quad
\end{align*}
so the tropical sign-sequence is again $(+,+,-,-,-,-,-,-)$, hence indeed Conjecture \ref{conj:dual_triviality} holds for this ${\bf m}$.

\vs

Now, suppose
$$
\varepsilon_{ij} = -3, \quad \varepsilon_{ji} = 1.
$$
By the mutation formula in eq.\eqref{eq:varepsilon_prime_formula} for the exchange matrices, one has $\varepsilon^{(\ell)}_{ij} = -3$ and $\varepsilon^{(\ell)}_{ji} = 1$ for $\ell = 0,2,4,6$, while $\varepsilon^{(\ell)}_{ij} = 3$ and  $\varepsilon^{(\ell)}_{ji} = -1$ for $\ell = 1,3,5$. One computes
\begin{align*}
& {\bf T}(X_i^{(0)}) = y_i, \quad
{\bf T}(X_j^{(1)}) = y_i y_j, \quad
{\bf T}(X_i^{(2)}) = y_i^2 y_j^3, \quad
{\bf T}(X_j^{(3)}) = y_i y_j^2, \quad
{\bf T}(X_i^{(4)}) = y_i y_j^3, \\
& {\bf T}(X_j^{(5)}) = y_j, \quad
{\bf T}(X_i^{(6)}) = y_i^{-1}, \quad
{\bf T}(X_j^{(7)}) = y_j^{-1}, \quad
\end{align*}
so the tropical sign-sequence is $(+,+,+,+,+,+,-,-)$.  In order to check Conjecture \ref{conj:dual_triviality}, assume now
$$
\varepsilon_{ij} = -1, \quad \varepsilon_{ji} = 3.
$$
By the mutation formula in eq.\eqref{eq:varepsilon_prime_formula} for the exchange matrices, one has $\varepsilon^{(\ell)}_{ij} = -1$ and $\varepsilon^{(\ell)}_{ji} = 3$ for $\ell = 0,2,4,6$, while $\varepsilon^{(\ell)}_{ij} = 1$ and  $\varepsilon^{(\ell)}_{ji} = -3$ for $\ell = 1,3,5$. One computes
\begin{align*}
& {\bf T}(X_i^{(0)}) = y_i, \quad
{\bf T}(X_j^{(1)}) = y_i^3 y_j, \quad
{\bf T}(X_i^{(2)}) = y_i^2 y_j, \quad
{\bf T}(X_j^{(3)}) = y_i^3 y_j^2, \quad
{\bf T}(X_i^{(4)}) = y_i y_j, \\
& {\bf T}(X_j^{(5)}) = y_j, \quad
{\bf T}(X_i^{(6)}) = y_i^{-1}, \quad
{\bf T}(X_j^{(7)}) = y_j^{-1}, \quad
\end{align*}
so the tropical sign-sequence is again $(+,+,+,+,+,+,-,-)$, hence indeed Conjecture \ref{conj:dual_triviality} holds for this ${\bf m}$.

\subsection{From double version identity to single version identity}
\label{subsec:from_double_to_single}

Let ${\bf m} = {\bf m}_r \circ \cdots \circ {\bf m}_1$ be a trivial cluster transformation of $\mathcal{D}$-seeds, with initial seed being $(\varepsilon,d,\{B_i,X_i\}_{i=1}^n)$ and each ${\bf m}_i$ is represented as a mutation $\mu_{k_i}$ or a seed automorphism $P_{\sigma_i}$. View this sequence ${\bf m}$ of mutations and seed automorphisms as being applied to the $\mathcal{X}$-seed $(\varepsilon,d,\{X_i\}_{i=1}^n)$, and suppose that it satisfies the statement of Conjecture \ref{conj:dual_triviality}, which in fact is a statement about $(\varepsilon,d)$, $k_i$'s, $d_i$'s, and $\sigma_i$'s (but not really about specific $X_i$'s or $B_i$'s). Then our arguments in \S\ref{subsec:collectin_the_linear_parts} apply, so eq.\eqref{eq:bf_c_product_is_identity} holds, and eq.\eqref{eq:bf_S_only} equals identity. Suppose also that $\pi^\hbar({\bf m}) = \pi^\hbar({\bf m}_1) \circ \cdots \circ \pi^\hbar({\bf m}_r)$ is a scalar operator. Then it follows that eq.\eqref{eq:quantum_dilog_factors_only} holds. Recall that the identity in eq.\eqref{eq:quantum_dilog_factors_only} is an equation having $2r'$ factors involving the quantum dilogarithm function, where $r'$ of them are results of functional calculus for the operators ${\bf x}^\hbar_{{\bf m}_i}$ and the others for $\til{\bf x}^\hbar_{{\bf m}_i}$. We will now `separate' these two groups of $r'$ factors, to obtain two identities. 
\begin{definition}
\label{def:operator-separation_property}
Consider self-adjoint operators $A_1,A_2,\ldots,B_2,B_2,\ldots$ on the Hilbert space $\mathscr{H}_\Gamma$. We say that $A_1,A_2,\ldots$ are \ul{\em strongly separated} with $B_1,B_2,\ldots$ if there exists a unitary operator $U : \mathscr{H}_\Gamma \to \mathscr{H}_\Gamma$ and a partition $\{1,\ldots,n\} = I_A \sqcup I_B$ of the indexing set  such that each of $U A_1 U^{-1}, \, U A_2 U^{-1}, \ldots$ is an $\mathbb{R}$-linear combination of the operators in $\{ {\bf p}^\hbar_i, {\bf q}^\hbar_i : i \in I_A \}$, while each of $U B_1 U^{-1}, \, U B_2 U^{-1}, \ldots$ is an $\mathbb{R}$-linear combination of the operators in $\{ {\bf p}^\hbar_i, {\bf q}^\hbar_i : i \in I_B \}$.

\vs

We say that the cluster transformation ${\bf m} = {\bf m}_r \circ \cdots \circ {\bf m}_1$ satisfies the \ul{\em operator-separation} property if the non-tilde operators ${\bf x}^\hbar_{{\bf m}_1}$, \ldots, ${\bf x}^\hbar_{{\bf m}_r}$ are strongly separated with the tilde operators $\til{\bf x}^\hbar_{{\bf m}_1}$, \ldots, $\til{\bf x}^\hbar_{{\bf m}_r}$.
\end{definition}

Later in the next section, we will show that this operator-separation property holds for each ${\bf m}$ that we will investigate. For now, suppose ${\bf m}$ satisfies the operator-separation property, and $U$ be an operator as in Def.\ref{def:operator-separation_property}; for convenience, let $I_A = \{1,\ldots,\ell\}$ and $I_B = \{\ell+1,\ldots,n\}$ be the partition of $\{1,\ldots,n\}$ in Def.\ref{def:operator-separation_property}. Now, conjugate by $U$ to both sides of eq.\eqref{eq:quantum_dilog_factors_only} to obtain the operator identity
\begin{align}
\label{eq:quantum_dilog_factors_only_identity}
\Phi^{\hbar_{k_1}}(\epsilon_1\, U {\bf x}^\hbar_{{\bf m}_1} U^{-1})^{\epsilon_1}
\, \Phi^{\hbar_{k_1}}(\epsilon_1\, U \til{\bf x}^\hbar_{{\bf m}_1} U^{-1})^{-\epsilon_1} \, \cdots \,
\Phi^{\hbar_{k_r}}(\epsilon_r\, U{\bf x}^\hbar_{{\bf m}_r}U^{-1})^{\epsilon_r}
\, \Phi^{\hbar_{k_r}}(\epsilon_r\, U\til{\bf x}^\hbar_{{\bf m}_r}U^{-1})^{-\epsilon_r} = c_{{\bf m}} \cdot {\rm Id}
\end{align}
on the Hilbert space $\mathscr{H}_\Gamma = L^2(\mathbb{R}^n, da_1\cdots da_n) = L^2(\mathbb{R}^\ell, da_1\cdots da_\ell) \otimes L^2(\mathbb{R}^{n-\ell}, da_{\ell+1} \cdots da_n)$, where $\otimes$ is the tensor product of Hilbert spaces, and $c_{\bf m} \in {\rm U}(1)$ is some scalar. We will now prove that the above operator identity yields two operator identities for each tensor factor, namely
\begin{align}
\label{eq:separated_identity1}
& \Phi^{\hbar_{k_1}}(\epsilon_1\, U {\bf x}^\hbar_{{\bf m}_1} U^{-1})^{\epsilon_1}
\, \cdots \,
\Phi^{\hbar_{k_r}}(\epsilon_r\, U{\bf x}^\hbar_{{\bf m}_r}U^{-1})^{\epsilon_r}
= c_{{\bf m}}' \cdot {\rm Id} \qquad \mbox{on} \quad L^2(\mathbb{R}^\ell, da_1\cdots da_\ell) \quad\mbox{and} \\\label{eq:separated_identity2}
& \Phi^{\hbar_{k_1}}(\epsilon_1\, U \til{\bf x}^\hbar_{{\bf m}_1} U^{-1})^{-\epsilon_1}
\, \cdots \,
\Phi^{\hbar_{k_r}}(\epsilon_r\, U \til{\bf x}^\hbar_{{\bf m}_r}U^{-1})^{-\epsilon_r}
= c_{{\bf m}}'' \cdot {\rm Id} \quad \mbox{on} \quad L^2(\mathbb{R}^{n-\ell}, da_{\ell+1}\cdots da_n)
\end{align}
hold for some scalars $c'_{\bf m}, c''_{\bf m} \in {\rm U}(1)$ such that
$$
c'_{\bf m} \, c''_{\bf m} = c_{\bf m};
$$
to be more precise, we should have denoted the operators $U {\bf x}^\hbar_{{\bf m}_i} U^{-1}$ and $U \til{\bf x}^\hbar_{{\bf m}_j} U^{-1}$ by different names because we are now regarding them as operators on the respective tensor factors.

\vs

To prove these separated identies, consider any elements $\psi(a_1,\ldots,a_\ell) \in L^2(\mathbb{R}^\ell,da_1\cdots da_\ell)$ and $\eta(a_{\ell+1},\ldots,a_n) \in L^2(\mathbb{R}^{n-\ell},da_{\ell+1}\cdots da_n)$, so that $(\psi\otimes \eta)(a_1,\ldots,a_n) := \psi(a_1,\ldots,a_\ell) \, \eta(a_{\ell+1},\ldots,a_n)$ belongs to $\mathscr{H}_\Gamma$; elements of $\mathscr{H}_\Gamma$ obtained like this form a dense subspace of $\mathscr{H}_\Gamma$. Denote the left-hand-side of eq.\eqref{eq:separated_identity1} as $A$, and the left-hand-side of eq.\eqref{eq:separated_identity2} as $B$. Each of $A,B$ can be understood as operators on $\mathscr{H}_\Gamma$, or on the respective tensor factor $L^2(\mathbb{R}^\ell,da_1\cdots da_\ell)$ or $L^2(\mathbb{R}^{n-\ell},da_{\ell+1}\cdots da_n)$, thanks to the separation assumption of Def.\ref{def:operator-separation_property}.  Hence, if we denote by $C$ the left-hand-side of eq.\eqref{eq:quantum_dilog_factors_only_identity}, we get
$$
(C(\psi\otimes \eta))(a_1,\ldots,a_n)
= (A\psi)(a_1,\ldots,a_\ell) \, (B\eta)(a_{\ell+1},\ldots,a_n)
$$
which equals
$$
c_{\bf m} \, (\psi\otimes \eta)(a_1,\ldots,a_n) = c_{\bf m} \, \psi(a_1,\ldots,a_\ell) \, \eta(a_{\ell+1},\ldots,a_n),
$$
by eq.\eqref{eq:quantum_dilog_factors_only_identity}.  This must hold for any $\psi$ and $\eta$, hence it follows that $A\psi$ must be a scalar $a_\psi$ times $\psi$ and $B\eta$ must be a scalar $b_\eta$ times $\eta$. A priori, these scalars may depend on $\psi$ and $\eta$, but we must have $c_{\bf m} = a_\psi \, b_\eta$. By fixing $\psi$ and varying $\eta$ one observes that $b_\eta$ must be a constant over $\eta$ (i.e. does not depend on $\eta$), and by fixing $\eta$ and varying $\psi$ one observes that $a_\psi$ is a constant over $\psi$. Hence the claims in eq.\eqref{eq:separated_identity1} and eq.\eqref{eq:separated_identity2} follow.

\vs

Now, view equations \eqref{eq:separated_identity1} and \eqref{eq:separated_identity2}  as operator identities on $\mathscr{H}_\Gamma$ again, and conjugate by $U^{-1}$. One then obtains the following two identities in $\mathscr{H}_\Gamma$:
\begin{align}
\label{eq:separated_identity3}
& \Phi^{\hbar_{k_1}}(\epsilon_1\, {\bf x}^\hbar_{{\bf m}_1})^{\epsilon_1}
\, \cdots \,
\Phi^{\hbar_{k_r}}(\epsilon_r\, {\bf x}^\hbar_{{\bf m}_r})^{\epsilon_r}
= c_{{\bf m}}' \cdot {\rm Id}, \\
\label{eq:separated_identity4}
& \Phi^{\hbar_{k_1}}(\epsilon_1\, \til{\bf x}^\hbar_{{\bf m}_1})^{-\epsilon_1}
\, \cdots \,
\Phi^{\hbar_{k_r}}(\epsilon_r\, \til{\bf x}^\hbar_{{\bf m}_r} )^{-\epsilon_r}
= c_{{\bf m}}'' \cdot {\rm Id},
\end{align}
with $c_{\bf m} = c'_{\bf m} \, c''_{\bf m}$. From either of these two identities, we will see that one can extract some general versions of operator identities of the non-compact quantum dilogarithm functions. Using these general version identities, we will in turn show that $(c''_{\bf m})^{-1} = c'_{\bf m}$, yielding the desired result $c_{\bf m}=1$. This is done in the next section.

\section{Triviality of the phase constants}
\label{sec:triviality_of_the_phase_constants}

We will assemble the ingredients collected so far and prove some operator identities of the non-compact quantum dilogarithm $\Phi^\hbar$, and finally show the main theorem, Thm.\ref{thm:main}.

\subsection{The $A_1$ identity}

Consider the trivial cluster transformation ${\bf m} = \mu_k \circ \mu_k$, applied to any $\mathcal{D}$-seed $\Gamma = \Gamma^{(0)}$. Recall eq.\eqref{eq:twice_flip_to_prove} of Prop.\ref{prop:the_rank_1_identity_for_intertwiners}. This constant $c_{A_1} \in {\rm U}(1)$ then passes to eq.\eqref{eq:quantum_dilog_factors_only}, for the tropical sign-sequence $\vec{\epsilon}$, which is $(+,-)$, as computed in \S\ref{subsec:triviality_of_linear_parts}; that is, $c_{\bf m} = c_{A_1}$. Equations \eqref{eq:separated_identity3}--\eqref{eq:separated_identity4} read $\Phi^{\hbar_k}({\bf x}^\hbar_{{\bf m}_1}) \, \Phi^{\hbar_k}(-{\bf x}^\hbar_{{\bf m}_2})^{-1} = c'_{\bf m}\cdot {\rm Id}$ and $\Phi^{\hbar_k}(\til{\bf x}^\hbar_{{\bf m}_1})^{-1} \, \Phi^{\hbar_k}(-\til{\bf x}^\hbar_{{\bf m}_2}) = c''_{\bf m} \cdot {\rm Id}$. From eq.\eqref{eq:bf_x_m_i_1}--\eqref{eq:til_bf_x_m_i_1} and Lem.\ref{lem:some_intertwining_equations_for_the_linear_part}, it follows ${\bf x}^\hbar_{{\bf m}_1} = {\bf x}^\hbar_{\Gamma;k}$, ${\bf x}^\hbar_{{\bf m}_2} = - {\bf x}^\hbar_{\Gamma;k}$, $\til{\bf x}^\hbar_{{\bf m}_1} = \til{\bf x}^\hbar_{\Gamma;k}$, and $\til{\bf x}^\hbar_{{\bf m}_2} = - {\bf x}^\hbar_{\Gamma;k}$. Hence it follows $c'_{\bf m} = c''_{\bf m}=1$, and therefore $c_{A_1} = c_{\bf m} = c'_{\bf m} \, c''_{\bf m} = 1$.

\subsection{The argument self-adjoint operators ${\bf x}^\hbar_{{\bf m}_i}$ and $\til{\bf x}^\hbar_{{\bf m}_i}$ for rank 2 identities}
\label{subsec:argument_operators}

Let ${\bf m} = {\bf m}_r \circ \cdots \circ {\bf m}_1$ be a sequence of mutations and seed automorphisms representing rank $2$ identities, namely the ones appearing in (S2) of \S\ref{sec:introduction}, or in Lem.\ref{lem:h_plus_2-gon_relations}. For the remainder of this section, I introduce an abbreviated notation for the operators ${\bf x}^\hbar_{\Gamma^{(\ell)};i}$ and $\til{\bf x}^\hbar_{\Gamma^{(\ell)};i}$ on the Hilbert space $\mathscr{H}_{\Gamma^{(\ell)}}$, for $\ell=0,1,\ldots,r-1$. The Hilbert spaces $\mathscr{H}_{\Gamma^{(\ell)}}$ for different $\ell$'s are related by composition of the unitary operators $({\bf K}')^{(\epsilon_{\ell-1})}_{\Gamma^{(\ell-1)}\mut{k_\ell} \Gamma^{(\ell)}}$, i.e. the linear part of the signed decomposition (Def.\ref{def:signed_two_parts}) of the unitary intertwiner for ${\bf m}_\ell$, where ${\bf m}_ \ell = \mu_{k_\ell}$ for each $\ell=1,2,...,r-1$ and $\vec{\epsilon}=(\epsilon_1,\epsilon_2,...)$ is the tropical sign-sequence computed in \S\ref{subsec:triviality_of_linear_parts}, for the relevant (trivial) cluster transformation ${\bf m}$. When writing equations for operators acting on different $\mathscr{H}_{\Gamma^{(\ell)}}$'s, we will now omit writing the operators $({\bf K}')^{(\epsilon_{\ell-1})}_{\Gamma^{(\ell-1)}\mut{k_\ell} \Gamma^{(\ell)}}$, and express the equal sign `$=$' as the symbol `$\sim$'. Then $\sim$ is an equivalence relation compatible with $=$; in particular, if $A$ and $B$ are operators on $\mathscr{H}_{\Gamma^{(\ell)}}$ for same $\ell$ and satisfies $A\sim B$, then $A=B$. In addition, 
$$
\mbox{the operators ${\bf x}^\hbar_{\Gamma^{(0)};k}$ and  $\til{\bf x}^\hbar_{\Gamma^{(0)};k}$ will be written as ${\bf x}^\hbar_k$ and $\til{\bf x}^\hbar_k$;}
$$
so, whenever there is no $\Gamma^{(\ell)}$ written, we mean we are dealing with $\Gamma=\Gamma^{(0)}$, i.e. operators on $\mathscr{H}_\Gamma$. For example, instead of writing an equation 
$$
{{\bf K}'}^{(\epsilon_1)}_{\Gamma^{(0)}\mut{k_1}\Gamma^{(1)}} \, {{\bf K}'}^{(\epsilon_2)}_{\Gamma^{(1)}\mut{k_2}\Gamma^{(2)}} \, {\bf x}^\hbar_{\Gamma^{(2)};j} \, ({{\bf K}'}^{(\epsilon_2)}_{\Gamma^{(1)}\mut{k_2}\Gamma^{(2)}})^{-1}\, ({{\bf K}'}^{(\epsilon_1)}_{\Gamma^{(0)}\mut{k_1}\Gamma^{(1)}})^{-1}
= {\bf x}^\hbar_{\Gamma^{(0)};i} + {\bf x}^\hbar_{\Gamma^{(0)};j},
$$
we would write
$$
{\bf x}^\hbar_{\Gamma^{(2)};j} \sim {\bf x}^\hbar_{\Gamma^{(0)};i} + {\bf x}^\hbar_{\Gamma^{(0)};j} = {\bf x}^\hbar_i + {\bf x}^\hbar_j
$$
in short. And, under this notation convention, equations \eqref{eq:bf_x_m_i_1}--\eqref{eq:til_bf_x_m_i_1} can be written as
\begin{align}
\label{eq:argument_operators_short}
{\bf x}^\hbar_{{\bf m}_\ell} \sim {\bf x}^\hbar_{\Gamma^{(\ell-1)};k_\ell},
\qquad
\til{\bf x}^\hbar_{{\bf m}_\ell} \sim \til{\bf x}^\hbar_{\Gamma^{(\ell-1)};k_\ell},
\end{align}
and the linear-part mutation formula in Lem.\ref{lem:some_intertwining_equations_for_the_linear_part} reads
\begin{align}
\label{eq:bf_K_prime_conjugation_on_bf_x_i_new}
& {\bf x}^\hbar_{\Gamma^{(\ell)};k} \sim \left\{
\begin{array}{ll}
{\bf x}^\hbar_{\Gamma^{(\ell-1)};k} + [\epsilon_\ell \, \varepsilon^{(\ell-1)}_{k,k_\ell}]_+ \, {\bf x}^\hbar_{\Gamma^{(\ell-1)};k_\ell}, & \mbox{if $k\neq k_\ell$}, \\
- {\bf x}^\hbar_{\Gamma^{(\ell-1)};k_\ell}, & \mbox{if $k=k_\ell$,}
\end{array}
\right. \\
\label{eq:bf_K_prime_conjugation_on_til_bf_x_i_new}
& \til{\bf x}^\hbar_{\Gamma^{(\ell)};k} \sim \left\{
\begin{array}{ll}
\til{\bf x}^\hbar_{\Gamma^{(\ell-1)};k} + [\epsilon_\ell \, \varepsilon^{(\ell-1)}_{k,k_\ell}]_+ \, \til{\bf x}^\hbar_{\Gamma^{(\ell-1)};k_\ell}, & \mbox{if $k\neq k_\ell$}, \\
- \til{\bf x}^\hbar_{\Gamma^{(\ell-1)};k_\ell}, & \mbox{if $k=k_\ell$.}
\end{array}
\right.
\end{align}

\vs

Now we prove the operator-separation property of Def.\ref{def:operator-separation_property}, for the cases (S2-2), (S2-3), and (S2-4) in \S\ref{sec:introduction}, or equivalently, for the cases $p=1,2$, or $3$ in Lem.\ref{lem:h_plus_2-gon_relations}; we will prove it by finding a unitary operator $U$ explicitly. The sequence ${\bf m}$ is one of $P_{(ij)} \circ \mu_i \circ \mu_j \circ \mu_i \circ \mu_j \circ \mu_i$, \quad $\mu_j \circ \mu_i \circ \mu_j \circ \mu_i \circ \mu_j \circ \mu_i$, \quad or \quad $\mu_j \circ \mu_i \circ \mu_j \circ \mu_i \circ \mu_j \circ \mu_i \circ \mu_j \circ \mu_i$. In particular, the indices $k_\ell$ for the mutations and in the argument operators in eq.\eqref{eq:argument_operators_short} are $i$ or $j$. By applying the above linear mutation formula (Lem.\ref{lem:some_intertwining_equations_for_the_linear_part}) repeatedly, we can observe that ${\bf x}^\hbar_{{\bf m}_\ell}$ is a $\mathbb{Z}$-linear combination of ${\bf x}^\hbar_i$ and ${\bf x}^\hbar_j$, while $\til{\bf x}^\hbar_{{\bf m}_\ell}$ is a $\mathbb{Z}$-linear combination of $\til{\bf x}^\hbar_i$ and $\til{\bf x}^\hbar_j$. So, in order to prove that ${\bf m}$ satisfies the operator-separation property, it suffices to prove that the operators ${\bf x}^\hbar_i$ and ${\bf x}^\hbar_j$ are strongly separated (Def.\ref{def:operator-separation_property}) with $\til{\bf x}^\hbar_i$ and $\til{\bf x}^\hbar_j$. Note  that $i\neq j$, and that $\varepsilon_{ij}$ and $\varepsilon_{ji}$ are nonzero; also, $n\ge 2$. Recall from Def.\ref{def:old_representation} that
$$
\textstyle {\bf x}^\hbar_i = \frac{d_i^{-1}}{2} {\bf p}^\hbar_i - \sum_{k=1}^n \varepsilon_{ik} {\bf q}^\hbar_k, \quad
\til{\bf x}^\hbar_i = \frac{d_i^{-1}}{2} {\bf p}^\hbar_i + \sum_{k=1}^n \varepsilon_{ik} {\bf q}^\hbar_k, \quad
{\bf x}^\hbar_j = \frac{d_j^{-1}}{2} {\bf p}^\hbar_j - \sum_{k=1}^n \varepsilon_{jk} {\bf q}^\hbar_k, \quad
\til{\bf x}^\hbar_j = \frac{d_j^{-1}}{2} {\bf p}^\hbar_j + \sum_{k=1}^n \varepsilon_{jk} {\bf q}^\hbar_k.
$$
As $\varepsilon_{ii}=0=\varepsilon_{jj}$, we may replace the sum sign $\sum_{k=1}^n$ with $\sum_{k\neq i}$ in the first two equations, and with $\sum_{k\neq j}$ in the latter two equations. From now on, we will conjugate the above self-adjoint operators by unitary operators $U_1, U_2,\ldots$ step by step, where conjugation by each $U_\ell$ takes an $\mathbb{R}$-linear combination of $\{ {\bf p}^\hbar_k, {\bf q}^\hbar_k : k=1,2,\ldots,n\}$ to an $\mathbb{R}$-linear combination of $\{ {\bf p}^\hbar_k, {\bf q}^\hbar_k : k=1,2,\ldots,n\}$ again. Let $U_1 := \mathcal{F}_i$ be the Fourier transform in the variable $a_i$ of $\mathscr{H}_\Gamma = L^2(\mathbb{R}^n, da_1 \cdots da_n)$, as defined in eq.\eqref{eq:F_i}. Conjugate the four operators ${\bf x}^\hbar_i,\til{\bf x}^\hbar_i,{\bf x}^\hbar_j,\til{\bf x}^\hbar_j$ by $U_1$; by eq.\eqref{eq:F_i_conjugation}, the results are indeed $\mathbb{R}$-linear combinations of ${\bf p}^\hbar_k$'s and ${\bf q}^\hbar_k$'s. In particular, $U_1 {\bf x}^\hbar_i U_1^{-1} = -2\pi^2 \hbar d_i^{-1} {\bf q}^\hbar_i - \sum_{k\neq i} \varepsilon_{ik} {\bf q}^\hbar_k$. Now let $U_2 := {\bf S}_{{\bf c}}$, with ${\bf c} = (c_{k\ell})_{k,\ell} \in {\rm SL}_\pm(n,\mathbb{R})$ given by
$$
\textstyle
c_{kk}=1, \quad \forall k, \qquad
c_{ki} = \frac{-\varepsilon_{ik}}{2\pi^2 \hbar d_i^{-1}}, \quad \forall k\neq i, \qquad
c_{k\ell}=0 \quad \mbox{otherwise}.
$$
Then by Lem.\ref{lem:bf_S_conjugation_on_bf_p_and_bf_q} it follows $U_2 U_1 {\bf x}^\hbar_i U_1^{-1} U_2^{-1} = - 2\pi^2\hbar d_i^{-1} {\bf q}^\hbar_i$. Let $U_3 = {\bf S}_{{\bf c}'}$ with ${\bf c}' \in {\rm SL}_\pm(n,\mathbb{R})$ being the diagonal matrix whose $i$-th diagonal entry is $(-2\pi^2\hbar d_i^{-1})^{-1}$ and $j$-th diagonal entry is $-2\pi^2\hbar d_i^{-1}$, while other diagonal entries are $1$. Then by Lem.\ref{lem:bf_S_conjugation_on_bf_p_and_bf_q} it follows $U_3U_2U_1 {\bf x}^\hbar_i U_1^{-1}U_2^{-1} U_3^{-1} = {\bf q}^\hbar_i$.  Meanwhile, note from eq.\eqref{eq:F_i_conjugation} that $U_1 {\bf x}^\hbar_j U_1^{-1} = \frac{d_j^{-1}}{2} {\bf p}^\hbar_j - \frac{\varepsilon_{ji}}{4\pi^2 \hbar} {\bf p}^\hbar_i - \sum_{k\neq j,i} \varepsilon_{jk} {\bf q}^\hbar_k$.   Observe that the entries of the inverse ${\bf c}^{-1} = (c^{k\ell})_{k,\ell}$ of ${\bf c}$ is given by $c_{kk}=1$, $\forall k$, $c_{ki} = \frac{\varepsilon_{ik}}{2\pi^2 \hbar d_i^{-1}}$, $\forall k \neq i$, and $c_{k\ell} =0$ otherwise. So, from Lem.\ref{lem:bf_S_conjugation_on_bf_p_and_bf_q} it follows that $U_2 U_1 {\bf x}^\hbar_j U_1^{-1} U_2^{-1} = \frac{d_j^{-1}}{2} ({\bf p}^\hbar_j + \frac{\varepsilon_{ij}}{2\pi^2\hbar d_i^{-1}} {\bf p}^\hbar_i) - \frac{\varepsilon_{ji}}{4\pi^2\hbar} {\bf p}^\hbar_i - \sum_{k\neq j,i} \varepsilon_{jk} {\bf q}^\hbar_k = \frac{d_j^{-1}}{2} {\bf p}^\hbar_j - \frac{\varepsilon_{ji}}{2\pi^2\hbar} {\bf p}^\hbar_i - \sum_{k\neq j,i} \varepsilon_{jk} {\bf q}^\hbar_k$, where we used $\frac{d_j^{-1}}{d_i^{-1}} \varepsilon_{ij} = \frac{\wh{\varepsilon}_{ij}}{d_i^{-1}} = \frac{-\wh{\varepsilon}_{ji}}{d_i^{-1}} = - \varepsilon_{ji}$. Then note $U_3 U_2 U_1 {\bf x}^\hbar_j U_1^{-1} U_2^{-1} U_3^{-1} = \frac{d_j^{-1}}{-4\pi^2\hbar d_i^{-1}} {\bf p}^\hbar_j + \varepsilon_{ji} d_i^{-1} {\bf p}^\hbar_i - \sum_{k\neq j,i} \varepsilon_{jk} {\bf q}^\hbar_k$. Now let $U_4 = \prod_{k\neq i,j} \mathcal{F}_k$. Then $U_4U_3U_2U_1 {\bf x}^\hbar_i U_1^{-1}U_2^{-1}U_3^{-1}U_4^{-1} = {\bf q}^\hbar_i$, while $U_4U_3U_2U_1 {\bf x}^\hbar_j U_1^{-1} U_2^{-1} U_3^{-1} U_4^{-1} = \frac{d_j^{-1}}{-4\pi^2 \hbar d_i^{-1}} {\bf p}^\hbar_j + \varepsilon_{ji} d_i^{-1} {\bf p}^\hbar_i - \sum_{k\neq j,i} \frac{\varepsilon_{jk}}{4\pi^2 \hbar} {\bf p}^\hbar_k$. Let $U_5 = {\bf S}_{{\bf c}''}$ where ${\bf c}''$ equals the identity matrix except at the $i$-th row, with the $(i,i)$-th entry being $1$, the $(i,k)$-th entry being $\frac{-1}{\varepsilon_{ji} d_i^{-1}} \cdot \frac{\varepsilon_{jk}}{4\pi^2 \hbar }$ for every $k\neq i,j$, and $(i,j)$-th entry being $\frac{-1}{\varepsilon_{ji} d_i^{-1}} \cdot \frac{d_j^{-1}}{4\pi^2 \hbar d_i^{-1}}$; this is possible because $\varepsilon_{ji} \neq 0$. Then $({\bf c}'')^{-1}$ equals the identity matrix except at the $i$-th row, with the $(i,i)$-th entry being $1$, and the $(i,k)$-th entry being $-1$ times the $(i,k)$-th entry of ${\bf c}''$ for each $k\neq i$. By Lem.\ref{lem:bf_S_conjugation_on_bf_p_and_bf_q} it follows that $U_5 U_4 U_3 U_2 U_1 {\bf x}^\hbar_i U_1^{-1} U_2^{-1} U_3^{-1} U_4^{-1} U_5^{-1} = {\bf q}^\hbar_i $ and $U_5 U_4 U_3 U_2 U_1 {\bf x}^\hbar_j U_1^{-1} U_2^{-1} U_3^{-1} U_4^{-1} U_5^{-1} = \varepsilon_{ji} d_i^{-1} {\bf p}^\hbar_i$. Let $U := U_5 U_4 U_3 U_2 U_1$. Conjugation by $U$ takes an $\mathbb{R}$-linear combination of ${\bf p}^\hbar_k$'s and ${\bf q}^\hbar_k$'s to an $\mathbb{R}$-linear combination of ${\bf p}^\hbar_k$'s and ${\bf q}^\hbar_k$'s, so one can write $U \til{\bf x}^\hbar_i U^{-1} = \sum_{k=1}^n (\alpha_k {\bf p}^\hbar_k + \beta_k {\bf q}^\hbar_k)$ and $U \til{\bf x}^\hbar_j U^{-1} = \sum_{k=1}^n (\alpha_k' {\bf p}^\hbar_k + \beta_k' {\bf q}^\hbar_k)$ for some real numbers $\alpha_k,\beta_k,\alpha_k',\beta_k'$. Note ${\bf x}^\hbar_i$ and ${\bf x}^\hbar_j$ strongly commute with the tilde operators $\til{\bf x}^\hbar_i$ and $\til{\bf x}^\hbar_j$; hence $U{\bf x}^\hbar_iU^{-1}$ and $U{\bf x}^\hbar_jU^{-1}$ strongly commute with $U\til{\bf x}^\hbar_iU^{-1}$ and $U\til{\bf x}^\hbar_j U^{-1}$. Computing the commutators using eq.\eqref{eq:Heisenberg_relations_of_p_and_q} and bilinearity, we have $[U{\bf x}^\hbar_i U^{-1}, \, U\til{\bf x}^\hbar_i U^{-1}] = [{\bf q}^\hbar_i, \sum_{k=1}^n (\alpha_k {\bf p}^\hbar_k + \beta_k {\bf q}^\hbar_k)] = - 2\pi {\rm i} \alpha_i \hbar \cdot {\rm Id}$, and similarly $[U{\bf x}^\hbar_i U^{-1}, U\til{\bf x}_j^\hbar U^{-1}] = -2\pi {\rm i} \alpha'_i \hbar \cdot {\rm Id}$, \, $[U{\bf x}^\hbar_j U^{-1}, U\til{\bf x}^\hbar_i U^{-1}] = \varepsilon_{ij} d_i^{-1} \, 2\pi {\rm i} \beta_i \hbar \cdot {\rm Id}$, and $[U{\bf x}^\hbar_j U^{-1}, U\til{\bf x}^\hbar_j U^{-1}] = \varepsilon_{ij} d_i^{-1} \, 2\pi {\rm i} \beta_i' \hbar \cdot {\rm Id}$. Since these commutators must be zero, it follows $\alpha_i=\alpha_i' = \beta_i = \beta'_i=0$, which shows that ${\bf x}^\hbar_i$ and ${\bf x}^\hbar_j$ are indeed strongly separated (Def.\ref{def:operator-separation_property}) with $\til{\bf x}^\hbar_i$ and $\til{\bf x}^\hbar_j$, via the unitary operator $U$. We summarize this as:
\begin{lemma}
\label{lem:operator-separation_property}
The sequence ${\bf m}$ corresponding to (S2-2), (S2-3), or (S2-4) of \S\ref{sec:introduction} satisfies the operator-separation property of Def.\ref{def:operator-separation_property}. \qed
\end{lemma}
Hence, for the cases (S2-2), (S2-3), and (S2-4), we can apply the arguments and results of \S\ref{subsec:from_double_to_single}.

\subsection{The $A_1 \times A_1$ identity}
\label{subsec:A1A1_identity}

Consider the trivial cluster transformation ${\bf m} = \mu_j \circ \mu_i \circ \mu_j \circ \mu_i$, applied to any $\mathcal{D}$-seed $\Gamma = \Gamma^{(0)}$ satisfying  $\varepsilon_{ij}=0$, with $i\neq j$. The constant $c_{A_1\times A_1} \in {\rm U}(1)$ appearing in eq.\eqref{eq:A1_times_A1_to_prove} of Prop.\ref{prop:rank_2_identities_for_intertwiners}.(2-1) passes to eq.\eqref{eq:quantum_dilog_factors_only}, for the tropical sign-sequence $\vec{\epsilon}$, which is $(+,+,-,-)$, as computed in \S\ref{subsec:triviality_of_linear_parts}; that is, $c_{\bf m} = c_{A_1\times A_1}$. The proof of the separation property which I gave in the last subsection does not work, so at the moment, we only have the non-separated identity in eq.\eqref{eq:quantum_dilog_factors_only}, but not necessarily the separated ones in eq.\eqref{eq:separated_identity3}--\eqref{eq:separated_identity4}. Note that eq.\eqref{eq:quantum_dilog_factors_only} reads
\begin{align}
\label{eq:quantum_dilog_factors_only_A1A1}
\Phi^{\hbar_i}({\bf x}^\hbar_{{\bf m}_1}) \Phi^{\hbar_i}(\til{\bf x}^\hbar_{{\bf m}_1})^{-1}
\Phi^{\hbar_j}({\bf x}^\hbar_{{\bf m}_2}) \Phi^{\hbar_j}(\til{\bf x}^\hbar_{{\bf m}_2})^{-1}
\Phi^{\hbar_i}(-{\bf x}^\hbar_{{\bf m}_3})^{-1} \Phi^{\hbar_i}(-\til{\bf x}^\hbar_{{\bf m}_3})
\Phi^{\hbar_j}(-{\bf x}^\hbar_{{\bf m}_4})^{-1}  \Phi^{\hbar_j}(-\til{\bf x}^\hbar_{{\bf m}_4})
= c_{\bf m} \cdot {\rm Id}.
\end{align}

\vs

By the mutation formula in eq.\eqref{eq:varepsilon_prime_formula} for the exchange matrices, one has $\varepsilon^{(\ell)}_{ij} = \varepsilon^{(\ell)}_{ji}=0$ for each $\ell = 0,1,2,3$. Keeping in mind $(k_1,k_2,k_3,k_4)=(i,j,i,j)$ and the equivalence $\sim$ of operators introduced in \S\ref{subsec:argument_operators}, from equations \eqref{eq:argument_operators_short}, \eqref{eq:bf_K_prime_conjugation_on_bf_x_i_new},  and \eqref{eq:bf_K_prime_conjugation_on_til_bf_x_i_new}, we compute the operators ${\bf x}^\hbar_{{\bf m}_\ell}$'s as follows:
\begin{align*}
& {\bf x}^\hbar_{{\bf m}_1} = {\bf x}^\hbar_{\Gamma^{(0)};i} = {\bf x}^\hbar_i, \qquad
{\bf x}^\hbar_{{\bf m}_2} \sim 
{\bf x}^\hbar_{\Gamma^{(1)};j} \sim {\bf x}^\hbar_{\Gamma^{(0)};j} + \underset{=0}{\ul{[\varepsilon^{(0)}_{ji}]_+}} \, {\bf x}^\hbar_{\Gamma^{(0)};i} 
= {\bf x}^\hbar_j, \\
& {\bf x}^\hbar_{{\bf m}_3} \sim 
{\bf x}^\hbar_{\Gamma^{(2)};i} 
\sim 
{\bf x}^\hbar_{\Gamma^{(1)};i} + \underset{=0}{\ul{ [\varepsilon^{(1)}_{ij}]_+ }} \, {\bf x}^\hbar_{\Gamma^{(1)};j} 
= {\bf x}^\hbar_{\Gamma^{(1)};i} \sim - {\bf x}^\hbar_{\Gamma^{(0)};i} = - {\bf x}^\hbar_i, \\
& {\bf x}^\hbar_{{\bf m}_4} \sim  
{\bf x}^\hbar_{\Gamma^{(3)};j} 
\sim {\bf x}^\hbar_{\Gamma^{(2)};j} + \underset{=0}{\ul{ [-\varepsilon^{(2)}_{ji}]_+ } } \, {\bf x}^\hbar_{\Gamma^{(2)};i} 
= {\bf x}^\hbar_{\Gamma^{(2)};j} \sim -{\bf x}^\hbar_{\Gamma^{(1)};j} 
\sim - {\bf x}^\hbar_j.
\end{align*}
For the tilde operators, note that in Lem.\ref{lem:some_intertwining_equations_for_the_linear_part}, the coefficients appearing in the result of the conjugation action of ${{\bf K}'}^{(\epsilon)}$ on the tilde operators are the same as those of the non-tilde operators. So by equations \eqref{eq:argument_operators_short}, \eqref{eq:bf_K_prime_conjugation_on_bf_x_i_new},  and \eqref{eq:bf_K_prime_conjugation_on_til_bf_x_i_new}, and by the above computation of ${\bf x}^\hbar_{{\bf m}_\ell}$'s, one obtains
$$
\til{\bf x}^\hbar_{{\bf m}_1} = \til{\bf x}^\hbar_i, \quad
\til{\bf x}^\hbar_{{\bf m}_2} = \til{\bf x}^\hbar_j, \quad
\til{\bf x}^\hbar_{{\bf m}_3} = - \til{\bf x}^\hbar_i, \quad
\til{\bf x}^\hbar_{{\bf m}_4} = - \til{\bf x}^\hbar_j,
$$
so eq.\eqref{eq:quantum_dilog_factors_only_A1A1} reads
$$
\Phi^{\hbar_i}({\bf x}^\hbar_i) \Phi^{\hbar_i}(\til{\bf x}^\hbar_i)^{-1}
\Phi^{\hbar_j}({\bf x}^\hbar_j) \Phi^{\hbar_j}(\til{\bf x}^\hbar_j)^{-1}
\Phi^{\hbar_i}({\bf x}^\hbar_i)^{-1} \Phi^{\hbar_i}(\til{\bf x}^\hbar_i)
\Phi^{\hbar_j}({\bf x}^\hbar_j)^{-1} \Phi^{\hbar_j}(\til{\bf x}^\hbar_j)
= c_{\bf m} \cdot {\rm Id}.
$$
Notice that all four operators ${\bf x}^\hbar_i = {\bf x}^\hbar_{\Gamma;i}$, ${\bf x}^\hbar_j = {\bf x}^\hbar_{\Gamma;j}$, $\til{\bf x}^\hbar_i = \til{\bf x}^\hbar_{\Gamma;i}$, and $\til{\bf x}^\hbar_j = \til{\bf x}^\hbar_{\Gamma;j}$ strongly commute with one another, as they satisfy the Weyl relations for the zero Heisenberg relations. So, in the left-hand-side of the above eq.\eqref{eq:quantum_dilog_factors_only_A1A1}, each factor commutes with each factor. Rearranging the factors, we observe that each factor in the left-hand-side gets canceled with some other factor. Hence it follows $c_{\bf m}=1$, and therefore $c_{A_1\times A_1} = c_{\bf m} =1$.

\subsection{The $A_2$ identity}
\label{subsec:A2}

Consider the trivial cluster transformation ${\bf m} = P_{(ij)} \circ \mu_i \circ \mu_j \circ \mu_i \circ \mu_j \circ \mu_i$, applied to any $\mathcal{D}$-seed $\Gamma = \Gamma^{(0)}$ satisfying $\varepsilon_{ij} = - \varepsilon_{ji} \in \{1,-1\}$. In particular one has $d_i=d_j$, so $\hbar_i = \hbar_j$. Let us first deal with the case when
$$
\varepsilon_{ij} = - \varepsilon_{ji} = 1.
$$
As computed in \S\ref{subsec:triviality_of_linear_parts}, the tropical sign-sequence is $\vec{\epsilon} = (+,+,-,-,-,0)$. The constant $c_{A_2} \in {\rm U}(1)$ in eq.\eqref{eq:A2_to_prove} of Prop.\ref{prop:rank_2_identities_for_intertwiners}.(2-2) passes to eq.\eqref{eq:quantum_dilog_factors_only}; that is, $c_{\bf m} = c_{A_2}$. Thanks to Lem.\ref{lem:operator-separation_property}, we get the separated identities obtained in \S\ref{subsec:from_double_to_single}, i.e. equations \eqref{eq:separated_identity3}--\eqref{eq:separated_identity4}, which read
\begin{align}
\label{eq:separated_identity3_A2_case1}
& \Phi^{\hbar_i}({\bf x}^\hbar_{{\bf m}_1}) \, \Phi^{\hbar_j}({\bf x}^\hbar_{{\bf m}_2}) \, \Phi^{\hbar_i}(-{\bf x}^\hbar_{{\bf m}_3})^{-1} \, 
\Phi^{\hbar_j}(-{\bf x}^\hbar_{{\bf m}_4})^{-1} \, 
\Phi^{\hbar_i}(-{\bf x}^\hbar_{{\bf m}_5})^{-1} = c'_{\bf m} \cdot {\rm Id},  \\
\label{eq:separated_identity4_A2_case1}
& \Phi^{\hbar_i} (\til{\bf x}^\hbar_{{\bf m}_1})^{-1} \, \Phi^{\hbar_j}(\til{\bf x}^\hbar_{{\bf m}_2})^{-1} \, \Phi^{\hbar_i}(-\til{\bf x}^\hbar_{{\bf m}_3}) \, 
\Phi^{\hbar_j}(-\til{\bf x}^\hbar_{{\bf m}_4}) \, 
\Phi^{\hbar_i}(-\til{\bf x}^\hbar_{{\bf m}_5}) = c''_{\bf m} \cdot {\rm Id}.
\end{align}
By the mutation formula in eq.\eqref{eq:varepsilon_prime_formula} for the exchange matrices, one has $\varepsilon^{(\ell)}_{ij} = - \varepsilon^{(\ell)}_{ji} = 1$ for $\ell = 0,2$ and $\varepsilon^{(\ell)}_{ij} = - \varepsilon^{(\ell)}_{ji} = -1$ for $\ell = 1,3$. Keeping in mind $(k_1,k_2,k_3,k_4,k_5) = (i,j,i,j,i)$, from equations \eqref{eq:argument_operators_short} and  \eqref{eq:bf_K_prime_conjugation_on_bf_x_i_new}, and under the notation convention introduced in \S\ref{subsec:argument_operators}, we get ${\bf x}^\hbar_{{\bf m}_1} = {\bf x}^\hbar_{\Gamma^{(0)};i} = {\bf x}^\hbar_i$, and 
\begin{align*}
\begin{array}{ll}
\qquad\quad {\bf x}^\hbar_{\Gamma^{(1)};i} \sim - {\bf x}^\hbar_{\Gamma^{(0)};i} = - {\bf x}^\hbar_i, & 
{\bf x}^\hbar_{{\bf m}_2} \sim {\bf x}^\hbar_{\Gamma^{(1)};j} 
\sim {\bf x}^\hbar_{\Gamma^{(0)};j} + \underset{=0}{\ul{ [\varepsilon^{(0)}_{ji}]_+}} \, {\bf x}^\hbar_{\Gamma^{(0)};i} = {\bf x}^\hbar_j,  \\
{\bf x}^\hbar_{{\bf m}_3} \sim 
{\bf x}^\hbar_{\Gamma^{(2)};i} 
\sim {\bf x}^\hbar_{\Gamma^{(1)};i}+\underset{=0}{\ul{ [\varepsilon^{(1)}_{ij}]_+ }}\, {\bf x}^\hbar_{\Gamma^{(1)};j} 
\sim - {\bf x}^\hbar_i, & \qquad\quad 
{\bf x}^\hbar_{\Gamma^{(2)};j} \sim - {\bf x}^\hbar_{\Gamma^{(1)};j} \sim - {\bf x}^\hbar_j, \\
\qquad\quad  {\bf x}^\hbar_{\Gamma^{(3)};i} \sim - {\bf x}^\hbar_{\Gamma^{(2)};i} \sim {\bf x}^\hbar_i, &
{\bf x}^\hbar_{{\bf m}_4} \sim 
{\bf x}^\hbar_{\Gamma^{(3)};j} 
\sim {\bf x}^\hbar_{\Gamma^{(2)};j}+\underset{=1}{\ul{ [-\varepsilon^{(2)}_{ji}]_+}} \, {\bf x}^\hbar_{\Gamma^{(2)};i} \sim  - {\bf x}^\hbar_i - {\bf x}^\hbar_j,
\\
{\bf x}^\hbar_{{\bf m}_5} \sim 
{\bf x}^\hbar_{\Gamma^{(4)};i} 
\sim {\bf x}^\hbar_{\Gamma^{(3)};i} + \underset{=1}{\ul{ [-\varepsilon^{(3)}_{ij}]_+}} \, {\bf x}^\hbar_{\Gamma^{(3)};j} 
\sim  -{\bf x}^\hbar_j, &
\end{array}
\end{align*}
where for each line, we used the previous line. In view of the comparison between equations \eqref{eq:bf_K_prime_conjugation_on_bf_x_i_new} and \eqref{eq:bf_K_prime_conjugation_on_til_bf_x_i_new}, the tilde version is parallel. Namely, we have
$$
\til{\bf x}^\hbar_{{\bf m}_1} = \til{\bf x}^\hbar_i, \qquad
\til{\bf x}^\hbar_{{\bf m}_2} = \til{\bf x}^\hbar_j, \qquad
\til{\bf x}^\hbar_{{\bf m}_3} = - \til{\bf x}^\hbar_i, \qquad
\til{\bf x}^\hbar_{{\bf m}_4} = - \til{\bf x}^\hbar_i - \til{\bf x}^\hbar_j, \qquad
\til{\bf x}^\hbar_{{\bf m}_5} = - \til{\bf x}^\hbar_j.
$$
Thus, from eq.\eqref{eq:separated_identity3_A2_case1}--\eqref{eq:separated_identity4_A2_case1} we get
\begin{align}
\label{eq:cluster_pentagon_identity1}
& \Phi^{\hbar_i}({\bf x}^\hbar_i) \, \Phi^{\hbar_j}({\bf x}^\hbar_j) \, \Phi^{\hbar_i}({\bf x}^\hbar_i)^{-1} \, \Phi^{\hbar_j}({\bf x}^\hbar_i + {\bf x}^\hbar_j)^{-1} \, \Phi^{\hbar_i}({\bf x}^\hbar_j)^{-1} = c_{{\bf m}}' \, \cdot {\rm Id}, \\
\label{eq:cluster_pentagon_identity2}
& \Phi^{\hbar_i}(\til{\bf x}^\hbar_i)^{-1} \, \Phi^{\hbar_j}(\til{\bf x}^\hbar_j)^{-1} \, \Phi^{\hbar_i}(\til{\bf x}^\hbar_i) \, \Phi^{\hbar_j}(\til{\bf x}^\hbar_i + \til{\bf x}^\hbar_j) \, \Phi^{\hbar_i}(\til{\bf x}^\hbar_j) = c_{{\bf m}}'' \, \cdot {\rm Id}.
\end{align}
We may now extract the following general operator identity for the non-compact quantum dilogarithm.
\begin{proposition}[pentagon identity for the non-compact quantum dilogarithm; \cite{FKV01} \cite{W} \cite{G08} ]
\label{prop:pentagon_identity_for_ncQD}
There exists a constant $c^\hbar_{\rm pent} \in {\rm U}(1)$ depending on $\hbar$ s.t. the equality of unitary operators
\begin{align}
\label{eq:general_pentagon_identity}
\Phi^\hbar({\bf P}^\hbar) \, \Phi^\hbar({\bf Q}^\hbar) \, \Phi^\hbar({\bf P}^\hbar)^{-1} \, \Phi^\hbar({\bf P}^\hbar+{\bf Q}^\hbar)^{-1} \, \Phi^\hbar({\bf Q}^\hbar)^{-1} = c_{\rm pent}^\hbar \, \cdot {\rm Id}
\end{align}
holds whenever ${\bf P}^\hbar, {\bf Q}^\hbar$ are self-adjoint operators on a separable Hilbert space that satisfy the Weyl relations for the Heisenberg relation
$$
[{\bf P}^\hbar, {\bf Q}^\hbar] = 2\pi {\rm i} \hbar \cdot {\rm Id}.
$$
\end{proposition}

{\it Proof.} Let us use our particular ${\bf m}$ and $\Gamma$. As found out in \S\ref{subsec:argument_operators} in the proof of Lem.\ref{lem:operator-separation_property}, there is a unitary operator $U : \mathscr{H}_\Gamma \to \mathscr{H}_\Gamma$ such that $U\, {\bf x}^\hbar_i \, U^{-1} = {\bf q}^\hbar_i$ and $U \, {\bf x}^\hbar_j \, U^{-1} = \varepsilon_{ji} d_i^{-1} {\bf p}^\hbar_i = -d_i^{-1} {\bf p}^\hbar_i$. Hence by conjugating by $U$ on both sides of eq.\eqref{eq:cluster_pentagon_identity1} we get
\begin{align}
\label{eq:cluster_pentagon_identity1_conjugated}
\textstyle 
\Phi^{\hbar_i}({\bf q}^\hbar_i) \, 
\Phi^{\hbar_i}(-d_i^{-1}{\bf p}^\hbar_i) \,  
\Phi^{\hbar_i}({\bf q}^\hbar_i)^{-1} \,
\Phi^{\hbar_i}({\bf q}^\hbar_i - d_i^{-1}{\bf p}^\hbar_i)^{-1} \,
\Phi^{\hbar_i}(-d_i^{-1}{\bf p}^\hbar_i)^{-1} = c'_{\bf m} \cdot {\rm Id},
\end{align}
where we used $\hbar_i = \hbar_j$. The left-hand-side is of the form $* \otimes {\rm Id}$ on $L^2(\mathbb{R}, da_i) \otimes L^2(\mathbb{R}^{n-1}, da_1 \cdots da_{i-1} \, da_{i+1} \cdots da_n) = L^2(\mathbb{R}^n,da_1 \cdots da_n)$, involving only the first tensor factor. Hence this identity can be viewed as holding on the Hilbert space $L^2(\mathbb{R},da_i)$, by a similar argument as in \S\ref{subsec:from_double_to_single}. Now, if ${\bf P}^{\hbar_i}, {\bf Q}^{\hbar_i}$ are any self-adjoint operators on a separable Hilbert space $\mathscr{H}'$ satisfying the Weyl relations for the Heisenberg relation $[{\bf P}^{\hbar_i}, {\bf Q}^{\hbar_i}] = 2\pi {\rm i} \hbar_i \cdot {\rm Id}$, then by the Stone-von Neumann theorem (Prop.\ref{prop:SvN}), one has $\mathscr{H}' = \bigoplus_{\ell=1}^N \mathscr{H}_\ell$, and for each $\ell$ there exists a unitary operator $U_\ell : \mathscr{H}_\ell \to L^2(\mathbb{R},da_i)$ such that $U_\ell \, ({\bf P}^{\hbar_i} \restriction \mathscr{H}_\ell) \, (U_\ell)^{-1} = {\bf p}^{\hbar_i}_i$ and $U_\ell \, ({\bf Q}^{\hbar_i} \restriction \mathscr{H}_\ell) \, (U_\ell)^{-1} = {\bf q}^{\hbar_i}_i$. Since ${\bf q}^\hbar_i$ and $-d_i^{-1} {\bf p}^\hbar_i$ are self-adjoint operators on a separable Hilbert space $L^2(\mathbb{R},da_i)$ satisfying the Weyl relations for $[{\bf q}^\hbar_i, - d_i^{-1} {\bf p}^\hbar_i] = 2\pi {\rm i} \hbar d_i^{-1} \cdot {\rm Id} = 2\pi{\rm i} \hbar_i \cdot {\rm Id}$, and since $L^2(\mathbb{R},da_i)$ is irreducible for these self-adjoint operators in the sense as in Prop.\ref{prop:SvN}.(4) as is well-known (see e.g. \cite{Hall}), there exists a unitary operator $U' : L^2(\mathbb{R},da_i) \to L^2(\mathbb{R},da_i)$ such that $U' {\bf q}^\hbar_i (U')^{-1} = {\bf p}^{\hbar_i}_i$ and $U' (-d_i^{-1}{\bf p}^\hbar_i) (U')^{-1} = {\bf q}^{\hbar_i}_i$. Thus it follows that the unitary operator $U'_\ell := (U')^{-1} U_\ell : \mathscr{H}_\ell \to L^2(\mathbb{R},da_i)$ satisfies $U'_\ell ({\bf P}^{\hbar_i} \restriction \mathscr{H}_\ell) (U'_\ell)^{-1} = {\bf q}^\hbar_i$ and $U'_\ell ({\bf Q}^{\hbar_i} \restriction\mathscr{H}_\ell) (U'_\ell)^{-1} = -d_i^{-1} {\bf p}^\hbar_i$. Note $\Phi^{\hbar_i}({\bf P}^{\hbar_i})$, $\Phi^{\hbar_i}({\bf Q}^{\hbar_i})$, and $\Phi^{\hbar_i}({\bf P}^{\hbar_i}+{\bf Q}^{\hbar_i})$ all preserve each $\mathscr{H}_\ell$, and $\Phi^{\hbar_i}({\bf P}^{\hbar_i})\restriction \mathscr{H}_\ell = \Phi^{\hbar_i}({\bf P}^{\hbar_i}\restriction \mathscr{H}_\ell)$, $\Phi^{\hbar_i}({\bf Q}^{\hbar_i})\restriction \mathscr{H}_\ell = \Phi^{\hbar_i}({\bf Q}^{\hbar_i}\restriction \mathscr{H}_\ell)$, and $\Phi^{\hbar_i}({\bf P}^{\hbar_i}+{\bf Q}^{\hbar_i})\restriction \mathscr{H}_\ell = \Phi^{\hbar_i}(({\bf P}^{\hbar_i}\restriction \mathscr{H}_\ell) +({\bf Q}^{\hbar_i}\restriction \mathscr{H}_\ell) )$. Therefore, 
\begin{align*}
\textstyle 
& \left( \Phi^{\hbar_i}({\bf P}^{\hbar_i}) \, 
\Phi^{\hbar_i}({\bf Q}^{\hbar_i}) \,  
\Phi^{\hbar_i}({\bf P}^{\hbar_i})^{-1} \,
\Phi^{\hbar_i}({\bf P}^{\hbar_i} + {\bf Q}^{\hbar_i})^{-1} \,
\Phi^{\hbar_i}({\bf Q}^{\hbar_i})^{-1} \right) \restriction \mathscr{H}_\ell \\
& = \Phi^{\hbar_i}({\bf P}^{\hbar_i}\restriction \mathscr{H}_\ell) \,
\Phi^{\hbar_i}({\bf Q}^{\hbar_i} \restriction \mathscr{H}_\ell) \, 
\Phi^{\hbar_i}({\bf P}^{\hbar_i} \restriction \mathscr{H}_\ell)^{-1} \, 
\Phi^{\hbar_i}(({\bf P}^{\hbar_i} \restriction \mathscr{H}_\ell) + ({\bf Q}^{\hbar_i}\restriction \mathscr{H}_\ell))^{-1} \,
\Phi^{\hbar_i}({\bf Q}^{\hbar_i}\restriction \mathscr{H}_\ell)^{-1},
\end{align*}
and conjugating by $U'_\ell$ yields the left-hand-side of eq.\eqref{eq:cluster_pentagon_identity1_conjugated}, which equals $c'_{\bf m} \cdot {\rm Id}$ on $L^2(\mathbb{R},da_i)$; this means that, before conjugating by $U'_\ell$, it itself equals $(U'_\ell)^{-1} (c'_{\bf m}\cdot {\rm Id}) (U'_\ell) = c'_{\bf m} \cdot {\rm Id}$ on $\mathscr{H}_\ell$. This holds for each $\ell=1,\ldots,N$, which implies that the operator $\Phi^{\hbar_i}({\bf P}^{\hbar_i}) \, 
\Phi^{\hbar_i}({\bf Q}^{\hbar_i}) \,  
\Phi^{\hbar_i}({\bf P}^{\hbar_i})^{-1} \,
\Phi^{\hbar_i}({\bf P}^{\hbar_i} + {\bf Q}^{\hbar_i})^{-1} \,
\Phi^{\hbar_i}({\bf Q}^{\hbar_i})^{-1}$ equals $c'_{\bf m} \cdot {\rm Id}$ on $\mathscr{H}'$. Note that, a priori, $c'_{\bf m}$ is a constant that may depend on $\hbar,i,j$; by what we just got, it may depend on $\hbar_i$, but does not depend on other information about $\hbar,i,j$. The desired result follows when we replace $\hbar_i$ by a general nonzero real number, which is denoted by $\hbar$ in the statement. \qed

\begin{remark}
As pointed out in \cite{FG09}, this pentagon equation was suggested in \cite{F95} and proved in \cite{FKV01}, \cite{W}, and \cite{G08}. The above proof of the present paper relies heavily on the double version pentagon identity proved in \cite[Thm.5.5]{FG09}, and the key point was to extract a single version identity from the double version identity.
\end{remark}

\vs

Coming back to our situation, note from eq.\eqref{eq:Heisenberg_relations} that ${\bf x}^\hbar_i$ and ${\bf x}^\hbar_j$ are self-adjoint operators on the separable Hilbert space $\mathscr{H}_\Gamma$ and satisfy the Weyl relations for $[{\bf x}^\hbar_i, {\bf x}^\hbar_j] = 2\pi {\rm i} \hbar \, \wh{\varepsilon}_{ij} \cdot {\rm Id} = 2\pi {\rm i} \hbar_j \cdot {\rm Id}$. Thus from Prop.\ref{prop:pentagon_identity_for_ncQD} it follows that $c'_{\bf m}$ in eq.\eqref{eq:cluster_pentagon_identity1} equals $c^{\hbar_j}_{\rm pent} = c^{\hbar_i}_{\rm pent}$.

\vs

We now write down a slight variant of eq.\eqref{eq:general_pentagon_identity}. Namely, we move $\Phi^\hbar({\bf P}^\hbar) \, \Phi^\hbar({\bf Q}^\hbar)$ to the rightmost, by multiplying $\Phi^\hbar({\bf Q}^\hbar)^{-1} \, \Phi^\hbar({\bf P}^\hbar)^{-1}$ to both sides from left and multiplying $\Phi^\hbar({\bf P}^\hbar) \, \Phi^\hbar({\bf Q}^\hbar)$ to both sides from right:
\begin{align}
\nonumber
\Phi^\hbar({\bf P}^\hbar)^{-1} \, \Phi^\hbar({\bf P}^\hbar+{\bf Q}^\hbar)^{-1} \, \Phi^\hbar({\bf Q}^\hbar)^{-1} \, \Phi^\hbar({\bf P}^\hbar) \, \Phi^\hbar({\bf Q}^\hbar) = c_{\rm pent}^\hbar \, \cdot {\rm Id}.
\end{align}
Taking the inverse of both sides, we get
\begin{align}
\label{eq:general_pentagon_identity_variant2}
\Phi^\hbar({\bf Q}^\hbar)^{-1} \, 
\Phi^\hbar({\bf P}^\hbar)^{-1} \,
\Phi^\hbar({\bf Q}^\hbar) \,
\Phi^\hbar({\bf P}^\hbar+{\bf Q}^\hbar) \,
\Phi^\hbar({\bf P}^\hbar) = (c_{\rm pent}^\hbar)^{-1} \, \cdot {\rm Id}.
\end{align}
Note from eq.\eqref{eq:Heisenberg_relations_involving_tildes} that $\til{\bf x}^\hbar_j$ and $\til{\bf x}^\hbar_i$ are self-adjoint operators on the separable Hilbert space $\mathscr{H}_\Gamma$ and satisfy the Weyl relations for $[\til{\bf x}^\hbar_i, \til{\bf x}^\hbar_j] = - 2\pi {\rm i} \hbar \, \wh{\varepsilon}_{ij} \cdot {\rm Id} = - 2\pi {\rm i} \hbar_j \cdot {\rm Id} = - 2\pi {\rm i} \hbar_i \cdot {\rm Id}$. Thus from eq.\eqref{eq:general_pentagon_identity_variant2} it follows that $c''_{\bf m} = (c^{\hbar_i}_{\rm pent})^{-1}$. Hence $c_{A_2} = c_{\bf m} = c'_{\bf m} \, c''_{\bf m} = c^{\hbar_i}_{\rm pent} \, (c^{\hbar_i}_{\rm pent})^{-1} = 1$, as desired.

\vs

It remains to deal with the opposite case
$$
\varepsilon_{ij} = - \varepsilon_{ji} = -1
$$
separately. As computed in \S\ref{subsec:triviality_of_linear_parts}, the tropical sign-sequence is $\vec{\epsilon} = (+,+,+,-,-,0)$ this time. The constant $c_{A_2} \in {\rm U}(1)$ in eq.\eqref{eq:A2_to_prove} of Prop.\ref{prop:rank_2_identities_for_intertwiners}.(2-2) passes to eq.\eqref{eq:quantum_dilog_factors_only}; that is, $c_{\bf m} = c_{A_2}$. Equations \eqref{eq:separated_identity3}--\eqref{eq:separated_identity4} read
\begin{align*}
& \Phi^{\hbar_i}({\bf x}^\hbar_{{\bf m}_1}) \, \Phi^{\hbar_j}({\bf x}^\hbar_{{\bf m}_2}) \, \Phi^{\hbar_i}({\bf x}^\hbar_{{\bf m}_3}) \, 
\Phi^{\hbar_j}(-{\bf x}^\hbar_{{\bf m}_4})^{-1} \, 
\Phi^{\hbar_i}(-{\bf x}^\hbar_{{\bf m}_5})^{-1} = c'_{\bf m} \cdot {\rm Id},  \\
& \Phi^{\hbar_i}(\til{\bf x}^\hbar_{{\bf m}_1})^{-1} \, \Phi^{\hbar_j}(\til{\bf x}^\hbar_{{\bf m}_2})^{-1} \, \Phi^{\hbar_i}(\til{\bf x}^\hbar_{{\bf m}_3})^{-1} \, 
\Phi^{\hbar_j}(-\til{\bf x}^\hbar_{{\bf m}_4}) \, 
\Phi^{\hbar_i}(-\til{\bf x}^\hbar_{{\bf m}_5})  = c''_{\bf m} \cdot {\rm Id}.
\end{align*}
We now need to compute the argument operators as before, using from equations \eqref{eq:argument_operators_short}, \eqref{eq:bf_K_prime_conjugation_on_bf_x_i_new},  and \eqref{eq:bf_K_prime_conjugation_on_til_bf_x_i_new}. Proof is completely parallel with what we already did, so here we just record the result. One has $\varepsilon^{(\ell)}_{ij} = - \varepsilon^{(\ell)}_{ji} = -1$ for $\ell = 0,2$ and $\varepsilon^{(\ell)}_{ij} = - \varepsilon^{(\ell)}_{ji} = 1$ for $\ell = 1,3$, with $(k_1,k_2,k_3,k_4,k_5)=(i,j,i,j,i)$. One computes
$$
{\bf x}^\hbar_{{\bf m}_1} = {\bf x}^\hbar_i, \qquad
{\bf x}^\hbar_{{\bf m}_2} = {\bf x}^\hbar_i + {\bf x}^\hbar_j, \qquad
{\bf x}^\hbar_{{\bf m}_1} ={\bf x}^\hbar_j, \qquad
{\bf x}^\hbar_{{\bf m}_1} =-{\bf x}^\hbar_i, \qquad
{\bf x}^\hbar_{{\bf m}_1} =-{\bf x}^\hbar_j,
$$
and similarly for the tilde operators. We thus get
\begin{align*}
& \Phi^{\hbar_i}({\bf x}^\hbar_i) \, \Phi^{\hbar_i}({\bf x}^\hbar_i+{\bf x}^\hbar_j) \, \Phi^{\hbar_i}({\bf x}^\hbar_j) \, 
\Phi^{\hbar_i}({\bf x}^\hbar_i)^{-1} \, 
\Phi^{\hbar_i}({\bf x}^\hbar_j)^{-1} = c'_{\bf m} \cdot {\rm Id},  \\
& \Phi^{\hbar_i}(\til{\bf x}^\hbar_i)^{-1} \, \Phi^{\hbar_i}(\til{\bf x}^\hbar_i+\til{\bf x}^\hbar_j)^{-1} \, \Phi^{\hbar_i}(\til{\bf x}^\hbar_j)^{-1} \, 
\Phi^{\hbar_i}(\til{\bf x}^\hbar_i) \, 
\Phi^{\hbar_i}(\til{\bf x}^\hbar_j) = c''_{\bf m} \cdot {\rm Id}.
\end{align*}
From eq.\eqref{eq:general_pentagon_identity_variant2}, move $\Phi^\hbar({\bf Q}^\hbar)^{-1} \, \Phi^\hbar({\bf P}^\hbar)^{-1}$ to the rightmost, by multiplying $\Phi^\hbar({\bf P}^\hbar) \, \Phi^\hbar({\bf Q}^\hbar)$ from left and multiplying $\Phi^\hbar({\bf Q}^\hbar)^{-1} \, \Phi^\hbar({\bf P}^\hbar)^{-1}$ from right; since 
$[{\bf x}^\hbar_j, {\bf x}^\hbar_i] = 2\pi {\rm i} \hbar \, \wh{\varepsilon}_{ji} \cdot {\rm Id} = 2\pi {\rm i} \hbar_i \cdot {\rm Id}$, this fits the first of the above two equations, yielding $c'_{\bf m} = (c^{\hbar_i}_{\rm pent})^{-1}$. From eq.\eqref{eq:general_pentagon_identity}, move $\Phi^\hbar({\bf P})^\hbar) \, \Phi^\hbar({\bf Q}^\hbar)$ to the rightmost, in a similar manner; since $[\til{\bf x}^\hbar_i, \til{\bf x}^\hbar_j] = -2\pi {\rm i} \hbar \, \wh{\varepsilon}_{ij} \cdot {\rm Id} = 2\pi {\rm i} \hbar_j \cdot {\rm Id} = 2\pi {\rm i} \hbar_i \cdot {\rm Id}$, this fits the second of the above two equations, yielding $c''_{\bf m} = c^{\hbar_i}_{\rm pent}$.  Hence $c_{A_2} = c_{\bf m} = c'_{\bf m} c''_{\bf m} = 1$.

\subsection{The $B_2$ identity}

The remaining subsections have similar structure as \S\ref{subsec:A2}, so detailed computations will often be omitted, with only the results written. Consider the cluster transformation ${\bf m} = \mu_j \circ \mu_i \circ \mu_j \circ \mu_i \circ \mu_j \circ \mu_i$, applied on the $\mathcal{D}$-seed $\Gamma = \Gamma^{(0)}$, satisfying $\varepsilon_{ij} = - 2 \varepsilon_{ji} = \pm 2$ or $\varepsilon_{ji} = - 2 \varepsilon_{ij} = \pm 2$; in total, there are four cases to deal with separately. Let us start with the case
$$
\varepsilon_{ij} = 2, \quad \varepsilon_{ji} = -1.
$$
The tropical sign-sequence is $\vec{\epsilon} = (+,+,-,-,-,-)$, as computed in \S\ref{subsec:triviality_of_linear_parts}. One has $d_j = 2d_i$, so $\hbar_i = 2 \hbar_j$.

\vs

By the mutation formula in eq.\eqref{eq:varepsilon_prime_formula} for the exchange matrices, one has $\varepsilon^{(\ell)}_{ij} = 2$ and $\varepsilon^{(\ell)}_{ji} = -1$ for $\ell = 0,2,4$, while $\varepsilon^{(\ell)}_{ij} = - 2$ and  $\varepsilon^{(\ell)}_{ji} = 1$ for $\ell = 1,3$. Keeping in mind $(k_1,k_2,k_3,k_4,k_5,k_6) = (i,j,i,j,i,j)$, from equations \eqref{eq:argument_operators_short},   \eqref{eq:bf_K_prime_conjugation_on_bf_x_i_new}, and \eqref{eq:bf_K_prime_conjugation_on_til_bf_x_i_new}, and under the notation convention of \S\ref{subsec:argument_operators}, we can obtain the following by similar computations as before:
$$
{\bf x}^\hbar_{{\bf m}_1} = {\bf x}^\hbar_i, \qquad
{\bf x}^\hbar_{{\bf m}_2} = {\bf x}^\hbar_j, \qquad
{\bf x}^\hbar_{{\bf m}_3} = - {\bf x}^\hbar_i, \qquad
{\bf x}^\hbar_{{\bf m}_4} = - {\bf x}^\hbar_i - {\bf x}^\hbar_j, \qquad
{\bf x}^\hbar_{{\bf m}_5} = -{\bf x}^\hbar_i - 2{\bf x}^\hbar_j, \qquad
{\bf x}^\hbar_{{\bf m}_6} = - {\bf x}^\hbar_j,
$$
and similarly for the tilde operators. Thus from eq.\eqref{eq:separated_identity3}--\eqref{eq:separated_identity4} one gets
\begin{align}
\label{eq:separated_identity3_B2_case1}
& \Phi^{\hbar_i}({\bf x}^\hbar_i) \, \Phi^{\hbar_j}({\bf x}^\hbar_j) \, \Phi^{\hbar_i}({\bf x}^\hbar_i)^{-1} \, \Phi^{\hbar_j}({\bf x}^\hbar_i +{\bf x}^\hbar_j)^{-1} \, \Phi^{\hbar_i}({\bf x}^\hbar_i + 2{\bf x}^\hbar_j)^{-1} \, \Phi^{\hbar_j}({\bf x}^\hbar_j)^{-1} = c'_{\bf m} \cdot {\rm Id}, \\
\label{eq:separated_identity4_B2_case1}
& \Phi^{\hbar_i}(\til{\bf x}^\hbar_i)^{-1} \, \Phi^{\hbar_j}(\til{\bf x}^\hbar_j)^{-1} \, \Phi^{\hbar_i}(\til{\bf x}^\hbar_i) \, \Phi^{\hbar_j}(\til{\bf x}^\hbar_i +\til{\bf x}^\hbar_j) \, \Phi^{\hbar_i}(\til{\bf x}^\hbar_i + 2\til{\bf x}^\hbar_j) \, \Phi^{\hbar_j}(\til{\bf x}^\hbar_j) = c''_{\bf m} \cdot {\rm Id}
\end{align}
Note from eq.\eqref{eq:Heisenberg_relations} that ${\bf x}^\hbar_i$ and ${\bf x}^\hbar_j$ satisfy the Weyl relations for $[{\bf x}^\hbar_i, {\bf x}^\hbar_j] = 2\pi {\rm i} \hbar \, \wh{\varepsilon}_{ij} \cdot {\rm Id} = 4\pi {\rm i} \hbar_j \cdot \rm {\rm Id}$, and from eq.\eqref{eq:Heisenberg_relations_involving_tildes} that $[\til{\bf x}^\hbar_i, \til{\bf x}^\hbar_j]= - 4\pi {\rm i} \hbar_j \cdot {\rm Id}$. Similar arguments as in the proof of Prop.\ref{prop:pentagon_identity_for_ncQD} applied to eq.\eqref{eq:separated_identity3_B2_case1}--\eqref{eq:separated_identity4_B2_case1} yield:

\begin{proposition}[hexagon identity for the non-compact quantum dilogarithm]
\label{prop:hexagon_identity_for_ncQD}
There exist constants $c^\hbar_{\rm hex}, \til{c}^{\,\hbar}_{\rm hex} \in {\rm U}(1)$ depending on $\hbar$ s.t. the equalities of unitary operators
\begin{align}
\label{eq:general_hexagon_identity}
& \Phi^{2\hbar}(2{\bf P}^\hbar) \, \Phi^\hbar({\bf Q}^\hbar) \, \Phi^{2\hbar}(2{\bf P}^\hbar)^{-1} \, \Phi^\hbar(2{\bf P}^\hbar+{\bf Q}^\hbar)^{-1} \, \Phi^{2\hbar}(2{\bf P}^\hbar + 2{\bf Q}^\hbar)^{-1} \, \Phi^\hbar({\bf Q}^\hbar)^{-1} = c^\hbar_{\rm hex} \, \cdot {\rm Id}, \\
\label{eq:general_hexagon_identity_tilde}
& \Phi^{2\hbar}(2{\bf P}^\hbar)^{-1} \, \Phi^\hbar(-{\bf Q}^\hbar)^{-1} \, \Phi^{2\hbar}(2{\bf P}^\hbar) \, \Phi^\hbar(2{\bf P}^\hbar-{\bf Q}^\hbar) \, \Phi^{2\hbar}(2{\bf P}^\hbar - 2{\bf Q}^\hbar) \, \Phi^\hbar(-{\bf Q}^\hbar) = \til{c}^{\,\hbar}_{\rm hex} \, \cdot {\rm Id}
\end{align}
hold whenever ${\bf P}^\hbar, {\bf Q}^\hbar$ are self-adjoint operators on a separable Hilbert space that satisfy the Weyl relations for the Heisenberg relation $[{\bf P}^\hbar, {\bf Q}^\hbar] = 2\pi {\rm i} \hbar \cdot {\rm Id}$. \qed
\end{proposition}
The above form of the hexagon identity for the non-compact quantum dilogarithm function $\Phi^\hbar$ has been hinted in the literature \cite{KN} \cite[\S3]{Ip}, including the previous versions of the present paper; I just gave a rigorous proof. By a simple manipulation, we can observe the following.
\begin{lemma}
\label{lem:hexagon_constant_inverse}
One has $\til{c}^{\,\hbar}_{\rm hex} = (c^\hbar_{\rm hex})^{-1}$.
\end{lemma}

{\it Proof.} Assume the situation as in Prop.\ref{prop:hexagon_identity_for_ncQD}. Using the functional relation of the quantum dilogarithm in eq.\eqref{eq:QD_identity_quadratic}, we have $\Phi^\hbar(-{\bf Q}^\hbar) \, \Phi^\hbar({\bf Q}^\hbar) = c_\hbar \, e^{({\bf Q}^\hbar)^2/(4\pi{\rm i} \hbar)}$, where $c_\hbar = e^{-\frac{\pi{\rm i}}{12}(\hbar+\hbar^{-1})}$. For convenience, denote by $U$ the unitary operator $e^{({\bf Q}^\hbar)^2/(4\pi{\rm i} \hbar)}$, which makes sense through the functional calculus of ${\bf Q}^\hbar$. Then $U^{-1} {\bf Q}^\hbar U = {\bf Q}^\hbar$, and I claim that the equality 
$$
U^{-1} {\bf P}^\hbar U = {\bf P}^\hbar + {\bf Q}^\hbar
$$
of self-adjoint operators holds. Thanks to the Stone-von Neumann theorem (Prop.\ref{prop:SvN}), it suffices to prove this in case when the underlying Hilbert space is $L^2(\mathbb{R}, da)$ and ${\bf P}^\hbar = {\bf p}^\hbar = 2\pi {\rm i} \hbar \frac{d}{da}$ and ${\bf Q}^\hbar = {\bf q}^\hbar = a$. In this case, $U = e^{a^2/(4\pi{\rm i}\hbar)}$, which is just a multiplication operator. It is then an easy exercise to show that the above equality holds when applied to any vector $\psi$ in the nice dense subspace of $L^2(\mathbb{R},da)$ defined in eq.\eqref{eq:D_Gamma}.

\vs

Now, in the left-hand-side of eq.\eqref{eq:general_hexagon_identity_tilde}, replace the factor $\Phi^\hbar(-{\bf Q}^\hbar)^{-1}$ by $c_\hbar^{-1} \, \Phi^\hbar({\bf Q}^\hbar) \,  U^{-1}$, and replace $\Phi^\hbar(-{\bf Q}^\hbar)$ by $c_\hbar \, U \, \Phi^\hbar({\bf Q}^\hbar)^{-1}$. Then the scalars $c_\hbar^{-1}$ and $c_\hbar$ cancel each other, and we get
\begin{align*}
& \Phi^{2\hbar}(2{\bf P}^\hbar)^{-1} \, \Phi^\hbar({\bf Q}^\hbar) \, U^{-1} \, \Phi^{2\hbar}(2{\bf P}^\hbar) \, \Phi^\hbar(2{\bf P}^\hbar-{\bf Q}^\hbar) \, \Phi^{2\hbar}(2{\bf P}^\hbar - 2{\bf Q}^\hbar) \,  U\, \Phi^\hbar({\bf Q}^\hbar)^{-1} \\
& = \Phi^{2\hbar}(2{\bf P}^\hbar)^{-1} \, \Phi^\hbar({\bf Q}^\hbar) \, \, \Phi^{2\hbar}(U^{-1}  (2{\bf P}^\hbar) U) \, \Phi^\hbar(U^{-1} ( 2{\bf P}^\hbar -  {\bf Q}^\hbar) U) \, \Phi^{2\hbar}(U^{-1} (2 {\bf P}^\hbar - 2{\bf Q}^\hbar) U) \, \Phi^\hbar({\bf Q}^\hbar)^{-1} \\
& = \Phi^{2\hbar}(2{\bf P}^\hbar)^{-1} \, \Phi^\hbar({\bf Q}^\hbar) \, \, \Phi^{2\hbar}(2{\bf P}^\hbar + 2{\bf Q}^\hbar) \, \Phi^\hbar(2{\bf P}^\hbar + {\bf Q}^\hbar) \, \Phi^{2\hbar}(2 {\bf P}^\hbar) \, \Phi^\hbar({\bf Q}^\hbar)^{-1} \\
& = \Phi^{2\hbar}(2{\bf P}^\hbar)^{-1} \, \ul{ \Phi^\hbar({\bf Q}^\hbar) \, \, \Phi^{2\hbar}(2{\bf P}^\hbar + 2{\bf Q}^\hbar) \, \Phi^\hbar(2{\bf P}^\hbar + {\bf Q}^\hbar) \, \Phi^{2\hbar}(2 {\bf P}^\hbar) \, \Phi^\hbar({\bf Q}^\hbar)^{-1} \, \Phi^{2\hbar}(2{\bf P}^\hbar)^{-1} } \,\Phi^{2\hbar}(2{\bf P}^\hbar),
\end{align*}
and this equals $\Phi^{2\hbar}(2{\bf P}^\hbar)^{-1} \,((c^\hbar_{\rm hex})^{-1} \cdot {\rm Id}) \, \Phi^{2\hbar}(2{\bf P}^\hbar)= (c^\hbar_{\rm hex})^{-1} \cdot {\rm Id}$, for the underlined part is merely the inverse of the left-hand-side of eq.\eqref{eq:general_hexagon_identity}. This proves $\til{c}^{\, \hbar}_{\rm hex} = (c^\hbar_{\rm hex})^{-1}$. \qed

\vs

Applying the general hexagon identities back to eq.\eqref{eq:separated_identity3_B2_case1}--\eqref{eq:separated_identity4_B2_case1}, we get $c'_{\bf m} = c^{\hbar_j}_{\rm hex}$ and $c''_{\bf m} = \til{c}^{\, \hbar_j}_{\rm hex} = (c^{\hbar_j}_{\rm hex})^{-1}$. Thus $c_{B_2} = c_{\bf m} = c'_{\bf m} c''_{\bf m} =1$.

\vs

We turn to the Langlands dual case
$$
\varepsilon_{ij} = 1, \quad \varepsilon_{ji} = -2.
$$
The tropical sign-sequence is still $\vec{\epsilon} = (+,+,-,-,-,-)$, as computed in \S\ref{subsec:triviality_of_linear_parts}. One has $d_i = 2d_j$, so $\hbar_j = 2 \hbar_i$.

\vs

By the mutation formula in eq.\eqref{eq:varepsilon_prime_formula} for the exchange matrices, one has $\varepsilon^{(\ell)}_{ij} = 1$ and $\varepsilon^{(\ell)}_{ji} = -2$ for $\ell = 0,2,4$, while $\varepsilon^{(\ell)}_{ij} = - 1$ and  $\varepsilon^{(\ell)}_{ji} = 2$ for $\ell = 1,3$. Keeping in mind $(k_1,k_2,k_3,k_4,k_5,k_6) = (i,j,i,j,i,j)$, from equations \eqref{eq:argument_operators_short},   \eqref{eq:bf_K_prime_conjugation_on_bf_x_i_new}, and \eqref{eq:bf_K_prime_conjugation_on_til_bf_x_i_new}, and under the notation convention of \S\ref{subsec:argument_operators}, one computes
$$
{\bf x}^\hbar_{{\bf m}_1} = {\bf x}^\hbar_i, \qquad
{\bf x}^\hbar_{{\bf m}_2} ={\bf x}^\hbar_j, \qquad
{\bf x}^\hbar_{{\bf m}_3} =-{\bf x}^\hbar_i, \qquad
{\bf x}^\hbar_{{\bf m}_4} =-2{\bf x}^\hbar_i - {\bf x}^\hbar_j, \qquad
{\bf x}^\hbar_{{\bf m}_5} =-{\bf x}^\hbar_i - {\bf x}^\hbar_j, \qquad
{\bf x}^\hbar_{{\bf m}_6} =-{\bf x}^\hbar_j,
$$
and similarly for the tilde operators. Thus from eq.\eqref{eq:separated_identity3}--\eqref{eq:separated_identity4} one gets
\begin{align*}
& \Phi^{\hbar_i}({\bf x}^\hbar_i) \, \Phi^{\hbar_j}({\bf x}^\hbar_j) \, \Phi^{\hbar_i}({\bf x}^\hbar_i)^{-1} \, \Phi^{\hbar_j}(2{\bf x}^\hbar_i +{\bf x}^\hbar_j)^{-1} \, \Phi^{\hbar_i}({\bf x}^\hbar_i + {\bf x}^\hbar_j)^{-1} \, \Phi^{\hbar_j}({\bf x}^\hbar_j)^{-1} = c'_{\bf m} \cdot {\rm Id}, \\
& \Phi^{\hbar_i}(\til{\bf x}^\hbar_i)^{-1} \, \Phi^{\hbar_j}(\til{\bf x}^\hbar_j)^{-1} \, \Phi^{\hbar_i}(\til{\bf x}^\hbar_i) \, \Phi^{\hbar_j}(2\til{\bf x}^\hbar_i +\til{\bf x}^\hbar_j) \, \Phi^{\hbar_i}(\til{\bf x}^\hbar_i + \til{\bf x}^\hbar_j) \, \Phi^{\hbar_j}(\til{\bf x}^\hbar_j) = c''_{\bf m} \cdot {\rm Id}.
\end{align*}
Note $\hbar_j = 2\hbar_i$, and we have $[{\bf x}^\hbar_i, {\bf x}^\hbar_j] = 2\pi{\rm i} \hbar \wh{\varepsilon}_{ij} \cdot {\rm Id} = 2\pi {\rm i} \hbar_j \cdot {\rm Id} = 4\pi{\rm i} \hbar_i \cdot {\rm Id}$ and $[\til{\bf x}^\hbar_i, \til{\bf x}^\hbar_j] = - 4\pi{\rm i} \hbar_i \cdot {\rm Id}$. In eq.\eqref{eq:general_hexagon_identity_tilde}, one can move $\Phi^{2\hbar}(2{\bf P}^\hbar)^{-1} \, \Phi^\hbar(-{\bf Q}^\hbar)^{-1}$ to the rightmost, by multiplying $\Phi^\hbar(-{\bf Q}^\hbar) \, \Phi^{2\hbar}(2{\bf P}^\hbar)$ to both sides from left and multiplying $\Phi^{2\hbar}(2{\bf P}^\hbar)^{-1} \, \Phi^\hbar(-{\bf Q}^\hbar)^{-1}$ to both sides from right. Then, take the inverse of both sides; this fits the first of the above two equations, yielding $(\til{c}^{\, \hbar_i}_{\rm hex})^{-1} = c'_{\bf m}$. In eq.\eqref{eq:general_hexagon_identity}, one can move $\Phi^{2\hbar}(2{\bf P}^\hbar) \, \Phi^\hbar({\bf Q}^\hbar)$  to the rightmost, by multiplying $\Phi^\hbar({\bf Q}^\hbar)^{-1} \, \Phi^{2\hbar}(2{\bf P}^\hbar)^{-1}$ from left and multiplying $\Phi^{2\hbar}(2{\bf P}^\hbar) \, \Phi^\hbar({\bf Q}^\hbar)$ from right, and then take the inverse of both sides; this fits the second of the above two equations, yielding $(c^{\hbar_i}_{\rm hex})^{-1} = c''_{\bf m}$. Hence $c_{B_2} = c_{\bf m} = c'_{\bf m} c''_{\bf m} =1$. For later use, we record these variants of the hexagon identities:

\begin{align}
\label{eq:general_hexagon_identity_variant1}
& \Phi^{2\hbar}(2{\bf P}^\hbar) \, \Phi^\hbar(2{\bf P}^\hbar-{\bf Q}^\hbar) \, \Phi^{2\hbar}(2{\bf P}^\hbar - 2{\bf Q}^\hbar) \, \Phi^\hbar(-{\bf Q}^\hbar) \, \Phi^{2\hbar}(2{\bf P}^\hbar)^{-1} \, \Phi^\hbar(-{\bf Q}^\hbar)^{-1} = (c^\hbar_{\rm hex})^{-1} \, \cdot {\rm Id}, \\
\label{eq:general_hexagon_identity_variant2}
& \Phi^\hbar(-{\bf Q}^\hbar) \,
\Phi^{2\hbar}(2{\bf P}^\hbar) \,
\Phi^\hbar(-{\bf Q}^\hbar)^{-1} \,
\Phi^{2\hbar}(2{\bf P}^\hbar - 2{\bf Q}^\hbar)^{-1} \,
\Phi^\hbar(2{\bf P}^\hbar-{\bf Q}^\hbar)^{-1} \,
\Phi^{2\hbar}(2{\bf P}^\hbar) ^{-1} = c^\hbar_{\rm hex} \, \cdot {\rm Id}, \\
\label{eq:general_hexagon_identity_variant3}
& \Phi^{2\hbar}(2{\bf P}^\hbar)^{-1} \, \Phi^\hbar(2{\bf P}^\hbar+{\bf Q}^\hbar)^{-1} \, \Phi^{2\hbar}(2{\bf P}^\hbar + 2{\bf Q}^\hbar)^{-1} \, \Phi^\hbar({\bf Q}^\hbar)^{-1} \, \Phi^{2\hbar}(2{\bf P}^\hbar) \, \Phi^\hbar({\bf Q}^\hbar)  = c^\hbar_{\rm hex} \, \cdot {\rm Id}, \\
\label{eq:general_hexagon_identity_variant4}
& \Phi^\hbar({\bf Q}^\hbar)^{-1} \, 
\Phi^{2\hbar}(2{\bf P}^\hbar) ^{-1} \,  
\Phi^\hbar({\bf Q}^\hbar) \,
\Phi^{2\hbar}(2{\bf P}^\hbar + 2{\bf Q}^\hbar) \,
\Phi^\hbar(2{\bf P}^\hbar+{\bf Q}^\hbar) \,
\Phi^{2\hbar}(2{\bf P}^\hbar)   = (c^\hbar_{\rm hex})^{-1} \, \cdot {\rm Id}.
\end{align}

\vs

Consider the case
$$
\varepsilon_{ij} = -2, \quad \varepsilon_{ji} = 1.
$$
The tropical sign-sequence is $\vec{\epsilon} = (+,+,+,+,-,-)$, as computed in \S\ref{subsec:triviality_of_linear_parts}. One has $d_j = 2d_i$, so $\hbar_i = 2 \hbar_j$.

\vs

By the mutation formula in eq.\eqref{eq:varepsilon_prime_formula} for the exchange matrices, one has $\varepsilon^{(\ell)}_{ij} = -2$ and $\varepsilon^{(\ell)}_{ji} = 1$ for $\ell = 0,2,4$, while $\varepsilon^{(\ell)}_{ij} = 2$ and  $\varepsilon^{(\ell)}_{ji} = -1$ for $\ell = 1,3$. Keeping in mind $(k_1,k_2,k_3,k_4,k_5,k_6) = (i,j,i,j,i,j)$, from equations \eqref{eq:argument_operators_short},   \eqref{eq:bf_K_prime_conjugation_on_bf_x_i_new}, and \eqref{eq:bf_K_prime_conjugation_on_til_bf_x_i_new}, and under the notation convention of \S\ref{subsec:argument_operators}, one computes
$$
{\bf x}^\hbar_{{\bf m}_1} = {\bf x}^\hbar_i, \qquad
{\bf x}^\hbar_{{\bf m}_2} = {\bf x}^\hbar_i + {\bf x}^\hbar_j, \qquad
{\bf x}^\hbar_{{\bf m}_3} = {\bf x}^\hbar_i + 2{\bf x}^\hbar_j, \qquad
{\bf x}^\hbar_{{\bf m}_4} = {\bf x}^\hbar_j, \qquad
{\bf x}^\hbar_{{\bf m}_5} = - {\bf x}^\hbar_i, \qquad
{\bf x}^\hbar_{{\bf m}_6} = - {\bf x}^\hbar_j,
$$
and similarly for the tilde operators. Thus from eq.\eqref{eq:separated_identity3}--\eqref{eq:separated_identity4} one gets
\begin{align*}
& \Phi^{\hbar_i}({\bf x}^\hbar_i) \, \Phi^{\hbar_j}({\bf x}^\hbar_i + {\bf x}^\hbar_j) \, \Phi^{\hbar_i}({\bf x}^\hbar_i+2{\bf x}^\hbar_j) \, \Phi^{\hbar_j}({\bf x}^\hbar_j) \, \Phi^{\hbar_i}({\bf x}^\hbar_i)^{-1} \, \Phi^{\hbar_j}({\bf x}^\hbar_j)^{-1} = c'_{\bf m} \cdot {\rm Id}, \\
& \Phi^{\hbar_i}(\til{\bf x}^\hbar_i)^{-1} \, \Phi^{\hbar_j}(\til{\bf x}^\hbar_i + \til{\bf x}^\hbar_j)^{-1} \, \Phi^{\hbar_i}(\til{\bf x}^\hbar_i+2\til{\bf x}^\hbar_j)^{-1} \, \Phi^{\hbar_j}(\til{\bf x}^\hbar_j)^{-1} \, \Phi^{\hbar_i}(\til{\bf x}^\hbar_i) \, \Phi^{\hbar_j}(\til{\bf x}^\hbar_j) = c''_{\bf m} \cdot {\rm Id}.
\end{align*}
Note $\hbar_i = 2\hbar_j$, $[{\bf x}^\hbar_i,{\bf x}^\hbar_j]=2\pi {\rm i} \hbar \wh{\varepsilon}_{ij} \cdot {\rm Id} = -4\pi{\rm i} \hbar_j \cdot {\rm Id}$, and $[\til{\bf x}^\hbar_i,\til{\bf x}^\hbar_j]= 4\pi{\rm i} \hbar_j \cdot {\rm Id}$. So the first equation falls into eq.\eqref{eq:general_hexagon_identity_variant1}, while the second into eq.\eqref{eq:general_hexagon_identity_variant3}; hence $c'_{\bf m} = (c^\hbar_{\rm hex})^{-1}$ and $c''_{\bf m} = c^\hbar_{\rm hex}$, hence $c_{B_2} = c_{\bf m} = c'_{\bf m} c''_{\bf m} = 1$.

\vs

The Langlands dual to the previous case
$$
\varepsilon_{ij} = -1, \quad \varepsilon_{ji} = 2.
$$
The tropical sign-sequence is $\vec{\epsilon} = (+,+,+,+,-,-)$, as computed in \S\ref{subsec:triviality_of_linear_parts}. One has $d_i = 2d_j$, so $\hbar_j = 2 \hbar_i$.

\vs

By the mutation formula in eq.\eqref{eq:varepsilon_prime_formula} for the exchange matrices, one has $\varepsilon^{(\ell)}_{ij} = -1$ and $\varepsilon^{(\ell)}_{ji} = 2$ for $\ell = 0,2,4$, while $\varepsilon^{(\ell)}_{ij} = 1$ and  $\varepsilon^{(\ell)}_{ji} = -2$ for $\ell = 1,3$. Keeping in mind $(k_1,k_2,k_3,k_4,k_5,k_6) = (i,j,i,j,i,j)$, from equations \eqref{eq:argument_operators_short},   \eqref{eq:bf_K_prime_conjugation_on_bf_x_i_new}, and \eqref{eq:bf_K_prime_conjugation_on_til_bf_x_i_new}, and under the notation convention of \S\ref{subsec:argument_operators}, one computes
$$
{\bf x}^\hbar_{{\bf m}_1} ={\bf x}^\hbar_i, \qquad
{\bf x}^\hbar_{{\bf m}_2} = 2{\bf x}^\hbar_i + {\bf x}^\hbar_j, \qquad
{\bf x}^\hbar_{{\bf m}_3} = {\bf x}^\hbar_i + {\bf x}^\hbar_j, \qquad
{\bf x}^\hbar_{{\bf m}_4} = {\bf x}^\hbar_j, \qquad
{\bf x}^\hbar_{{\bf m}_5} = - {\bf x}^\hbar_i, \qquad
{\bf x}^\hbar_{{\bf m}_6} = - {\bf x}^\hbar_j,
$$
and similarly for the tilde operators. Thus from eq.\eqref{eq:separated_identity3}--\eqref{eq:separated_identity4} one gets
\begin{align*}
& \Phi^{\hbar_i}({\bf x}^\hbar_i) \, \Phi^{\hbar_j}(2{\bf x}^\hbar_i + {\bf x}^\hbar_j) \, \Phi^{\hbar_i}({\bf x}^\hbar_i+{\bf x}^\hbar_j) \, \Phi^{\hbar_j}({\bf x}^\hbar_j) \, \Phi^{\hbar_i}({\bf x}^\hbar_i)^{-1} \, \Phi^{\hbar_j}({\bf x}^\hbar_j)^{-1} = c'_{\bf m} \cdot {\rm Id}, \\
& \Phi^{\hbar_i}(\til{\bf x}^\hbar_i)^{-1} \, \Phi^{\hbar_j}(2\til{\bf x}^\hbar_i + \til{\bf x}^\hbar_j) ^{-1} \, \Phi^{\hbar_i}(\til{\bf x}^\hbar_i+\til{\bf x}^\hbar_j)^{-1} \, \Phi^{\hbar_j}(\til{\bf x}^\hbar_j)^{-1} \, \Phi^{\hbar_i}(\til{\bf x}^\hbar_i) \, \Phi^{\hbar_j}(\til{\bf x}^\hbar_j) = c''_{\bf m} \cdot {\rm Id}.
\end{align*}
Note $\hbar_j = 2\hbar_i$, $[{\bf x}^\hbar_j,{\bf x}^\hbar_i] = 2\pi {\rm i} \hbar \wh{\varepsilon}_{ji} \cdot {\rm Id} = 4\pi {\rm i} \hbar_i \cdot {\rm Id}$, and $[\til{\bf x}^\hbar_j, \til{\bf x}^\hbar_i] = - 4\pi {\rm i} \hbar_i \cdot {\rm Id}$. So one notes that the first equation corresponds to the inverse of eq.\eqref{eq:general_hexagon_identity}, yielding $c'_{\bf m} = (c^\hbar_{\rm hex})^{-1}$, while the second equation to the inverse of eq.\eqref{eq:general_hexagon_identity_tilde}, yielding $c''_{\bf m} = (\til{c}^{\,\hbar}_{\rm hex})^{-1} = c^\hbar_{\rm hex}$. Hence $c_{B_2} = c_{\bf m} = c'_{\bf m} c''_{\bf m} =1$.

\subsection{The $G_2$ identity}

Consider the cluster transformation ${\bf m} = \mu_j \circ \mu_i \circ \mu_j \circ \mu_i \circ \mu_j \circ \mu_i \circ \mu_j \circ \mu_i$, applied on the $\mathcal{D}$-seed $\Gamma = \Gamma^{(0)}$, satisfying $\varepsilon_{ij} = - 3 \varepsilon_{ji} = \pm 3$ or $\varepsilon_{ji} = - 3 \varepsilon_{ij} = \pm 3$; in total, there are four cases to deal with separately. Let us start with the case
$$
\varepsilon_{ij} = 3, \quad \varepsilon_{ji} = -1.
$$
The tropical sign-sequence is $\vec{\epsilon} = (+,+,-,-,-,-,-,-)$, as computed in \S\ref{subsec:triviality_of_linear_parts}. One has $d_j = 3d_i$, so $\hbar_i = 3 \hbar_j$.

\vs

By the mutation formula in eq.\eqref{eq:varepsilon_prime_formula} for the exchange matrices, one has $\varepsilon^{(\ell)}_{ij} = 3$ and $\varepsilon^{(\ell)}_{ji} = -1$ for $\ell = 0,2,4,6$, while $\varepsilon^{(\ell)}_{ij} = - 3$ and  $\varepsilon^{(\ell)}_{ji} = 1$ for $\ell = 1,3,5$. Keeping in mind $(k_1,k_2,k_3,k_4,k_5,k_6,k_7,k_8) = (i,j,i,j,i,j,i,j)$, from equations \eqref{eq:argument_operators_short},   \eqref{eq:bf_K_prime_conjugation_on_bf_x_i_new}, and \eqref{eq:bf_K_prime_conjugation_on_til_bf_x_i_new}, and under the notation convention of \S\ref{subsec:argument_operators}, one computes
\begin{align*}
\begin{array}{llll}
{\bf x}^\hbar_{{\bf m}_1} = {\bf x}^\hbar_i, 
& {\bf x}^\hbar_{{\bf m}_2} = {\bf x}^\hbar_j,
& {\bf x}^\hbar_{{\bf m}_3} = -{\bf x}^\hbar_i,
& {\bf x}^\hbar_{{\bf m}_4} = -{\bf x}^\hbar_i - {\bf x}^\hbar_j, \\
{\bf x}^\hbar_{{\bf m}_5} = -2{\bf x}^\hbar_i - 3{\bf x}^\hbar_j, 
& {\bf x}^\hbar_{{\bf m}_6} =-{\bf x}^\hbar_i - 2{\bf x}^\hbar_j, 
& {\bf x}^\hbar_{{\bf m}_7} =-{\bf x}^\hbar_i - 3{\bf x}^\hbar_j,
& {\bf x}^\hbar_{{\bf m}_8} =-{\bf x}^\hbar_j, 
\end{array}
\end{align*}
and similarly for the tilde operators. Thus from eq.\eqref{eq:separated_identity3}--\eqref{eq:separated_identity4} one gets
\begin{align*}
& \Phi^{\hbar_i}({\bf x}^\hbar_i) \, \Phi^{\hbar_j}({\bf x}^\hbar_j) \, \Phi^{\hbar_i}({\bf x}^\hbar_i)^{-1} \, \Phi^{\hbar_j}({\bf x}^\hbar_i +{\bf x}^\hbar_j)^{-1} \, \Phi^{\hbar_i}(2{\bf x}^\hbar_i + 3{\bf x}^\hbar_j)^{-1} \, \Phi^{\hbar_j}({\bf x}^\hbar_i+2{\bf x}^\hbar_j)^{-1} \, \Phi^{\hbar_i}({\bf x}^\hbar_i+3{\bf x}^\hbar_j)^{-1} \, \Phi^{\hbar_j}({\bf x}^\hbar_j)^{-1}  = c'_{\bf m} \cdot {\rm Id}, \\
& \Phi^{\hbar_i}(\til{\bf x}^\hbar_i)^{-1} \, \Phi^{\hbar_j}(\til{\bf x}^\hbar_j)^{-1} \, \Phi^{\hbar_i}(\til{\bf x}^\hbar_i) \, \Phi^{\hbar_j}(\til{\bf x}^\hbar_i +\til{\bf x}^\hbar_j) \, \Phi^{\hbar_i}(2\til{\bf x}^\hbar_i + 3\til{\bf x}^\hbar_j) \, \Phi^{\hbar_j}(\til{\bf x}^\hbar_i+2\til{\bf x}^\hbar_j) \, \Phi^{\hbar_i}(\til{\bf x}^\hbar_i+3\til{\bf x}^\hbar_j) \, \Phi^{\hbar_j}(\til{\bf x}^\hbar_j)  = c''_{\bf m} \cdot{\rm Id}.
\end{align*}
Note $\hbar_i = 3\hbar_j$, and ${\bf x}^\hbar_i$ and ${\bf x}^\hbar_j$ satisfy the Weyl relations for $[{\bf x}^\hbar_i, {\bf x}^\hbar_j] = 2\pi {\rm i} \hbar \, \wh{\varepsilon}_{ij} \cdot {\rm Id} = 6 \pi {\rm i} \hbar_j \cdot \rm {\rm Id}$, while $\til{\bf x}^\hbar_i$ and $\til{\bf x}^\hbar_j$ satisfy the Weyl relations for $[\til{\bf x}^\hbar_i, \til{\bf x}^\hbar_j] = - 6\pi{\rm i} \hbar_j \cdot {\rm Id}$. Similar arguments as in the proof of Prop.\ref{prop:pentagon_identity_for_ncQD} applied to these two equations yield:

\begin{proposition}[octagon identity for the non-compact quantum dilogarithm]
\label{prop:octagon_identity_for_ncQD}
There exists a constant $c^\hbar_{\rm oct} \in {\rm U}(1)$ depending on $\hbar$ s.t. the equalities of unitary operators
\begin{align}
\label{eq:general_octagon_identity}
& \Phi^{3\hbar}(3{\bf P}^\hbar) \,\Phi^\hbar({\bf Q}^\hbar) \, \Phi^{3\hbar}(3{\bf P}^\hbar)^{-1} \, \Phi^\hbar(3{\bf P}^\hbar+ {\bf Q}^\hbar)^{-1} \, \Phi^{3\hbar}(6{\bf P}^\hbar+3{\bf Q}^\hbar)^{-1} \, \Phi^\hbar (3{\bf P}^\hbar +2{\bf Q}^\hbar)^{-1} \Phi^{3\hbar}(3{\bf P}^\hbar + 3{\bf Q}^\hbar)^{-1} \, \Phi^\hbar({\bf Q}^\hbar)^{-1} \\
\nonumber
& = c^\hbar_{\rm oct} \, \cdot {\rm Id}, \\
\label{eq:general_octagon_identity_tilde}
& \Phi^{3\hbar}(3{\bf P}^\hbar)^{-1} \,\Phi^\hbar(-{\bf Q}^\hbar)^{-1} \, \Phi^{3\hbar}(3{\bf P}^\hbar) \, \Phi^\hbar(3{\bf P}^\hbar - {\bf Q}^\hbar) \, \Phi^{3\hbar}(6{\bf P}^\hbar-3{\bf Q}^\hbar) \, \Phi^\hbar (3{\bf P}^\hbar -2{\bf Q}^\hbar)\Phi^{3\hbar}(3{\bf P}^\hbar - 3{\bf Q}^\hbar) \, \Phi^\hbar(-{\bf Q}^\hbar) \\
\nonumber
& = (c^\hbar_{\rm oct})^{-1} \, \cdot {\rm Id}
\end{align}
hold whenever ${\bf P}^\hbar, {\bf Q}^\hbar$ are self-adjoint operators on a separable Hilbert space that satisfy the Weyl relations for the Heisenberg relation $[{\bf P}^\hbar, {\bf Q}^\hbar] = 2\pi {\rm i} \hbar \cdot {\rm Id}$.
\end{proposition}

A priori, the right-hand-side of eq.\eqref{eq:general_octagon_identity_tilde} should be written as e.g. $\til{c}^{\, \hbar}_{\rm oct} \cdot {\rm Id}$. By similar argument as in the proof of Lem.\ref{lem:hexagon_constant_inverse}, one can show that $\til{c}^{\, \hbar}_{\rm oct} = (c^\hbar_{\rm oct})^{-1}$. Namely, in eq.\eqref{eq:general_octagon_identity_tilde}, replace $\Phi^\hbar(-{\bf Q}^\hbar)^{-1}$ by $c_\hbar^{-1} \Phi^\hbar({\bf Q}^\hbar) U^{-1}$ and $\Phi^\hbar(-{\bf Q}^\hbar)$ by $c_\hbar U \Phi^\hbar({\bf Q}^\hbar)^{-1}$, where $U = e^{({\bf Q}^\hbar)^2/(4\pi {\rm i} \hbar)}$, so that $U^{-1} {\bf P}^\hbar U = {\bf P}^\hbar + {\bf Q}^\hbar$ and $U^{-1} {\bf Q}^\hbar U = {\bf Q}^\hbar$. Then eq.\eqref{eq:general_octagon_identity_tilde} becomes
$$
\Phi^{3\hbar}(3{\bf P}^\hbar)^{-1} \,\Phi^\hbar({\bf Q}^\hbar) \, \Phi^{3\hbar}(3{\bf P}^\hbar+ 3{\bf Q}^\hbar) \, \Phi^\hbar(3{\bf P}^\hbar +2 {\bf Q}^\hbar) \, \Phi^{3\hbar}(6{\bf P}^\hbar+3{\bf Q}^\hbar) \, \Phi^\hbar (3{\bf P}^\hbar +{\bf Q}^\hbar)\Phi^{3\hbar}(3{\bf P}^\hbar) \, \Phi^\hbar({\bf Q}^\hbar)^{-1} = \til{c}^{\, \hbar}_{\rm oct} \cdot {\rm Id}.
$$
Now move the leftmost factor $\Phi^{3\hbar}(3{\bf P}^\hbar)^{-1}$ to the rightmost, by multiplying $\Phi^{3\hbar}(3{\bf P}^\hbar)$ to both sides from left and multiplying $\Phi^{3\hbar}(3{\bf P}^\hbar)^{-1}$ to both sides from right. Then take the inverse of both sides; then the resulting left-hand-side equals the left-hand-side of eq.\eqref{eq:general_octagon_identity_tilde}, while the resulting right-hand-side is $(\til{c}^{\, \hbar}_{\rm oct})^{-1} \cdot {\rm Id}$, hence proving $c^\hbar_{\rm oct} = (\til{c}^{\, \hbar}_{\rm oct})^{-1}$ as desired.

\vs

The above form of the octagon identity for the non-compact quantum dilogarithm function $\Phi^\hbar$ has been hinted in the literature \cite{KN} \cite[\S3]{Ip}, including the previous versions of the present paper; I just gave a rigorous proof. Applying the general octagon identities to our specific identities, we get $c'_{\bf m} = c^{\hbar_j}_{\rm oct}$ and $c''_{\bf m} = (c^{\hbar_j}_{\rm oct})^{-1}$, hence $c_{G_2} = c_{\bf m} = c'_{\bf m} c''_{\bf m} =1$.

\vs

For later use, let us write down some variants of the octagon identities. One can take eq.\eqref{eq:general_octagon_identity}, move the leftmost part $\Phi^{3\hbar}(3{\bf P}^\hbar) \Phi^\hbar({\bf Q}^\hbar)$ to the rightmost, and then also can take the inverse of both sides. Or, one can take eq.\eqref{eq:general_octagon_identity_tilde}, move the leftmost part $\Phi^{3\hbar}(3{\bf P}^\hbar)^{-1} \Phi^\hbar(-{\bf Q}^\hbar)^{-1}$ to the rightmost, and then also can take the inverse of both sides. These are:
\begin{align}
\label{eq:general_octagon_identity_variant1}
& \Phi^{3\hbar}(3{\bf P}^\hbar)^{-1} \, \Phi^\hbar(3{\bf P}^\hbar+ {\bf Q}^\hbar)^{-1} \, \Phi^{3\hbar}(6{\bf P}^\hbar+3{\bf Q}^\hbar)^{-1} \, \Phi^\hbar (3{\bf P}^\hbar +2{\bf Q}^\hbar)^{-1} \Phi^{3\hbar}(3{\bf P}^\hbar + 3{\bf Q}^\hbar)^{-1} \, \Phi^\hbar({\bf Q}^\hbar)^{-1} \, \Phi^{3\hbar}(3{\bf P}^\hbar) \,\Phi^\hbar({\bf Q}^\hbar) \\
\nonumber
& = c^\hbar_{\rm oct} \, \cdot {\rm Id}, \\
\label{eq:general_octagon_identity_variant2}
& \Phi^\hbar({\bf Q}^\hbar)^{-1} \,
\Phi^{3\hbar}(3{\bf P}^\hbar)^{-1} \,
\Phi^\hbar({\bf Q}^\hbar) \,
\Phi^{3\hbar}(3{\bf P}^\hbar + 3{\bf Q}^\hbar) \,
\Phi^\hbar (3{\bf P}^\hbar +2{\bf Q}^\hbar) \,
\Phi^{3\hbar}(6{\bf P}^\hbar+3{\bf Q}^\hbar) \,
\Phi^\hbar(3{\bf P}^\hbar+ {\bf Q}^\hbar) \, \Phi^{3\hbar}(3{\bf P}^\hbar) \\
\nonumber
& = (c^\hbar_{\rm oct})^{-1} \, \cdot {\rm Id}, \\
\label{eq:general_octagon_identity_variant3}
& \Phi^{3\hbar}(3{\bf P}^\hbar) \, \Phi^\hbar(3{\bf P}^\hbar - {\bf Q}^\hbar) \, \Phi^{3\hbar}(6{\bf P}^\hbar-3{\bf Q}^\hbar) \, \Phi^\hbar (3{\bf P}^\hbar -2{\bf Q}^\hbar)\Phi^{3\hbar}(3{\bf P}^\hbar - 3{\bf Q}^\hbar) \, \Phi^\hbar(-{\bf Q}^\hbar) \, \Phi^{3\hbar}(3{\bf P}^\hbar)^{-1} \,\Phi^\hbar(-{\bf Q}^\hbar)^{-1} \\
\nonumber
& = (c^\hbar_{\rm oct})^{-1} \, \cdot {\rm Id}, \\
\label{eq:general_octagon_identity_variant4}
& \Phi^\hbar(-{\bf Q}^\hbar) \,
\Phi^{3\hbar}(3{\bf P}^\hbar) \,
\Phi^\hbar(-{\bf Q}^\hbar)^{-1} \,
\Phi^{3\hbar}(3{\bf P}^\hbar - 3{\bf Q}^\hbar)^{-1} \,
\Phi^\hbar (3{\bf P}^\hbar -2{\bf Q}^\hbar)^{-1} \,
\Phi^{3\hbar}(6{\bf P}^\hbar-3{\bf Q}^\hbar)^{-1} \,
\Phi^\hbar(3{\bf P}^\hbar - {\bf Q}^\hbar)^{-1} \,
\Phi^{3\hbar}(3{\bf P}^\hbar)^{-1} \\
\nonumber
& = c^\hbar_{\rm oct} \, \cdot {\rm Id}.
\end{align}

\vs

Then we turn to the Langlands dual case
$$
\varepsilon_{ij} = 1, \quad \varepsilon_{ji} = -3.
$$
The tropical sign-sequence is still $\vec{\epsilon} = (+,+,-,-,-,-,-,-)$, as computed in \S\ref{subsec:triviality_of_linear_parts}. One has $d_i = 3d_j$, so $\hbar_j = 3 \hbar_i$.

\vs

By the mutation formula in eq.\eqref{eq:varepsilon_prime_formula} for the exchange matrices, one has $\varepsilon^{(\ell)}_{ij} = 1$ and $\varepsilon^{(\ell)}_{ji} = -3$ for $\ell = 0,2,4,6$, while $\varepsilon^{(\ell)}_{ij} = - 1$ and  $\varepsilon^{(\ell)}_{ji} = 3$ for $\ell = 1,3,5$. Keeping in mind $(k_1,k_2,k_3,k_4,k_5,k_6,k_7,k_8) = (i,j,i,j,i,j,i,j)$, from equations \eqref{eq:argument_operators_short},   \eqref{eq:bf_K_prime_conjugation_on_bf_x_i_new}, and \eqref{eq:bf_K_prime_conjugation_on_til_bf_x_i_new}, and under the notation convention of \S\ref{subsec:argument_operators}, one computes
\begin{align*}
\begin{array}{llll}
{\bf x}^\hbar_{{\bf m}_1} = {\bf x}^\hbar_i, 
& {\bf x}^\hbar_{{\bf m}_2} = {\bf x}^\hbar_j, 
& {\bf x}^\hbar_{{\bf m}_3} = - {\bf x}^\hbar_i, 
& {\bf x}^\hbar_{{\bf m}_4} = - 3{\bf x}^\hbar_i - {\bf x}^\hbar_j, \\
{\bf x}^\hbar_{{\bf m}_5} = -2{\bf x}^\hbar_i - {\bf x}^\hbar_j,
& {\bf x}^\hbar_{{\bf m}_6} = - 3{\bf x}^\hbar_i - 2{\bf x}^\hbar_j, 
& {\bf x}^\hbar_{{\bf m}_7} = - {\bf x}^\hbar_i - {\bf x}^\hbar_j, 
& {\bf x}^\hbar_{{\bf m}_8} = - {\bf x}^\hbar_j, 
\end{array}
\end{align*}
and similarly for the tilde operators. Thus from eq.\eqref{eq:separated_identity3}--\eqref{eq:separated_identity4} one gets
\begin{align*}
& \Phi^{\hbar_i}({\bf x}^\hbar_i) \, 
\Phi^{\hbar_j}({\bf x}^\hbar_j) \, 
\Phi^{\hbar_i}({\bf x}^\hbar_i)^{-1} \, \Phi^{\hbar_j}(3{\bf x}^\hbar_i +{\bf x}^\hbar_j)^{-1} \, \Phi^{\hbar_i}(2{\bf x}^\hbar_i + {\bf x}^\hbar_j)^{-1} \, \Phi^{\hbar_j}(3{\bf x}^\hbar_i + 2 {\bf x}^\hbar_j)^{-1} \, \Phi^{\hbar_i}({\bf x}^\hbar_i + {\bf x}^\hbar_j)^{-1} \,
 \Phi^{\hbar_j}({\bf x}^\hbar_j)^{-1} = c'_{\bf m} \cdot {\rm Id}, \\
& \Phi^{\hbar_i}(\til{\bf x}^\hbar_i)^{-1} \, 
\Phi^{\hbar_j}({\bf x}^\hbar_j)^{-1} \, \Phi^{\hbar_i}(\til{\bf x}^\hbar_i) \, 
\Phi^{\hbar_j}(3\til{\bf x}^\hbar_i +\til{\bf x}^\hbar_j) \, \Phi^{\hbar_i}(2\til{\bf x}^\hbar_i + \til{\bf x}^\hbar_j) \, \Phi^{\hbar_j}(3\til{\bf x}^\hbar_i + 2 \til{\bf x}^\hbar_j) \, \Phi^{\hbar_i}(\til{\bf x}^\hbar_i + \til{\bf x}^\hbar_j) \, \Phi^{\hbar_j}(\til{\bf x}^\hbar_j) = c''_{\bf m} \cdot {\rm Id}.
\end{align*}
Note $\hbar_j = 3\hbar_i$, $[{\bf x}^\hbar_i, {\bf x}^\hbar_j] = 2\pi {\rm i} \hbar \wh{\varepsilon}_{ij} \cdot {\rm Id} = 2\pi {\rm i} \hbar_j \cdot {\rm Id} = 6 \pi {\rm i} \hbar_j\cdot {\rm Id}$, and $[\til{\bf x}^\hbar_i, \til{\bf x}^\hbar_j] = - 6\pi {\rm i} \hbar_i \cdot {\rm Id}$. Hence the first equation fits to eq.\eqref{eq:general_octagon_identity_variant4}, yielding $c'_{\bf m} = c^{\hbar_i}_{\rm oct}$, while the second equation fits to eq.\eqref{eq:general_octagon_identity_variant2}, yielding $c''_{\bf m} = (c^{\hbar_i}_{\rm oct})^{-1}$. So $c_{G_2} = c_{\bf m} = c'_{\bf m} c''_{\bf m} =1$.

\vs

Consider now the case
$$
\varepsilon_{ij} = -3, \quad \varepsilon_{ji} = 1.
$$
The tropical sign-sequence is $\vec{\epsilon} = (+,+,+,+,+,+,-,-)$, as computed in \S\ref{subsec:triviality_of_linear_parts}. One has $d_j = 3d_i$, so $\hbar_i = 3 \hbar_j$.

\vs

By the mutation formula in eq.\eqref{eq:varepsilon_prime_formula} for the exchange matrices, one has $\varepsilon^{(\ell)}_{ij} = -3$ and $\varepsilon^{(\ell)}_{ji} = 1$ for $\ell = 0,2,4,6$, while $\varepsilon^{(\ell)}_{ij} = 3$ and  $\varepsilon^{(\ell)}_{ji} = -1$ for $\ell = 1,3,5$. Keeping in mind $(k_1,k_2,k_3,k_4,k_5,k_6,k_7,k_8) = (i,j,i,j,i,j,i,j)$, from equations \eqref{eq:argument_operators_short},   \eqref{eq:bf_K_prime_conjugation_on_bf_x_i_new}, and \eqref{eq:bf_K_prime_conjugation_on_til_bf_x_i_new}, and under the notation convention of \S\ref{subsec:argument_operators}, one computes
\begin{align*}
\begin{array}{llll}
{\bf x}^\hbar_{{\bf m}_1} = {\bf x}^\hbar_i,
& {\bf x}^\hbar_{{\bf m}_2} = {\bf x}^\hbar_i + {\bf x}^\hbar_j, 
& {\bf x}^\hbar_{{\bf m}_3} = 2{\bf x}^\hbar_i + 3{\bf x}^\hbar_j, 
& {\bf x}^\hbar_{{\bf m}_4} = {\bf x}^\hbar_i + 2{\bf x}^\hbar_j, \\
{\bf x}^\hbar_{{\bf m}_5} = {\bf x}^\hbar_i + 3{\bf x}^\hbar_j, 
& {\bf x}^\hbar_{{\bf m}_6} = {\bf x}^\hbar_j, 
& {\bf x}^\hbar_{{\bf m}_7} = - {\bf x}^\hbar_i, 
& {\bf x}^\hbar_{{\bf m}_8} = - {\bf x}^\hbar_j, 
\end{array}
\end{align*}
and similarly for the tilde operators. Thus from eq.\eqref{eq:separated_identity3}--\eqref{eq:separated_identity4} one gets
\begin{align*}
& \Phi^{\hbar_i}({\bf x}^\hbar_i) \, \Phi^{\hbar_j}({\bf x}^\hbar_i + {\bf x}^\hbar_j) \, \Phi^{\hbar_i}(2{\bf x}^\hbar_i+3{\bf x}^\hbar_j) \, \Phi^{\hbar_j}({\bf x}^\hbar_i + 2{\bf x}^\hbar_j) \, \Phi^{\hbar_i}({\bf x}^\hbar_i+3{\bf x}^\hbar_j) \, \Phi^{\hbar_j}({\bf x}^\hbar_j) \, \Phi^{\hbar_i}({\bf x}^\hbar_i)^{-1} \, \Phi^{\hbar_j}({\bf x}^\hbar_j)^{-1} = c'_{\bf m} \cdot {\rm Id}, \\
& \Phi^{\hbar_i}(\til{\bf x}^\hbar_i)^{-1} \, \Phi^{\hbar_j}(\til{\bf x}^\hbar_i + \til{\bf x}^\hbar_j)^{-1} \, \Phi^{\hbar_i}(2\til{\bf x}^\hbar_i+3\til{\bf x}^\hbar_j)^{-1} \, \Phi^{\hbar_j}(\til{\bf x}^\hbar_i + 2\til{\bf x}^\hbar_j)^{-1} \, \Phi^{\hbar_i}(\til{\bf x}^\hbar_i+3\til{\bf x}^\hbar_j)^{-1} \, \Phi^{\hbar_j}(\til{\bf x}^\hbar_j)^{-1} \, \Phi^{\hbar_i}(\til{\bf x}^\hbar_i) \, \Phi^{\hbar_j}(\til{\bf x}^\hbar_j) = c''_{\bf m} \cdot {\rm Id}.
\end{align*}
Note $\hbar_i = 3\hbar_j$, $[{\bf x}^\hbar_i, {\bf x}^\hbar_j] = 2\pi {\rm i} \hbar \wh{\varepsilon}_{ij} \cdot {\rm Id} = -6\pi {\rm i} \hbar_j \cdot {\rm Id}$, and $[\til{\bf x}^\hbar_i, \til{\bf x}^\hbar_j] = 6\pi {\rm i} \hbar_j \cdot {\rm Id}$. Hence the first equation fits to eq.\eqref{eq:general_octagon_identity_variant3}, yielding $c'_{\bf m} = (c^{\hbar_i}_{\rm oct})^{-1}$, while the second equation fits to eq.\eqref{eq:general_octagon_identity_variant1}, yielding $c''_{\bf m} = c^{\hbar_i}_{\rm oct}$. So $c_{G_2} = c_{\bf m} = c'_{\bf m} c''_{\bf m} =1$.

\vs
The Langlands dual to the previous case
$$
\varepsilon_{ij} = -1, \quad \varepsilon_{ji} = 3.
$$
The tropical sign-sequence is $\vec{\epsilon} = (+,+,+,+,+,+,-,-)$, as computed in \S\ref{subsec:triviality_of_linear_parts}. One has $d_i = 3d_j$, so $\hbar_j = 3 \hbar_i$.

\vs

By the mutation formula in eq.\eqref{eq:varepsilon_prime_formula} for the exchange matrices, one has $\varepsilon^{(\ell)}_{ij} = -1$ and $\varepsilon^{(\ell)}_{ji} = 3$ for $\ell = 0,2,4,6$, while $\varepsilon^{(\ell)}_{ij} = 1$ and  $\varepsilon^{(\ell)}_{ji} = -3$ for $\ell = 1,3,5$. Keeping in mind $(k_1,k_2,k_3,k_4,k_5,k_6,k_7,k_8) = (i,j,i,j,i,j,i,j)$, from equations \eqref{eq:argument_operators_short},   \eqref{eq:bf_K_prime_conjugation_on_bf_x_i_new}, and \eqref{eq:bf_K_prime_conjugation_on_til_bf_x_i_new}, and under the notation convention of  \S\ref{subsec:argument_operators}, one computes
\begin{align*}
\begin{array}{llll}
{\bf x}^\hbar_{{\bf m}_1} = {\bf x}^\hbar_i,
& {\bf x}^\hbar_{{\bf m}_2} = 3{\bf x}^\hbar_i + {\bf x}^\hbar_j, 
& {\bf x}^\hbar_{{\bf m}_3} = 2{\bf x}^\hbar_i + {\bf x}^\hbar_j, 
& {\bf x}^\hbar_{{\bf m}_4} = 3{\bf x}^\hbar_i + 2{\bf x}^\hbar_j, \\
{\bf x}^\hbar_{{\bf m}_5} = {\bf x}^\hbar_i + {\bf x}^\hbar_j, 
& {\bf x}^\hbar_{{\bf m}_6} = {\bf x}^\hbar_j, 
& {\bf x}^\hbar_{{\bf m}_7} = -{\bf x}^\hbar_i, 
& {\bf x}^\hbar_{{\bf m}_8} = - {\bf x}^\hbar_j, 
\end{array}
\end{align*}
and similarly for the tilde operators. Thus from eq.\eqref{eq:separated_identity3}--\eqref{eq:separated_identity4} one gets
\begin{align*}
& \Phi^{\hbar_i}({\bf x}^\hbar_i) \, \Phi^{\hbar_j}(3{\bf x}^\hbar_i + {\bf x}^\hbar_j) \, \Phi^{\hbar_i}(2{\bf x}^\hbar_i+{\bf x}^\hbar_j) \, \Phi^{\hbar_j}(3{\bf x}^\hbar_i + 2{\bf x}^\hbar_j) \, \Phi^{\hbar_i}({\bf x}^\hbar_i+{\bf x}^\hbar_j) \, \Phi^{\hbar_j}({\bf x}^\hbar_j) \, \Phi^{\hbar_i}({\bf x}^\hbar_i)^{-1} \, \Phi^{\hbar_j}({\bf x}^\hbar_j)^{-1} = c'_{\bf m} \cdot {\rm Id}, \\
& \Phi^{\hbar_i}(\til{\bf x}^\hbar_i)^{-1} \, \Phi^{\hbar_j}(3\til{\bf x}^\hbar_i + \til{\bf x}^\hbar_j)^{-1} \, \Phi^{\hbar_i}(2\til{\bf x}^\hbar_i+\til{\bf x}^\hbar_j)^{-1} \, \Phi^{\hbar_j}(3\til{\bf x}^\hbar_i + 2\til{\bf x}^\hbar_j)^{-1} \, \Phi^{\hbar_i}(\til{\bf x}^\hbar_i+\til{\bf x}^\hbar_j)^{-1} \, \Phi^{\hbar_j}(\til{\bf x}^\hbar_j)^{-1} \, \Phi^{\hbar_i}(\til{\bf x}^\hbar_i) \, \Phi^{\hbar_j}(\til{\bf x}^\hbar_j) = c''_{\bf m} \cdot {\rm Id}.
\end{align*}
Note $\hbar_j = 3\hbar_i$, $[{\bf x}^\hbar_i, {\bf x}^\hbar_j] = 2\pi {\rm i} \hbar \wh{\varepsilon}_{ij} \cdot {\rm Id} = -2\pi {\rm i} \hbar_j \cdot {\rm Id} = - 6\pi {\rm i} \hbar_i \cdot {\rm Id}$, and $[\til{\bf x}^\hbar_i, \til{\bf x}^\hbar_j] = 6\pi {\rm i} \hbar_i \cdot {\rm Id}$. Hence the first equation fits to the inverse of eq.\eqref{eq:general_octagon_identity}, yielding $c'_{\bf m} = (c^{\hbar_i}_{\rm oct})^{-1}$, while the second equation fits to the inverse of eq.\eqref{eq:general_octagon_identity_tilde}, yielding $c''_{\bf m} = c^{\hbar_i}_{\rm oct}$. So $c_{G_2} = c_{\bf m} = c'_{\bf m} c''_{\bf m} =1$.

\section{Further research}

\vspace{-2mm}

The present paper dealt with some sequences ${\bf m}$ of mutations and seed automorphisms which start and end at a same seed. As mentioned before, classification of all such sequences is still an open problem, and the answer may or may not depend on whether it is for $\mathcal{A}$-, $\mathcal{X}$-, or $\mathcal{D}$-seeds. In \cite{KY}, some more examples of such sequences which are not consequences of the ones dealt with in the present paper are investigated; a shortest such sequence is of length $32$, applied to a rank $6$ seed coming from an ideal triangulation of a twice-punctured torus. I suggest trying to extend the results of the present paper to such sequences. The starting point of our arguments was Prop.\ref{thm:FG_constant}, which is \cite[Thm.5.4]{FG09}, saying that the composition of intertwiners for such a sequence is a scalar operator. Although Fock and Goncharov claimed the statement only for the sequences dealt with in the present paper, their proof in \cite{FG09} seems to apply to more general seed-trivial sequences; proof in \cite{FG09} uses the algebraic counterpart proved in \cite{BZ05} which does not restrict just to rank 2 identities. The next step would then be to show that Conjecture \ref{conj:dual_triviality}, which is about Langlands dual seeds and the tropical sign-sequence, holds for these more general sequences. Ideally, one should come up with a non-computational proof, but as a first step, one can check it e.g. for the above length 32 sequence by computation. One by-product would be rigorous proofs of new operator identities for $\Phi^\hbar$, in the style of Prop.\ref{prop:hexagon_identity_for_ncQD} and Prop.\ref{prop:octagon_identity_for_ncQD}. Lastly, one should show that  the arguments of \S\ref{sec:triviality_of_the_phase_constants} work for such sequences, i.e. that the constant coming from the non-tilde operator identity cancels that from the tilde operator identity. If one can come up with a more conceptual way of doing so, then it would also help shorten the present paper too.


\vspace{-3mm}


\begin{thebibliography}{GHKK17}





\bibitem[B01]{B01} E. W. Barnes, {\it Theory of the double gamma function}, Phil. Trans. Roy. Soc. A{\bf 196} (1901) 265--388.


\bibitem[BZ05]{BZ05} A. Berenstein and A. Zelevinsky, {\it Quantum cluster algebras}, Adv. Math. {\bf 195}(2) (2005) 405-455. \quad arXiv:math.QA/0404446

\bibitem[DWZ10]{DWZ} H. Derksen, J. Weyman, and A. Zelevinsky, {\it Quivers with potentials and their representations II: Applications to cluster algebras}, J. Amer. Math. Soc. {\bf 23} no.3 (2010), 749--790.


\bibitem[F95]{F95} L. D. Faddeev, {\it Discrete Heisenberg-Weyl group and modular group}, Lett. Math. Phys. {\bf 34} (1995) 249--254

\bibitem[FK94]{FK94} L. D. Faddeev and R. M. Kashaev, {\it Quantum dilogarithm}, Modern Phys. Lett. A{\bf 9} (1994) 427--434

\bibitem[FKV01]{FKV01} L. D. Faddeev, R. M. Kashaev, and A. Y. Volkov, {\it Strongly coupled quantum discrete Liouville theory, I: Algebraic approach and duality}, Commun. Math. Phys. {\bf 219} (2001) 199--219


\bibitem[FG07]{FG07} V. V. Fock and A. B. Goncharov, ``Cluster ensembles, quantization and the dilogarithm II: The intertwiner" in {\it Manin's Festschrift}, Birkh\"auser, 2007. \quad arXiv:math.QA/0702398 \quad MR2641183

\bibitem[FG09]{FG09} V. V. Fock and A. B. Goncharov, {\it The quantum dilogarithm and representations of the quantum cluster varieties}, Invent. Math. {\bf 175} (2009) 223--286.


\bibitem[FG09b]{FG09b} V. V. Fock and A. B. Goncharov, {\it Cluster ensembles, quantization and the dilogarithm}, Ann. Sci. \'Ec. Norm. Sup\'er. (4) {\bf 42}(6) (2009) 865--930. \quad arXiv:math.AG/0311245


\bibitem[FZ02]{FoZ02} S. Fomin and A. Zelevinsky, {\it Cluster Algebras I: Foundations}, J. Am. Math. Soc. {\bf 15}(2) (2002) 497--529.


\bibitem[FZ07]{FoZ07} S. Fomin and A. Zelevinsky, {\it Cluster Algebras IV: Coefficients}, Compos. Math. {\bf 143} (2007) 112-164.



\bibitem[FK12]{FK12} I. B. Frenkel and H. Kim, {\it Quantum Teichm\"uller space from the quantum plane}, Duke Math. J. {\bf 161} no.2 (2012) 305--366.

\bibitem[FS10]{FS10} L. Funar and V. Sergiescu, {\it Central extensions of the Ptolemy-Thompson group and quantized Teichm\"uller theory}, J. Topol. {\bf 3} (2010) 29--62.

\bibitem[G08]{G08} A. B. Goncharov, ``Pentagon relation for the quantum dilogarithm and quantized $\mcal{M}^{\rm cyc}_{0,5}$'' in {\it Geometry and Dynamics of Groups and Spaces}, Progr. Math. {\bf 265}, Birkh\"auser, Basel, 2008, pp.415--428.

\bibitem[GHKK17]{GHKK} M. Gross, P. Hacking, S. Keel, and M. Kontsevich, {\it Canonical bases for cluster algebras}, J. Amer. Math. Soc. {\bf 31} no.2 (2018), 497--608. \quad arXiv:1411.1394

\bibitem[H13]{Hall} B. C. Hall, {\it Quantum theory for mathematicians}, GTM {\bf 267}, Springer, New York, 2013. 

\bibitem[I15]{Ip} I. C. Ip, {\it Positive representations of split real quantum groups: the universal R operator}, Int. Math. Res. Not. IMRN 2015 no.1, 240--287. \quad arXiv:1212.5149



\bibitem[KN11]{KN} R. M. Kashaev and T. Nakanishi, {\it Classical and quantum dilogarithm identities}, SIGMA Symmetry Integrability Geom. Methods Appl. {\bf 7} (2011), Paper 102, 29 pages. \quad arXiv:1104.4630





\bibitem[K11]{Keller} Keller, ``On cluster theory and quantum dilogarithm identities'' in {\it Representations of algebras and related topics}, EMS Ser. Congr. Rep., Eur. Math. Soc., Zürich, 2011, pp.85–-116. \quad arXiv:1102.4148



\bibitem[K12]{K12} H. Kim, {\it The dilogarithmic central extension of the Ptolemy-Thompson group via the Kashaev quantization}, Adv. Math. {\bf 293} (2016), 529--588. \quad arXiv:1211.4300


\bibitem[K16l]{K16} H. Kim, {\it Phase constants in the Fock-Goncharov quantization of cluster varieties: long version}, arXiv:1602.00361v2

\bibitem[KS]{KS} H. Kim and C. Scarinci, {\it Quantization of moduli spaces of three-dimensional spacetimes} (tentative), in preparation.

\bibitem[KY]{KY} H. Kim and M. Yamazaki, {\it Comments on exchange graphs in cluster algebras}, to appear in Exp. Math., arXiv:1612.00145


\bibitem[RS80]{RS70} Michael Reed and Barry Simon, {\it Methods of Modern Mathematical Physics I: Functional Analysis.}, Revised and enlarged ed., Academic Press, San Diego, 1980.

\bibitem[S12]{Schmudgen} K. Schm\"udgen, {\it Unbounded self-adjoint operators on Hilbert space}, GTM {\bf 265}, Springer, Dordrecht, 2012.

\bibitem[v31]{vN31} J. von Neumann, {\it Die Eindeutigkeit der Schr\"odingerschen Operatoren}, Mathematische Annalen {\bf 104} (1931) 570--578. 


\bibitem[W00]{W} S. L. Woronowicz, {\it Quantum exponential function}, Rev. Math. Phys. {\bf 12} no.6 (2000) 873--920.

\bibitem[X14]{X14} B. Xu, {\it Central extension of mapping class group via Chekhov-Fock quantization}, J. Geom. Phys. {\bf 110} (2016), 9--24.



\end{thebibliography}
\end{document}